\documentclass[12pt]{article}
\usepackage{amsfonts,amsmath,amssymb,latexsym}
\usepackage{eucal}
\usepackage{color}
\usepackage{geometry}

\newtheorem{theorem}{Theorem}[section]
\newtheorem{proposition}[theorem]{Proposition}
\newtheorem{corollary}[theorem]{Corollary}
\newtheorem{lemma}[theorem]{Lemma}
\newtheorem{definition}{Definition}[section]
\newtheorem{example}{Example}[section]
\newtheorem{remark}{Remark}[section]

\numberwithin{equation}{section}

\begin{document}

\title{Mean curvature flow solitons in the presence of conformal vector fields}
\author{ L. J. Al\'\i as\thanks{This research is a result of the activity developed within the framework of the Programme in Support of Excellence Groups of the Regi\'on de Murcia, Spain, by Fundaci\'on S\'eneca, Science and Technology Agency of the Regi\'on de Murcia.
Partially supported by MINECO/FEDER project reference MTM2015-65430-P and Fundaci\'on S\'eneca project reference 19901/GERM/15, Spain.},  J. H. de Lira\thanks
{Partially supported by CNPq Produtividade em Pesquisa Grant $\#$  302067/2014-0 and FUNCAP/CNPq/PRONEX Grant ``N\'ucleo de An\'alise Geom\'etrica e Aplica\c c\~oes'' $\#$ 09.01.00/11.}, M. Rigoli}
\date{}
\maketitle

\tableofcontents

\section{Introduction}


One of the main research topics on mean curvature flow is the study of self-similar solutions particularly as models of singularities after suitable rescalings of the flow  \cite{huisken}, \cite{huisken-2}, \cite{white}, \cite{colding}, \cite{colding2}, \cite{colding3}. From another point of view, self-shrinkers, self-expanders and translating solitons can be considered as weighted minimal submanifolds for some suitably chosen densities in Euclidean space, see for instance \cite{morgan}, \cite{ilmanen}, \cite{Lott}, \cite{zhou} and references therein. In other terms, they are minimal submanifolds with respect to a Riemannian metric conformal to the Euclidean one. Finally, as pointed out in \cite{colding}  self-shrinkers are  critical hypersurfaces for the entropy functional. 




The subject experienced an increasing activity after the seminal paper by Colding and Minicozzi \cite{colding} that inspired an impressive amount of work on existence and classification problems, rigidity results, stability and spectral properties \cite{colding2}, \cite{colding3}, \cite{white2}, \cite{sesum}, \cite{CL}, \cite{cheng-peng}, \cite{DX}, \cite{martin2} only to mention a few contributions most related to the present work.  

However,  to the best of our knowledge, there is no systematic investigation of mean curvature flow solitons in Riemannian ambient spaces in spite of some relevant contributions as, for instance, \cite{hunger}, \cite{Futaki}, \cite{ortega}, \cite{yamamoto}.  The aim of this paper is to introduce a notion of mean curvature flow soliton general enough to encompass target spaces of constant sectional curvature, Riemannian products or, in increasing generality, warped product spaces. In fact, the definition we propose only refers to some vector field $X$ on the target manifold. A number of examples described in Sections \ref{ssf} and \ref{fundamental} indicates how this general notion  recovers those already existing in the literature. As expected the definition below is motivated by the self-similarity of certain special solutions of the mean curvature flow with respect to the flow generated by the vector field. 


\begin{definition}\label{main-soliton-definition}
An isometric immersion $\psi:M^m\to \bar M^{n+1}$ is a \emph{mean
curvature flow soliton}  with respect to $X\in \Gamma(T\bar M)$  if
\begin{equation}
\label{solitonA-2-intro}c\, X^\perp = {\bf H}
\end{equation}
along $\psi$ for some constant $c\in \mathbb{R}$. 
\end{definition}

Euclidean self-shrinkers and self-expanders correspond  to the choice of $X$ as the position vector field  in $\bar M = \mathbb{R}^{n+1}$ and, respectively, constants $c<0$ and $c>0$. In the same way, translating solitons in Riemannian products $\bar M=\mathbb{R}\times P$ are mean curvature flow solitons with respect to the parallel vector field that generates translations along the factor $\mathbb{R}$. 
In all these cases, $X$ is a closed conformal vector field. Restricting ourselves to this class of vector fields we obtain  as a consequence of the definition above the tensor equation

\begin{equation}
\label{h-ric}
II_{-{\bf H}} + \frac{c}{2}\pounds_{X^\top} g = c\varphi g
\end{equation}
where $g$ is the Riemannian metric in the soliton and $\varphi = \frac{1}{n+1}{\rm div} X$. This equation is structurally close to that defining Ricci solitons and it enables us to recover, for instance, a Hamilton type identity and related consequences in a quite general context. Because of these similarities and other further properties we expect that the machinery of maximum principles, recently applied by one of us   to the study of complete (non-compact) Ricci solitons and to Euclidean self-shrinkers \cite{catino}, \cite{BPR} should provide an effective tool of investigation also in this context. This is indeed the case  under the assumption that $X$ is a conformal closed vector field. In this case we identify some natural geometric quantities that satisfy  elliptic equations or differential inequalities in a simple and manageable form  for which the weak maximum principle is valid as explained in detail in Section \ref{wmp}.
As it is well known the existence of a closed conformal vector field imposes  (in general local) restrictions on the geometry of $\bar M$ analyzed in the work of Montiel \cite{montiel}. It turns out that is not a severe loss of generality to  restrict ourselves to manifolds $\bar M$ given by warped product spaces. 
However, when possible we state our results in greater generality.

In sum, in some sense our most important contribution is to fix a unified, geometrically natural and analitically treatable notion of self-similar solutions of the mean curvature flow in the presence of a conformal or parallel vector field $X$ in in the ambient space. We now describe the general plan of the paper; in doing so we also indicate  some results that we feel to be of a certain interest.

\

\noindent{\bf Plan of the paper.}
The paper is organized as follows. In Section \ref{ssf} we introduce the general definition of self-similar mean curvature flow in a Riemannian manifold $\bar M^{n+1}$ endowed with a vector field $X$. We illustrate the definition with a number of examples showing also, for instance in the Euclidean space with $X$ the position vector field, that the new definition agrees with those known in the literature. Further examples enlightening the dependence on $X$ are given in hyperbolic space, the sphere, Riemannian and warped product spaces and so on. In Proposition \ref{soliton-prop} by way of equations we characterize the notion of self-similar mean curvature flow. For $X$ conformal the result points the formal similarity with Ricci solitons mentioned above and that will be reconsidered in Section \ref{fundamental} and analyzed in some aspects in Section \ref{wmp}. 

We formally introduce the notion of mean curvature soliton in Definition \ref{soliton-definition} in Section \ref{fundamental}. Note that in our definition the mean curvature vector field is non-normalized. Starting from equation (\ref{solitonA-2-intro}) we basically reprove Proposition \ref{soliton-prop} but from a different perspective. We then give further examples in several interesting geometric settings $\bar M$.

Section \ref{warped} is devoted to the detailed description of the special case of warped product spaces of the form $I\times_hP^n$ with $h:I\to \mathbb{R}^+=(0,+\infty)$, $I\subset \mathbb{R}$ an open interval. From this section on the vector field $X$ on $\bar M$ will always be assumed closed and conformal. Next, in Section \ref{sec5} we introduce a function which will be essential in our study, namely
\[
\eta(x) = \int^{t(\psi(x))}_{t_0} h(\tau)\, {\rm d}\tau,\quad x\in M,
\]
for some fixed $t_0\in I$, where $t$ denotes the natural coordinate in the factor $I$. With the aid of the previous results, in Proposition \ref{laplace-eta} we establish the first two basic equations
\begin{align}
& \Delta \eta = m\varphi + \frac{1}{c}|{\bf H}|^2 \label{delta-intro-1}\\
& \Delta_{-c\eta}\eta = m\varphi + c|X|^2 \label{delta-intro-2}
\end{align}
 that we shall use for the study of the geometry of mean curvature flow solitons. The choices $\bar M=I\times_hP^n$ and $X=h(t)\partial_t$ reveals necessary to give to equations (\ref{delta-intro-1}) and (\ref{delta-intro-2}) a form that will be possible to treat with the analytical tools in part introduced in the Section \ref{wmp}. Indeed, the aim of Section \ref{wmp} is to introduce and describe the weak maximum principle basically for operators of the form $\Delta_f=\Delta-\langle\nabla f,\nabla\,\cdot\rangle$ for some $f\in C^1(M)$. We recall some sufficient conditions for its validity, as well as its equivalent \textit{open} version. We  extend our considerations also to parabolicity viewed as a stronger form of the weak maximum principle. We then show some simple and geometrically significant conditions to guarantee that the weak maximum principle holds for a certain operator on a complete mean curvature flow soliton. This is the content of Theorem \ref{cond-wmp}. The proof of this result is based on an upper weighted volume estimate of geodesic balls. Other estimates of this type can also be obtained via the similarity of mean curvature flow solitons with Ricci solitons established by (\ref{h-ric}). At this respect we refer the reader to Propositions \ref{ode-eta} to \ref{vol-quad}.

Section \ref{sec7} contains the first applications of the weak maximum principle to mean curvature flow solitons. We begin this section by characterizing those immersions with image in a leaf of $\bar M=I\times_hP^n$ that are mean curvature flow solitons with respect to the vector field $X=h(t)\partial_t$. This will justify the study of the function 
\begin{equation}
\zeta(t)=mh'(t)+ch^2(t)
\end{equation}
 that will frequently appear in our considerations and that will be called the \textit{soliton function}.

We then give some first simple applications of the weak maximum principle and parabolicity for the operator $\Delta_{-c\eta}$. Our results in this section are all based on the analysis of the equations (\ref{delta-intro-1}) and (\ref{delta-intro-2}) in the particular case of warped product spaces. The aimed result is to show that, under various different assumptions, the image of a mean curvature flow soliton is a leaf of the natural foliation of $\bar M=I\times_hP^n$ induced by $X$. We also give some height estimates, for example those contained in Corollary \ref{hyp-sphere}.

Sections \ref{char} and \ref{simons} constitute the core of the paper. In Proposition  \ref{laplacian-H2} we introduce a new fundamental equation satisfied by $|{\bf H}|^2$ in complete generality, that is, for a mean curvature flow soliton with respect to a generic closed conformal vector field $X$ on the target space. However, in order to be able to deal with the complexity of equation (\ref{deltaH2}), we restrict ourselves to the case $\bar M = I\times_h P$, a warped product space with fiber $P$ of constant sectional curvature. In this setting, we deduce some rigidity results when $|{\bf H}|^2$ is integrable (in a weighted measure), see Theorem  \ref{theorem30}.  Moreover, Theorem \ref{theorem30.1} extends to our much more general setting results proved by  Cao and Li \cite{CL} and Ding and Xin \cite{DX} for self-shrinkers in the Euclidean space. See also Remark \ref{remark8.1} and what follows for a characterization result in the same spirit.

In Subsection \ref{subsectionEinstein} we deal with codimension one mean curvature flow solitons in the particular case when  $\bar M$ is an Einstein space. We then obtain some vanishing results for $|{\bf H}|$ both with the help of the weak maximum principle or with a second technique developed in \cite{PRSgreenbook}. Theorem \ref{gap-E} generalizes a well known result for self-shrinkers of $\mathbb{R}^{n+1}$ due to Cao and Li \cite{CL}. Corollary \ref{cor-gap-E1} covers the codimension $1$ case in $\mathbb{R}^{n+1}$ extending \cite{CL} and it will be given also in a general codimension case but with more stringent assumptions in Corollary \ref{gap-Fbis}.

Section 9 deals with applications of a Simons' type formula for mean curvature flow solitons in a space of constant sectional curvature. Again we produce a sort of rigidity results such as Theorem \ref{huisken-cl} for the compact case and Theorem \ref{huisken-complete} for the non-compact case for which we prove the vanishing of the tensor
\[
V = \nabla H \otimes A - H \nabla A
\]
with $A$ the Weingarten operator. In favorable circumstances, for instance for $\bar M = \mathbb{R}^{n+1}$, the vanishing of $V$ gives rise to a classification result in terms of special examples of mean curvature flow solitons via the work of Mart\'\i n and collaborators \cite{martin2}. It is worth to point out that for general ambient spaces much has to be done in this direction. Particularly interesting is Theorem \ref{gap-G} together with the observation that assumption (\ref{growth-Hp}) can be substituted with (\ref{growth-Hp.ex})  implied by an $L^p$ condition on $H$ and $A$ as in assumption (\ref{H2p-bis}) of Corollary \ref{gap-self}. The section ends with three results proving that the soliton is totally umbilical under different various assumptions. 

Section 10 shows the non-existence of translating soliton graphs confined in appropriate regions of a product space $\mathbb{R}\times P$. See Theorem \ref{th-transl-2}. Some examples and a further elaborated discussion on translating solitons in general Riemannian products may be found at \cite{ML}.

In the next section we present, in the codimension one case for simplicity, how the equation characterizing mean curvature flow solitons can be recovered as the Euler-Lagrange equation of an appropriate weighted volume functional. In this way, we justify the notion of stable as well as both  finite and infinite index mean curvature flow solitons. We then study relations between the volume growth of geodesic balls and other geometric quantities via the eigenvalues of the weighted Laplacian naturally associated to a soliton immersed into a warped product $I\times_h P$; see for instance Corollary \ref{cor-G2}. We end this paper with some sufficient conditions for the soliton to have infinite index. We emphasize the fact that the techniques used in this section are quite different from the applications of the maximum principle in the previous sections.


\section{Self-similar mean curvature flows}
\label{ssf}

Let $M^m$ and $\bar M^{n+1}$ be Riemannian
manifolds. Given $\omega_* < 0 <\omega^*$, we consider a
differentiable map
\begin{equation}
\Psi:(\omega_*,\omega^*)\times M\to \bar M
\end{equation}
such that $\Psi_{\tau} = \Psi(\tau,\cdot)$ is an immersion, for all $\tau\in (\omega_*,\omega^*)$. We denote
 $\psi=\Psi_0$. The submanifolds $\Psi_{\tau}(M),\,\tau\in
(\omega_*,\omega^*),$ are evolving by their mean curvature vector
field if
\begin{equation}
\frac{d\Psi}{d\tau}=\Psi_* \frac{\partial}{\partial \tau} = {\bf
H},
\end{equation}
 where
\begin{equation}
{\bf H} =\bigg(\sum_{i=1}^{m} 
\bar\nabla_{\Psi_{\tau*}{\sf e}_i} \Psi_{\tau*} {\sf e}_i\bigg)^\perp
\end{equation}
is the (non-normalized) mean curvature vector of $\Psi_\tau$. Here and in what follows $\perp$ indicates the projection onto the normal bundle; the local  tangent frame $\{{\sf
e}_i\}_{i=1}^m$ is orthonormal with respect to the metric induced in $M$ by $\Psi_\tau$. The notations $\bar{g}=\langle\cdot, \cdot\rangle$ and $\bar\nabla$ stand for the Riemannian metric and connection in $\bar M$, respectively.

For a given vector field $X$ on $\bar M$, we set $\Phi:(\Omega_*, \Omega^*)\times \bar M \to \bar M$ to denote the
flow generated by $X$ defined in the maximal interval $(\Omega_*,
\Omega^*)$. Let  $s$ be the flow parameter in $\Phi$ and define
\begin{equation}
\widetilde \Psi_{\tau} (x) =\widetilde \Psi(\tau, x)=
\Phi^{-1}(\sigma(\tau),\Psi_{\tau}(x)), \quad x\in M,
\end{equation}
where $\sigma:(\omega_*,\omega^*)\to (\Omega_*, \Omega^*)$ is a
reparametrization of the flow lines of $X$ of the form
\[
s=\sigma(\tau).
\]
Equivalently we can write
\begin{equation}
\label{psiphi} \Psi(\tau,x) =
\Phi(\sigma(\tau),\widetilde\Psi(\tau,x)), \quad (\tau,x)\in (\omega_*,\omega^*)\times M .
\end{equation}

\begin{definition}
\label{selfsimilar} Let $\bar M^{n+1}$ be a Riemannian manifold
endowed with a  vector field $X\in \Gamma(T\bar
M)$. Given an $m$-dimensional Riemannian manifold $M^m$ we say that
a mean curvature flow $\Psi:(\omega_*, \omega^*)\times M \to \bar
M$ is \emph{self-similar} if  there exists an isometric immersion
$\psi:M \to \bar M$ and a reparametrization $\sigma:(\omega_*, \omega^*)\to (\Omega_*, \Omega^*)$ of the flow lines of $X$ such that
\begin{equation}
\Psi_{\tau}(M) =\Phi_{\sigma(\tau)}(\psi(M)),
\end{equation}
for all $\tau\in (\omega_*,\omega^*)$, where
 $\Phi:(\Omega_*,
\Omega^*)\times \bar M\to\bar M$ is the flow generated by $X$. In other terms, $\widetilde\Psi_\tau(M) = \psi(M)$, for all $\tau\in (\omega_*, \omega^*)$.
\end{definition}

Although the above definition does not require special properties of $X$, in the examples below, used to illustrate the concepts, we deal with a conformal vector field $X$ and later on we will restrict ourselves to a closed conformal vector field. This choice is due to the desire of simplifying equations, as we will see. For the moment, let us recall that $X$ is said to be conformal if 
\begin{equation}
\label{killing-1} \pounds_{X}\bar  g = 2\bar\varphi\, \bar g,
\end{equation}
is satisfied, where
\begin{equation}
\bar\varphi = \frac{1}{n+1}\,\textrm{div}_{\bar M} X.
\end{equation}
Furthermore, see section \ref{warped}, the conformal vector field $X$ is \emph{closed} if the $1$-form metrically equivalent to $X$ is closed. This amounts to the validity of the condition
\begin{equation}
\label{closedconf}
\bar\nabla_U X = \bar\varphi U, \,\, \mbox{ for all } \,\, U\in\Gamma(T\bar M).
\end{equation}
We refer the reader to Section \ref{warped} for further details on the geometry of $\bar M$ in the presence of conformal closed vector fields. 

We now consider some examples to illustrate Definition \ref{selfsimilar}.

\begin{example}
\label{spheres-example} Choose $\bar M = \mathbb{R}^{n+1}$ with
the Euclidean metric expressed outside the origin $o=(0,\ldots,0)\in \mathbb{R}^{n+1}$ in polar coordinates $(t, \theta)
\in (0,\infty) \times \mathbb{S}^n$ as $dt^2+ t^2 d\theta^2$, where $d\theta^2$ denotes the standard metric in $\mathbb{S}^n$. Let 
\[
X(x) = x,
\]
for  $x\in \mathbb{R}^{n+1}$ and
\[
\Phi(s,x) =e^s x
\]
for $(s,x) \in \mathbb{R}\times \mathbb{R}^{n+1}$, where $s=\log
t$. We note that 
\[
\Phi_*\frac{\partial }{\partial s}\Big|_{(s,x)} 
= X(\Phi(s,x))
\]
and that in terms  of polar coordinates one has
\[
X= t\partial_t.
\]
Fix $M=\mathbb{S}^n$ and  let $\psi:\mathbb{S}^n\to \mathbb{R}^{n+1}$ be the standard inclusion map. Choose a reparametrization
$\sigma:(\omega_*, \omega^*)\to (-\infty, \infty)$ of the form
\[
s = \sigma(\tau) = \log t(\tau)
\]
and define $\Psi: (\omega_*, \omega^*)\times \mathbb{S}^n \to \mathbb{R}^{n+1}$ as 
\[
\Psi (\tau, x) = \Phi(\sigma(\tau), \psi(x)) = t(\tau) \psi(x).
\]
This defines a mean curvature flow if and only if
\[
\frac{\partial\Psi}{\partial\tau} =
-\frac{n}{t(\tau)}\,\partial_t\big|_{\Psi(\tau,x)}
\]
for $(\tau, x)\in (\omega_*, \omega^*)\times\mathbb{S}^n$. However we have
\[
\frac{\partial \Psi}{\partial \tau}=\frac{\partial
\Phi}{\partial s}\frac{ds}{dt}\frac{dt}{d\tau}=
t(\tau)\frac{1}{t(\tau)}\frac{dt}{d\tau}\,\partial_t\big|_\Psi = \frac{dt}{d\tau}\,\partial_t\big|_\Psi 
\]
from which we infer
\[
\frac{dt}{d\tau} = -\frac{n}{t(\tau)}\cdot
\]
Solving this equation we find
\[
t(\tau) = \sqrt{c-2n\tau}
\]
for $c>0$ and  $\tau\in (-\infty, c/2n)$. We conclude that
\[
\Psi(\tau, x)= \sqrt{c-2n\tau}\,\psi(x),
\]
for $(\tau, x)\in (-\infty, c/2n)\times\mathbb{S}^n$, with
\[
\Psi(0,x) =\sqrt c\,\psi(x), \quad x\in\mathbb{S}^n,
\]
is a \emph{self-similar} mean curvature flow in $\mathbb{R}^{n+1}$
according to Definition \ref{selfsimilar}. It recovers the
well-known classical Euclidean example of \emph{self-shrinker} evolving concentric
spheres. The \emph{self-expander} evolving concentric spheres example correspond to invert the orientation of the flow lines, that is, to consider the vector field $X= -t\partial_t$.
\end{example}

\begin{example}
\label{warped-example} In general, consider a warped product  of
the form $\bar M = I\times_h P$, where $I$ is an open interval in
$\mathbb{R}$ and $P$ is a $n$-dimensional Riemannian manifold. We define a warped
metric in $\bar M$  expressed in coordinates $(t, x)
\in I \times P$ as $dt^2+ h^2(t) g_0(x)$, where $h$ is a positive
smooth function in $I$ and $g_0$ stands for the
metric in $P$. We refer the reader to Section \ref{warped} for
further details on warped product geometry.

In this setting, the  vector field $X= h(t)\partial_t$ is conformal and closed {\rm (}see equation {\rm (\ref{warped-conformal})} in Section \ref{warped}{\rm )} and the flow parameter is given by
\begin{equation}
\label{ts} s =\int\frac{dt}{h(t)}
\end{equation}
up to an additive integration constant. Indeed, the flow lines are described in terms of coordinates $(t,x)\in I\times P$ as
\[
\Phi(s,x) = (t(s), x)
\]
what implies that
\[
\Phi_* \frac{\partial}{\partial s} = \frac{dt}{ds} \partial_t |_\Phi = \frac{dt}{ds} \frac{1}{h(t)} X(\Phi(s,x)).
\]
Hence we have
\[
\frac{ds}{dt}\Big|_t =\frac{1}{h(t)}, \quad t\in I. 
\]
Having fixed $t_0\in I$, we define a
mean curvature flow in $\bar M$ by
\begin{equation}
\label{mcf-warped} \frac{d\Psi}{d\tau}\Big|_{(\tau, x)}
=-n\frac{h'(t(\tau))}{h(t(\tau))}\,\partial_t\big|_{\Psi(\tau,x)},
\end{equation}
where $(\tau, x)\in (\omega_*,\omega^*)\times P_{t_0}$ and $'$
indicates derivative with respect to $t$. Here $M=P_{t_0}$ denotes
the leaf $\{t_0\}\times P$ in $\bar M$ whereas the open interval
$(\omega_*,\omega^*)$ is the maximal interval of definition of the
solution of the equation
\[
\frac{dt}{d\tau}\Big|_\tau = -n\frac{h'(t(\tau))}{h(t(\tau))}.
\]
Indicating with $\tau(t)$ the inverse function of $t(\tau)$, one has
\[
n\tau(t)=-\int^t_{c} \frac{h(r)}{h'(r)}\,dr + n \tau(c),
\]
for some constant $c\in I$. It follows that the map $\Psi$ given
by
\[
\Psi(\tau, x) = \Phi(s(t(\tau)),x),
\]
for $(\tau, x)\in (\omega_*,\omega^*)\times P_{t_0}$, satisfies
${\rm (\ref{mcf-warped})}$ above with initial condition
\[
\Psi(0,x) =\Phi(s(c),P_{t_0}).
\]
We conclude that $\Psi$ defines a self-similar mean
curvature flow in the warped space $\bar M$. This example extends Example
\ref{spheres-example} in the general  setting  of warped product spaces.
Indeed, $\Psi$ is a mean curvature flow of the totally umbilical
leaves $P_t$. In particular, in the Euclidean case, where $P=\mathbb{S}^n$, we
have $h(t)=t$ and $P_t =\mathbb{S}^n (t)$.
\end{example}

\begin{example}\label{horosphere-example}
As a special case of the previous example, we consider the
hyperbolic space $\mathbb{H}^{n+1}$ of constant sectional curvature $-1$, described as the warped
product space $\mathbb{R}\times_h \mathbb{R}^n$, where $h(t)=e^t$, $t\in
\mathbb{R}$. In this particular description, the leaves of the natural foliation $t\mapsto \{t\}\times \mathbb{R}^n$ are
horospheres. We have $X=e^t\partial_t$ and the flow parameter is chosen to be
\[
s =C-\frac{1}{e^{t}},\quad t\in\mathbb{R},
\]
for some constant $C$ and $s\in (-\infty, C)$. We may
define a mean curvature flow by setting
\[
\frac{d\Psi}{d\tau}\Big|_{(\tau, x)} = -n\,\partial_t\big|_{\Psi(\tau,x)}
\]
for $(\tau, x)\in \mathbb{R}\times P_{t_0}$, where
$P_{t_0}=\{t_0\}\times\mathbb{R}^n$ for some $t_0\in\mathbb{R}$. This gives a self-similar mean curvature flow 
\[
\Psi(\tau,x) =\Phi(s(t(\tau)), x),
\]
for $(\tau, x)\in \mathbb{R}\times P_{t_0}$, where
\[
t(\tau)=c-n\tau,
\]
for some constant $c$. In this case the initial condition is
$\Psi(0,P_{t_0})=P_{c}$ and the time-slices are horospheres.
\end{example}

\begin{example}\label{geodsphere-example}
Having fixed an origin $o\in \mathbb{H}^{n+1}$, we now consider the model of  $\mathbb{H}^{n+1}\backslash\{o\}$ as the warped
product $(0,\infty)\times_h \mathbb{S}^n$, where $h(t)=\sinh t$,
$t\in (0,\infty)$. In this particular description, the leaves of the corresponding natural foliation are
geodesic spheres. We have $X=\sinh t\partial_t$ and the flow
parameter is
\[
s =\log\Big(\frac{e^t-1}{e^t+1}\Big)+C,\quad t>0,
\]
for some constant $C$ and $s\in \mathbb{R}$. In this model, we define a self-similar 
mean curvature flow by
\[
\Psi(\tau,x) =\Phi(s(t(\tau)), x),
\]
for $(\tau, x)\in (-\infty, c/n)\times P_{t_0}$, where
\[
\cosh t(\tau)=e^{c-n\tau},
\]
for a constant $c$. One concludes that with this choice of the
function $t=t(\tau)$ the flow  has  initial condition
$\Psi(0,P_{t_0})=P_{c'}$, where $c'=\log(e^c+\sqrt{e^{2c}-1})$
and whose time-slices are geodesic spheres.
\end{example}

\begin{example}\label{hyp-example}
Foliating  $\mathbb{H}^{n+1}$ by equidistant hypersurfaces yields a different warped product model of it, namely, 
as the product $\mathbb{R}\times_h \mathbb{H}^n$, where $h(t)=\cosh t$,
$t\in \mathbb{R}$.  With these choices  we have $X=\cosh t\,\partial_t$ and the flow
parameter is
\[
s =2 \arctan e^t+C,
\]
for some constant $C$ and $s\in \mathbb{R}$. We now define a self-similar 
 mean curvature flow by
\[
\Psi(\tau,x) =\Phi(s(t(\tau)), x),
\]
for $(\tau, x)\in \mathbb{R}\times P_{t_0}$, where
\[
\sinh t(\tau)=e^{c-n\tau},
\]
for a constant $c$. One concludes that with this choice of the
function $t=t(\tau)$ one indeed defines
\emph{two} self-similar
mean curvature flows in $\mathbb{H}^{n+1}$ with
initial conditions $\Psi(0,P_{t_0})=P_{c'}$, where $c'=\log(e^c\pm\sqrt{e^{2c}+1})$
and whose time-slices are equidistant hypersurfaces.
\end{example}

\begin{example}
\label{spheres} 
The Euclidean sphere  $\mathbb{S}^{n+1}$ may also be described, outside a pair of antipodal points $\{o,-o\}\subset \mathbb{S}^{n+1}$, as the warped
product  space $(0,\pi)\times_h \mathbb{S}^n$, with $h(t)=\sin t$,
$t\in (0, \pi)$. In this particular description, the leaves of the corresponding natural foliation are
geodesic spheres and we set $X=\sin t\,\partial_t$.
We may define a self-similar
 mean curvature flow by
\[
\Psi(\tau,x) =\Phi(s(t(\tau)), x),
\]
for $(\tau, x)\in (-\infty, -c/n)\times P_{t_0}$, where
\[
\cos t(\tau)=e^{c+n\tau},
\]
for a constant $c$. With these choices this flow has initial condition $\Psi(0,P_{t_0})=P_{c'}$, where $c'=\log(e^c+\sqrt{e^{2c}-1})$ and geodesic spheres as time-slices.
\end{example}

\begin{example}\label{cone-example} {\rm (Cone manifolds, \cite{Futaki}.)}  Suppose that $\bar M$ is a {\rm cone
manifold} of the form $\bar M = (0,\infty)\times_t P$ according to the
notation in Example \ref{warped-example} above. Then the warped product
metric is expressed as $dt^2+ t^2 g_0$, where $t$ is the
natural parameter in $(0,\infty)$ and  $g_0$ stands
for the metric in $P$. The mean curvature flow in cone manifolds had been previously studied in {\rm \cite{Futaki}}, where one can find a detailed description of parabolic rescalings and the behavior of type I singularities.
\end{example}

\begin{example}\label{product-example} {\rm (Translating self-similar flows.)} We consider the case of a parallel vector field. Then
$\bar\varphi=0$. In particular, suppose that $\bar M$ is a Riemannian
product of the form $\bar M = I\times P$ according to the notation
in Example \ref{warped-example} above. Then the product metric is
expressed as $dt^2+ g_0$, where $t$ is the natural parameter
in $I\subset \mathbb{R}$ and  $g_0$ stands for the metric in
$P$.
In this setting, the  vector field $X=\partial_t$ is 
parallel. Moreover the flow parameter is given by
\begin{equation}
\label{ts-bis} s =t
\end{equation}
up to an additive integration constant. Fixed $t_0\in I$, we may define a
mean curvature flow in $\bar M$ by
\[
\Psi(\tau, x) = (c, x), \quad \tau\in \mathbb{R},
\]
for any $x\in P_c=\{c\}\times P$.  Less trivial examples may be given by evolving translating graphs: consider a function $u: (\omega_*, \omega^*)\times M  \to I$ and define $\Psi: (\omega_*, \omega^*)\times M \to I$ as
\[
\Psi(\tau, x) = (u(\tau, x), x).
\]
This defines a mean curvature flow if and only if $u$ satisfies the quasilinear parabolic equation
\begin{equation}
\label{pde-parabolic}
\frac{\partial u}{\partial \tau} =W \, {\rm div}_P \Big(\frac{\nabla^P u}{W}\Big),
\end{equation}
with
\[
W = \sqrt{1+|\nabla^P u|^2}
\]
where $\nabla^P$ and ${\rm div}_P$ are, respectively, the Riemannian connection and divergence in  $(P, g_0)$. This notion of translating soliton has been extensively studied in Euclidean spaces, see for instance {\rm \cite{huisken2}, \cite{AltschulerWu}, \cite{Clutterbuck},  \cite{martin}, \cite{martin2}, \cite{nguyen}, \cite{shariyari}, \cite{xin}, \cite{wang}} only to quote a few examples of the vast literature on the subject. Our definition is the natural setting to these special flows in Riemannian products $I \times P$ with $I\subset \mathbb{R}$.
\end{example}

Next, we present some fundamental consequences of Definition \ref{selfsimilar} that motivates the notion of mean curvature flow soliton in a general geometric
setting. 

\begin{proposition}
\label{soliton-prop} Let $\Psi:(\omega_*, \omega^*)\times M \to
\bar M$ be a self-similar mean curvature flow with respect to some  vector field $X\in \Gamma(T\bar M)$. Then for all
$\tau\in (\omega_*, \omega^*)$ there exists a constant $c_\tau$
such that
\begin{equation}
\label{solitonA} c_\tau X =c_\tau\Psi_{\tau *}T + {\bf H},
\end{equation}
where ${\bf H}$ is the mean curvature vector of $\Psi_\tau= \Psi(\tau, \cdot\,)$ and
$T\in\Gamma(TM)$ is the pull-back by $\Psi_\tau$ of the tangential
component of $X$. Furthermore, 
\begin{equation}
\label{solitonC} \bar\nabla^\perp {\bf H}+c_\tau\, II(T,\cdot)=0,
\end{equation}
where $II$ is the second fundamental form of $\Psi_\tau$ and $\bar\nabla^\perp$ is its normal connection.  
Moreover, if $X$ is conformal
\begin{equation}
\label{solitonB} II_{-\mathbf{H}} + \frac{c_\tau}{2}\pounds_T g =
c_\tau\varphi g, \quad \varphi = \bar\varphi \circ\Psi_\tau
\end{equation}
where $g$ is the metric induced in $M$ by $\Psi_\tau$ and $II_{-\bf
H}$ is its second fundamental form in the direction of $-{\bf H}$. 
\end{proposition}

\noindent \emph{Proof.} Differentiating both sides in
(\ref{psiphi}) with respect to $\tau$ we obtain
\begin{eqnarray}
\frac{d\Psi}{d\tau}\Big|_{(\tau,x)}&=&
\frac{\partial\Phi}{\partial s} \Big|_{(\sigma
(\tau),\widetilde\Psi(\tau, x))}\frac{d\sigma}{d\tau}\Big|_{\tau}
+\Phi_{\sigma(\tau)*}(\widetilde\Psi(\tau, x))\frac{d\widetilde\Psi}{d\tau}\Big|_{(\tau, x)}\nonumber\\
& = & X(\Phi(\sigma(\tau), \widetilde\Psi(\tau,
x)))\frac{d\sigma}{d\tau}\Big|_{\tau}+\Phi_{\sigma(\tau)*}(\widetilde\Psi(\tau,
x))
\frac{d\widetilde\Psi}{d\tau}\Big|_{(\tau, x)}\nonumber\\
& = & X(\Psi(\tau,
x))\frac{d\sigma}{d\tau}\Big|_{\tau}+\Phi_{\sigma(t)*}(\widetilde\Psi(\tau,
x))\frac{d\widetilde\Psi}{d\tau}\Big|_{(\tau, x)},\label{XT}
\end{eqnarray}
where $\Phi_\sigma = \Phi(\sigma, \cdot).$ Since $\Psi$ is a
self-similar mean curvature flow with respect to $X$ there exists an isometric
immersion $\psi:M\to \bar M$ such that $\Psi(0,\,\cdot\,)=\psi$
and
\begin{equation}
\label{soliton-definition-2} \widetilde\Psi_{\tau} (M) = \psi(M),
\end{equation}
for all $\tau\in (\omega_*,\omega^*)$. That is,
\[
\Psi_\tau (M) = \Phi(\sigma(\tau), \psi (M)).
\]
Equation (\ref{soliton-definition-2})  implies that
\[
\frac{d\widetilde\Psi}{d\tau}\Big|_{(\tau, x)} \in
T_{\widetilde\Psi (\tau, x)} \psi(M)
\]
and
\[
\Phi_{\sigma(\tau)*}(\widetilde\Psi(\tau,
x))\frac{d\widetilde\Psi}{d\tau}\Big|_{(\tau, x)} \in T_{\Psi
(\tau, x)} \Psi_\tau(M).
\]
We conclude that for all $\tau\in (\omega_*,\omega^*)$ the
tangential component of $\frac{d\sigma}{d\tau} X$ onto
$\Psi_\tau(M)$ is given by
\[
X^\top(\Psi(\tau,x))\frac{d\sigma}{d\tau}\Big|_{\tau} =
-\Phi_{\sigma(\tau)*}(\widetilde\Psi(\tau,x))\frac{d\widetilde\Psi}{d\tau}\Big|_{(\tau,x)}
\]
where the superscript $\top$ denotes tangential projection.
We note that the  expression
\begin{equation}
c_\tau \Psi_{\tau*}T(\tau,x) =
-\Phi_{\sigma(\tau)*}(\widetilde\Psi(\tau,x))\frac{d\widetilde\Psi}{d\tau}\Big|_{(\tau,x)}.
\end{equation}
defines a vector field $T(\tau,\cdot)\in \Gamma(TM)$,
for each $\tau\in (\omega_*,\omega^*)$, where
\[
c_\tau = \frac{d\sigma}{d\tau}\Big|_\tau.
\]
We conclude from (\ref{XT}) that
\begin{equation}
\label{soliton-quase}
{\bf H}|_{\Psi(t,x)} =
c_\tau X(\Psi(t,x))
-c_\tau\Psi_{\tau*}T(\tau,x).
\end{equation}
Then, we rewrite (\ref{soliton-quase}) in the form
\begin{equation}
\label{soliton-tau}
c_\tau X|_{\Psi_\tau} = 
c_\tau\Psi_{\tau *}T_\tau   + {\bf H}|_{\Psi_\tau} 
\end{equation}
where $T_\tau(x)= T(\tau,x)$.

Next, for a fixed $\tau$  one computes along the immersion
$\Psi_\tau$ with induced metric $g_\tau$ obtaining
\begin{eqnarray*}
\pounds_T  g_\tau (U,V)  = \langle \nabla_U T, V\rangle +
\langle U, \nabla_VT\rangle
\end{eqnarray*}
where $\nabla$ is the induced connection in $(M, g_\tau)$. Taking traces yields
\begin{eqnarray*}
\textrm{tr}_{g_\tau} \pounds_T g_\tau = 2\,\textrm{div}_{g_\tau}T.
\end{eqnarray*}
Using (\ref{soliton-tau}) one obtains
\begin{equation}
\label{nablaX} c_\tau\bar\nabla_{\Psi_* U} X = c_\tau\Psi_{*}\nabla_U T
+ c_\tau(\bar\nabla_{\Psi_* U} \Psi_*T)^\perp +\bar\nabla_{\Psi_*U}
\mathbf{H}.
\end{equation}
Taking the normal projection in both sides one has
\begin{equation*}
\label{nablaX-bis} c_\tau(\bar\nabla_{\Psi_* U} X)^\perp = c_\tau (\bar\nabla_{\Psi_* U} \Psi_*T)^\perp +\bar\nabla^\perp_{\Psi_*U} \mathbf{H}.
\end{equation*}
If $X$ is closed and conformal we have from (\ref{closedconf}) above that
\begin{equation*}
\bar\nabla_{\Psi_* U} X = \varphi\, \Psi_* U, \,\, \mbox{ for all }\,\, U\in \Gamma(TM)
\end{equation*}
where $\varphi = \bar\varphi\circ\Psi_\tau$  which yields
\begin{equation*}
c_\tau (\bar\nabla_{\Psi_* U} \Psi_*T)^\perp
+\bar\nabla^\perp_{\Psi_*U} {\bf H}=0,
\end{equation*}
that is, 
\begin{equation}
c_\tau\, II(T,U) + \bar\nabla_{\Psi_* U}^\perp {\bf H} =0.
\end{equation}
Next, denoting by $II_{-\mathbf{H}}$ the second fundamental form of
$\Psi_\tau$ in the opposite direction of the mean curvature vector field
$\mathbf{H}$, one deduces from  (\ref{killing-1}),  that is, from the fact that $X$ is conformal, 
\begin{eqnarray*}
& & 2c_\tau \varphi \langle \Psi_*U, \Psi_*V\rangle = c_\tau\langle \bar\nabla_{\Psi_*U} X, \Psi_*V\rangle + c_\tau\langle\Psi_*U, \bar\nabla_{\Psi_*V} X \rangle \\
&  &\,\, =  c_\tau\langle \nabla_U T, V\rangle +c_\tau \langle U, \nabla_V T\rangle + \langle\bar\nabla_{\Psi_* U}\mathbf{H},\Psi_* V\rangle + \langle\Psi_* U, \bar\nabla_{\Psi_* V}\mathbf{H}\rangle\\
&  &\,\, =  c_\tau\langle \nabla_U T, V\rangle + c_\tau\langle U, \nabla_V T\rangle - \langle\bar\nabla_{\Psi_* U}\Psi_* V, \mathbf{H}\rangle - \langle \mathbf{H}, \bar\nabla_{\Psi_* V}\Psi_* U\rangle\\
& & \,\,=c_\tau\,\pounds_T  g (U,V) +2 II_{\mathbf{-H}}(U,V),
\end{eqnarray*}
where we have omitted the subscript $\tau$ for the sake of brevity.  We
then have proved that $\Psi_\tau$ satisfies the soliton equation
\begin{equation}
\label{soliton-2} II_{-\mathbf{H}} + \frac{c_\tau}{2}\pounds_T g =
c_\tau\bar\varphi\, g.
\end{equation}
This completes the proof of Proposition \ref{soliton-prop}. \hfill $\square$

\section{Generalized mean curvature flow solitons}\label{fundamental}

Motivated by the above geometric setting, we define a general notion of mean curvature flow soliton  with respect to a given vector field $X\in \Gamma(T\bar M)$ as follows.

\begin{definition}\label{soliton-definition}
An isometric immersion $\psi:M^m\to \bar M^{n+1}$ is a \emph{mean
curvature flow soliton}  with respect to $X\in \Gamma(T\bar M)$  if
\begin{equation}
\label{solitonA-2}c\, X^\perp = {\bf H}
\end{equation}
along $\psi$ for some constant $c\in \mathbb{R}$. 
\end{definition}
With a slight abuse of notation, we also say that
the submanifold  $\psi(M)$ itself is the mean curvature flow
soliton {\rm (}with respect to the vector field $X${\rm )}. If $m=n$, that is, for codimension $1$, the condition
becomes
\begin{equation}
\label{soliton-scalar} H=c\,\langle X, N\rangle,
\end{equation}
where the mean curvature $H$, with respect to the local normal vector field $N$ along $\psi$, is given by
\begin{equation}
\label{HHN}
{\bf H} = HN.
\end{equation}

We observe that in case $X$ is a conformal vector field on $\bar M$, equation (\ref{solitonA-2}) is enough to deduce the following important consequences that we have considered in Proposition \ref{soliton-prop} in a different setting. 



\begin{proposition}
\label{soliton-prop-2} Let $\psi:M^m\to \bar M^n$ be a mean curvature
flow soliton with respect to a conformal vector field  $X\in \Gamma(T\bar M)$. Then along $\psi$ we have
\begin{equation}
\label{solitonB-2} II_{-\mathbf{H}} + \frac{c}{2}\pounds_{T} g =
c\,\varphi\, g
\end{equation}
where $g$ is the metric induced in $M$ by $\psi$ and $II_{-\bf H}$
is its second fundamental form in the direction of $-{\bf H}$. Here   the vector field $T$ is defined by $\psi_*T = X^\top$ and 
\begin{equation}
\label{barvarphi}
\varphi = \frac{1}{n+1}\, {\rm div}_{\bar M} X\circ\psi =\bar\varphi\circ\psi.
\end{equation}
Furthermore, if $X$ is also closed
\begin{equation}
\label{solitonC-2} \bar\nabla^\perp {\bf H}+c\, II(T,\cdot)=0,
\end{equation}
where $II$ is the second fundamental tensor of $\psi$ and $\bar\nabla^\perp$ is its normal connection.
\end{proposition}

\noindent \emph{Proof.}  Using (\ref{solitonA-2}) by a direct computation we have for any tangent vector fields $U, V\in \Gamma(TM)$
\begin{eqnarray*}
& & 2c\varphi\, g (U, V)=c\langle \bar\nabla_{\psi_* U} X, \psi_* V\rangle+c\langle \bar\nabla_{\psi_* V} X, \psi_* U\rangle\\
& &\,\, = c\langle \bar\nabla_{\psi_* U} \psi_* T, \psi_* V\rangle+c\langle \bar\nabla_{\psi_* V} \psi_* T, \psi_* U\rangle + \langle \bar\nabla_{\psi_* U} {\bf H}, \psi_* V\rangle+\langle \bar\nabla_{\psi_* V} {\bf H}, \psi_* U\rangle\\
& &\,\, = c\langle \nabla_{U} T, V\rangle + c\langle \nabla_{V} T, U\rangle+ 2 II_{-\bf H} (U,V).
\end{eqnarray*}
Hence,
\[
II_{-\bf H} +\frac{c}{2}\pounds_{T} g = c\varphi g.
\]
Now, if $X$ is also closed,  using (\ref{closedconf}) one has
\begin{eqnarray*}
& & 0 =c ( \varphi\,  \psi_* U)^\perp = c(\bar\nabla_{\psi_* U} X)^\perp  = c (\bar\nabla_{\psi_*U}\psi_*T)^\perp + c(\bar\nabla_{\psi_* U} X^\perp)^\perp\\
& & \,\, = c\,II (T, U) + (\bar\nabla_{\psi_*U}{\bf H})^\perp = c\, II (T, U) + \nabla_{\psi_*U}^\perp{\bf H}
\end{eqnarray*}
what concludes the proof of Proposition \ref{soliton-prop-2}. \hfill $\square$

\begin{remark}
\label{ricci-mcf} There is a strict similarity between  {\rm (\ref{solitonB-2})} and the equation defining Ricci solitons. Indeed this latter can be expressed in the form
\begin{equation}
\label{ric-sol}
{\rm Ric}_M + \frac{1}{2}\pounds_U g = \lambda g,
\end{equation}
for some $U\in \Gamma(TM)$ and $\lambda\in \mathbb{R}$. Note that if the Ricci soliton is gradient, that is, $U =\nabla f$ for some $f\in C^\infty(M)$ we deduce from  {\rm (\ref{ric-sol})}  the fundamental equation
\begin{equation}
\label{nabla-scalar}
\frac{1}{2}\nabla R  = {\rm Ric}_M (\nabla f, \cdot)^\sharp
\end{equation}
where $R = {\rm tr}_g {\rm Ric}_M$ is the scalar curvature and $\sharp$ denotes the musical isomorphism.  In {\rm (\ref{solitonB-2})} the tensor $II_{-{\bf H}}$ plays the role of ${\rm Ric}_M$. Its trace is $-|{\bf H}|^2$ and from {\rm (\ref{solitonC-2})}  we immediately deduce 
\begin{equation}
\label{solitonD-2}
\frac{1}{2}\nabla |{\bf H}|^2 =c\, II_{-{\bf H}} (T, \cdot)^\sharp
\end{equation}
which has the same structure of {\rm (\ref{nabla-scalar})}. 
\end{remark}

\begin{remark} In case of codimension $1$  we shall characterize {\rm (\ref{soliton-scalar})}  in Section \ref{variational} from a variational point of view al least in case $\bar M$ has the structure of a warped product $\bar M = I\times_h P$.  In fact    {\rm (\ref{soliton-scalar})}  characterizes weighted minimal hypersurfaces  with respect to  a weight depending on $X=h(t)\partial_t$ in the notation of Example \ref{warped-example} above.
\end{remark}

We now give some examples.

\begin{example}\label{spheres-example-2}
Consider the Euclidean self-similar mean curvature flow $\Psi$ of
concentric spheres in $\bar M = \mathbb{R}^{n+1}$ given in Example
\ref{spheres-example}. We fix a particular value 
$\bar\tau$ of the parameter $\tau$ and denote
\[
\frac{d\sigma}{d\tau}\Big|_{\bar\tau}=c.
\]
Then,
\[
c=\frac{ds}{dt}\Big|_{t(\bar\tau)}\frac{dt}{d\tau}\Big|_{\bar\tau}=
-\frac{n}{t^2(\bar\tau)}
\]
Hence $\bar\tau$ is given implicitly  by the radius
\[
t(\bar\tau)=\sqrt{-n/c}.
\]
We conclude that the sphere $\mathbb{S}^n(\sqrt{-n
/c})\subset\mathbb{R}^{n+1}$ with radius $\sqrt{-n/c}$ is a mean
curvature soliton with respect to the radial vector field $X=t\partial_t$ according to Definition \ref{soliton-definition}. This
 agrees with the usual notion of self-shrinker established in
the literature, {\rm \cite{Ecker-book}}.
\end{example}

\begin{example}
\label{warped-example-2} In a warped product space $\bar M = I\times_h P$,
consider the self-similar mean curvature flow $\Psi$ defined in
Example \ref{warped-example} above with $X= h(t)\partial_t$. Fixed a 
particular value $\bar\tau\in (\omega_*, \omega^*)$ of the
parameter $\tau$ we denote
\[
\frac{d\sigma}{d\tau}\Big|_{\bar\tau}=c.
\]
Hence,
\[
c=\frac{ds}{dt}\Big|_{t(\bar\tau)}\frac{dt}{d\tau}\Big|_{\bar\tau}=
-n\frac{h'(t(\bar\tau))}{h^2(t(\bar\tau))}
\]
and implies that $\bar\tau$ is implicitly given by the condition
\begin{equation}
\label{h-soliton} c h^2(t(\bar\tau))+nh'(t(\bar\tau))=0.
\end{equation}
It follows that the leaf $P_{t(\bar\tau)}$ corresponding to
$t=t(\bar\tau)$ is a mean curvature flow soliton  with respect to the vector field $X=h(t)\partial_t$ according to
Definition \ref{soliton-definition}. In fact we note that {\rm (\ref{h-soliton})} is exactly the scalar soliton
equation {\rm (\ref{soliton-scalar})} since in this case
\[
H =-n\frac{h'(t)}{h(t)}
\]
and
\[
\langle X, N\rangle|_{P_t} = h(t)
\]
for a given slice $P_t$ in $I\times_h P$.
\end{example}

\begin{example}\label{horosphere-example-2}
For the particular case of the hyperbolic space $\mathbb{H}^{n+1}$
described as the warped product $\mathbb{R}\times_h \mathbb{R}^n$,
where $h(t)=e^t$, $t\in \mathbb{R}$, let $\Psi$ be the
self-similar mean curvature flow in Example
\ref{horosphere-example}.  In this case, the horosphere $P_{t(\bar\tau)}$ with $\bar\tau$ given implicitly by
\[
t(\bar\tau)=\log (-n/c)
\]
is a mean curvature soliton with respect to the vector field $X=e^t\partial_t$ according to
Definition \ref{soliton-definition}.
\end{example}

\begin{example}\label{geodsphere-example-2}
Now, we consider the self-similar mean curvature flow $\Psi$ in
$\mathbb{H}^{n+1}\backslash\{o\}$ as the warped product
$(0,\infty)\times_h \mathbb{S}^n$, where $h(t)=\sinh t$, $t\in
(0,\infty)$ defined in Example \ref{geodsphere-example}.  In this case, the geodesic sphere $P_{t(\bar\tau)}$ with $\bar\tau$ given implicitly by
\[
\cosh t(\bar\tau)=\frac{1}{2c}(-n\pm\sqrt{n^2+4c^2}).
\]
is a mean curvature soliton with respect to $X=\sinh t\partial_t$  according to
Definition \ref{soliton-definition}.
\end{example}

\begin{example} {\rm (Solitons in cone manifolds, \cite{Futaki}.)}  Suppose that $\bar M$ is a cone
manifold of the form $\bar M = (0,\infty)\times_t P$ as in Example \ref{cone-example} above. Then it follows from Example \ref{warped-example-2} that the hypersurface $P_{\sqrt{-\frac{n}{c}}} = \{\sqrt {-n/c}\}\times P$ is a mean curvature flow soliton.  In {\rm \cite{Futaki}}, the authors prove that a suitable sequence of parabolic rescalings of the mean curvature flow around a type $I_{c}$ singularity converge to a soliton. This kind of singularity is an extension of the usual definition of the type $I$ singularity to cone manifolds. A key tool for the proof is a variant of the classical Huisken's monotonicity formula valid in this context. 
\end{example}

\begin{example}\label{translating-solitons} {\rm (Translating solitons.)} In case $X$ is a parallel vector field we consider the setting in Example \ref{product-example}.  We can suppose that $\bar M$ is a Riemannian
product of the form $\bar M = I\times P$ with $X = \partial_t$. A mean curvature flow soliton, in this case named a \emph{translating soliton}, can be described non-parametrically as the graph $\Gamma_u$ of a solution $u: P\to I$ of the partial differential equation
\begin{equation}
\label{pde-translating}
{\rm div}_P \bigg(\frac{\nabla^P u}{W}\bigg) = \frac{c}{W},
\end{equation}
for some constant $c\in \mathbb{R}$, with $W=\sqrt{1+|\nabla^Pu|^2}$. Some examples of such translating solitons in this more general context have been presented and characterized in {\rm \cite{ML}}. 
\end{example}

\section{Conformal fields and warped product spaces}\label{warped}

From now on, we suppose that
 $X$ is a  \emph{closed} conformal  vector field on $\bar M$. This means that its metrically
equivalent $1$-form is closed. As a consequence of (\ref{killing-1}) in this case we have, besides conformality,
\begin{equation}
\label{closedconf-bis}
\bar\nabla_U  X =\bar\varphi\, U,
\end{equation}
for any vector field $U\in \Gamma(T\bar M)$.  In particular, 
\begin{equation}
\bar\nabla |X| = \frac{\bar\varphi}{|X|}X
\end{equation}
whenever $|X|\neq 0$. We suppose that there are no singular points of $X$ in
$\bar M$ by replacing $\bar M$ with a proper open subset of it,
if necessary. 

It turns out that  if $\Phi:(\Omega_*, \Omega^*)\times \bar M \to \bar M$ is the flow generated by $X$, then $\Phi_s = \Phi(s, \cdot):\bar M \to \bar M$, $s\in (\Omega_*, \Omega^*)$, is a conformal map in the sense that  there
exists a smooth positive function $\lambda:(\Omega_*,
\Omega^*)\times\bar M\to\mathbb{R}$ such that
\begin{equation}
\label{conformal} \Phi^*_s \bar g|_x = \lambda^2(s,x)\, \bar g|_x
\end{equation}
for all $x\in \bar M$. It follows from (\ref{killing-1}) that
\begin{equation}
\label{philambda} \bar\varphi(\Phi(s,x)) =\lambda(s,x)\,\partial_s
\lambda (s,x),
\end{equation}
for $(s,x)\in (\Omega_*, \Omega^*)\times\bar M$.  For our purposes, it is convenient to parameterize the flow $\Phi$ by fixing initial conditions on the fixed leaf $P$, that is, we consider $\Phi$ as
a global chart of $\bar M$ of the form
\begin{equation}
\Phi:(\Omega_*, \Omega^*)\times P\to \bar M.
\end{equation}
Having fixed this map, the hypersurface $P$, identified with the slice
$\{0\}\times P\subset (\Omega_*, \Omega^*)\times P$, is preserved
by the flow. In general, the integral leaves are given by
$P_s:=\Psi_s(P)$ and identified with the slice $\{s\}\times P$.
The principal
curvatures of an integral leaf $P_s$ with respect to $-X/|X|$ are
given by $\varphi/|X|$ and its mean curvature is 
\begin{equation}
\label{warped-H} \mathcal{H}(\Phi(s,x))=
-n\frac{\varphi}{|X|}\frac{X}{|X|}\Big|_{\Phi(s,x)}.
\end{equation}
We then consider the particular case when the conformal factor
depends only on the flow parameter, that is,
\[
\lambda=\lambda(s).
\]
In this case, each leaf $P_s$ is homothetic to $P$. This
particular case corresponds to warped product spaces. More precisely,
given the change of variables
\[
t=\int h(s)\,ds,
\]
we can describe $\bar M$ as a product $I\times P$ with warped
Riemannian metric given by
\[
\textrm{d}t^2 + h^2(t) \, g_0,
\]
where $g_0$ is the metric in $P$ and
\[
h(t)=|X|(s(t)).
\]
In this case we have
\begin{equation}
\label{conformal-warped}
X =h(t) \partial_t
\end{equation}
and
\begin{equation}
\label{warped-conformal}
\bar\nabla_{U}X = h'(t) U,
\end{equation}
that is, $\varphi=h'$ in this case.

\begin{remark}
Examples \ref{spheres-example}-\ref{product-example} in Section \ref{ssf} are particular cases of warped product spaces including spaces of constant sectional curvature, cones and Riemannian products. Other well-known examples are Riemannian Schwzarschild and Reissner-N\"ordstrom spaces of General Relativity.
\end{remark}

For later use, we observe that in the case where $P$, the integral leaf identified to $\{0\}\times
P$ through the flow $\Phi$,  is a Riemannian manifold with constant
sectional curvature $\kappa$, the Riemann curvature tensor $\bar R$ of $\bar M = I\times_h P$ is 
\begin{eqnarray}
& & \langle \bar R(U,V)W,Z\rangle = 
\bigg(\frac{h'^2}{h^2}-\frac{\kappa}{h^2}\bigg)\big(\langle \hat{U},\hat{W}\rangle \langle \hat{V},\hat{Z}\rangle - \langle \hat{V},\hat{W}\rangle \langle \hat{U},\hat{Z}\rangle\big)\nonumber\\
& & \,\, + \frac{h''}{h}\langle U, \partial_t\rangle\langle W, \partial_t\rangle \langle \hat{Z}, \hat{V}\rangle- \frac{h''}{h}\langle U, \partial_t\rangle\langle Z, \partial_t\rangle \langle \hat{W},\hat{V}\rangle \nonumber\\
& &\,\,\,\,- \frac{h''}{h}\langle V, \partial_t\rangle\langle W,
\partial_t\rangle  \langle \hat{Z}, \hat{U}\rangle+ \frac{h''}{h}\langle
V, \partial_t\rangle\langle Z,
\partial_t\rangle \langle \hat{W},
\hat{U}\rangle,
\label{riemann-warped}
\end{eqnarray}
where, given a vector field $U\in \Gamma(T\bar M)$ one
denotes
\[
\hat{U} = U -\langle U, \partial_t\rangle \partial_t .
\]

\section{Some fundamental elliptic equations}
\label{sec5}

In  case of a warped product space target $\bar M^{n+1} = I \times_h P^n$  we define the function
\begin{equation}
\label{inth} \hat\eta(t) =\int^{t}_{t_0} h(s) \,ds, \quad (t_0, t)\subset I,
\end{equation}
where $t_0\in I$ is arbitrarily fixed. Let $\pi:\bar M \to I$ be the projection $\pi(t,p) =t$, for $(t,p)\in I\times P$.  Given an isometric  immersion $\psi: M^m \to \bar M^{n+1}$ we define $\eta: M \to \mathbb{R}$ as the composition
\begin{equation}
\eta(x) = \hat\eta  (\pi\circ\psi(x)), \quad x\in M. 
\end{equation}
We also define $\varphi: M \to \mathbb{R}$ by
\begin{equation}
\label{phipsi}
\varphi(x) =  (\bar\varphi \circ \psi)(x), \,\, \mbox{ for }\,\, x\in M.
\end{equation}
where 
\begin{equation}
\bar\varphi = \frac{1}{n+1}\, {\rm div}_{\bar M} X
\end{equation}
is the function in (\ref{closedconf}) for a closed conformal field $X$. In particular, for warped product spaces 
$\bar M=I\times_hP$ with the choice $X=h(t)\partial_t$, we have $\bar\varphi = h'\circ\pi$ and $\varphi = h'(\pi\circ\psi)$. 

Given a function $\zeta\in C^1(M)$, in what follows we shall indicate by $\Delta_\zeta$ the operator acting on, say, $u\in C^2(M)$ as
\[
\Delta_{\zeta} u=
\Delta u-\langle \nabla\zeta, \nabla u\rangle,
\]
with $\nabla$ and $\Delta$, respectively, the Riemannian connection and Laplace-Beltrami operator in $M$.
\begin{proposition}
\label{laplace-eta} Let $\psi:M^m\to \bar M^{n+1}$ be a mean curvature flow soliton in $\bar M=I\times_h P$ with respect to the conformal vector field $X=h(t)\partial_t$. Then
\begin{equation}
\label{delta-eta} \Delta \eta =m\varphi +\frac{1}{c} | {\bf H}|^2 = mh'(\pi\circ\psi)+\frac{1}{c}|{\bf H}|^2
\end{equation}
and
\begin{equation}
\label{delta-eta-X} \Delta_{-c\eta}\eta = m \varphi +c|X|^2 = mh'(\pi\circ\psi)+ch^2 (\pi\circ\psi)
\end{equation}
on $M$.
\end{proposition}


\noindent \emph{Proof.} For the sake of simplicity, from now on we will identify the tangent spaces $T_x M$ and $T_{\psi(x)}\psi(M)$ for any $x\in M$. For any $U\in TM$ we have
\[
\langle \nabla\eta, U\rangle = h\langle \partial_t, U\rangle=
\langle X, U\rangle
\]
and therefore
\begin{equation}
\label{Luis16}
\nabla\eta= X^\top,
\end{equation}
where, from now on, ${}^\top$ and ${}^\perp$ denote, respectively, tangential and normal projections. Furthermore
\begin{eqnarray*}
\langle \nabla_U \nabla\eta, V\rangle =\langle \bar\nabla_U X,
V\rangle-\langle \bar\nabla_U  X^\perp, V\rangle =\varphi\,\langle U, V\rangle + \langle
II(U,V),X^\perp\rangle.
\end{eqnarray*}
Taking traces with respect to the induced metric in $M$ we obtain
\begin{equation}
\Delta \eta = m\varphi +\langle {\bf H}, X^\perp\rangle.
\end{equation}
Using the soliton equation (\ref{solitonA}) we conclude that
\begin{equation}
\Delta\eta = m\varphi +c|X^\perp|^2 = m\varphi +\frac{1}{c}|{\bf H}|^2 ,
\end{equation}
where $\varphi = \bar\varphi\circ\psi$. Now from
\begin{equation}
\langle \nabla \eta, X\rangle = |X^\top|^2,
\end{equation}
one gets
\begin{equation}
\Delta_{-c\eta} \eta = \Delta\eta +c\langle \nabla \eta, X\rangle =
m\varphi +c|X^\perp|^2 +c|X^\top|^2 = m\varphi +c|X|^2.
\end{equation}
Since $\varphi = h'(\pi\circ\psi)$ and $|X|_{\psi}=h(\pi\circ\psi)$, this finishes the proof of Proposition \ref{laplace-eta}.\hfill $\square$

\begin{remark}
We note that $\hat \eta$ is invertible. Hence if $\eta = \hat\eta(\pi\circ\psi)$ is constant then $\pi\circ\psi$ is constant and from {\rm (\ref{delta-eta-X})} we deduce that $mh'(\pi\circ\psi) +c h^2(\pi\circ\psi)\equiv 0$ on $M$ so that $\psi(M)$  is contained in a leaf $\{\bar t\}\times P$ with $\bar t$ implicitly given by $mh'(\bar t) + ch^2(\bar t)=0$. 
\end{remark}

\section{Weak maximum principle for solitons}\label{wmp}

One of the main analytical tools used in this paper is the weak maximum principle either for the  operator
\begin{equation}
\label{delta-u-w}
\Delta_U = \Delta - \langle U, \nabla\, \cdot\,\rangle
\end{equation}
for some $U\in \Gamma(TM)$, or 
\begin{equation}
\label{delta-nabla-f}
\Delta_f= \Delta - \langle \nabla f, \nabla \,\cdot \,\rangle
\end{equation}
for some $f\in C^1(M)$, where $\Delta$ and $\nabla$ are, respectively, the Beltrami-Laplace operator in $M$ and the Riemannian connection. A general discussion on the weak maximum principle for a very large class of operators can be found in \cite{AMR}; however, for the convenience of the reader and to introduce also some new results, we will present a few basic concepts and properties  restricting ourselves to the above type of operators acting on $C^2(M)$.  Of course the regularity of $U$ and  $f$ and the space of functions acted upon by the operator can be respectively relaxed and enlarged.

\begin{definition}
\label{super-harmonic}
Let $M$ be a {\rm (not necessarily complete)} Riemannian manifold and let $U\in\Gamma(TM)$ {\rm (}respectively, $f\in C^1(M)${\rm )}. We say that the weak maximum principle holds for the operator $\Delta_U$  {\rm (}respectively, $\Delta_f${\rm )} on $M$ if
for any $u\in C^2(M)$ with $u^* = \sup_M u <\infty$ and for each $\gamma < u^*$ we have
\[
\inf_{\Omega_\gamma}\Delta_U u \le 0,  
\,\,\,  {\rm \big(}respectively\,\,\, \inf_{\Omega_\gamma}\Delta_f u \le 0 {\rm )}, 
\]
where
\[
\Omega_\gamma =\{x\in M: u(x)>\gamma\}.
\]
\end{definition}

From Theorem 3.11 of \cite{PRS} we know that  the validity of the weak maximum principle for the operator $\Delta_f$ is  equivalent to its stochastic completeness, a probabilistic concept.   Therefore, as a consequence of a classical result of Khas'minskii,  a sufficient (but in this case also necessary  \cite{MV}) condition for the validity of the weak maximum principle in a complete manifold is  the existence of a function $v\in C^2(M)$ satisfying
\begin{equation}
\label{K-test}
\begin{cases}
\Delta_f v \le Av  & \mbox{ on } \quad M\backslash K,\\
v(x) \to +\infty   & \mbox{ as } \quad x\to\infty \,\, \mbox{ in }\,\, M,
\end{cases}
\end{equation} 
for some constant $A\in \mathbb{R}$ and compact set $K\subset M$. We note that the first condition in (\ref{K-test}) can be substituted, for instance, with $\Delta_f v \le A$ on $M\backslash K$. By Theorem 3.1  of \cite{AMR} the same result extends to the operator $\Delta_U$ for any $U\in \Gamma(TM)$. 

We also observe that stochastic completeness does not imply neither is implied by geodesic completeness. For instance $\mathbb{R}^2\backslash\{0\}$ with the flat metric is stochastically complete with respect to $\Delta$ but obviously not complete. For the converse it is enough to consider a model manifold with radial Ricci curvature decreasing sufficiently fast. In this respect, see Proposition 2.3 of \cite{AMR}. 

We will often use stochastic completeness a a natural assumption in some of our results. For a sound discussion on it we refer to Grigor'yan \cite{Gri}, and to \cite[Chapter 3]{PRS} and \cite[Chapter 4]{AMR} for its relation with the weak maximum principle. For the operator $\Delta_f$ the validity of the weak maximum principle, always in the complete case, can also be detected by the growth of the volume of geodesic balls as expressed  in a very general setting in Theorem 4.1 of \cite{AMR}.

Occasionally, for instance in Theorem \ref{hs-D}, we shall also use an equivalent form of the weak maximum principle that, following \cite{AMiR},  we call the open weak maximum principle. The latter can be stated as follows

\begin{definition}
\label{owmp}
We say that the open weak maximum principle holds for the operator $\Delta_U$ {\rm (}respectively, $\Delta_f${\rm )} on $M$ if for each $F\in C^0(\mathbb{R})$, for each open set $\Omega \subset M$ with $\partial\Omega \neq\emptyset$ and for each $v\in C^0(\bar\Omega)\cap C^2(\Omega)$ satisfying
\begin{equation}
\label{open-max-p}
\begin{cases}
\Delta_U v\ge F(v) & {\rm (}respectively,\,  \Delta_f v \ge F(v){\rm )} \quad \mbox{\rm on } \quad \Omega\\
\sup_\Omega v <+\infty, &
\end{cases}
\end{equation} 
we have that either
\[
\sup_\Omega v = \sup_{\partial\Omega} v
\]
or
\[ 
F(\sup_\Omega v)\le 0.
\]
\end{definition}

For the equivalence of the two forms of the weak maximum principle see for instance \cite{AMiR}.

It is apparent that the open form of the weak maximum principle is reminiscent of Ahlfors' parabolicity criterion for Riemann surfaces \cite{A}, Theorem 6c. In fact parabolicity for $\Delta_U$ (or $\Delta_f$), in the sense of the validity of a Liouville type theorem for bounded above $\Delta_U$-subharmonic functions  (or  $\Delta_f$-subharmonic functions), can be expressed as a stronger form of the weak maximum principle. Indeed, 

\begin{definition}
\label{parabolic} A Riemannian manifold $M$ is strongly parabolic with respect to 
 $\Delta_U$ {\rm (}respectively, $\Delta_f${\rm )} if for any non-constant $u\in C^2(M)$ with $u^*= \sup_M u <\infty$ and for each $\gamma < u^*$ we have
\[
\inf_{\Omega_\gamma} \Delta_U u <0,\,\,\,
{\rm (}\mbox{respectively } \,\,\,
\inf_{\Omega_\gamma} \Delta_f u < 0 {\rm )},
\]
with $\Omega_\gamma$ as above.
\end{definition}

In the case of our operators $\Delta_U$ and $\Delta_f$, the three forms of parabolicity, that is, Ahlfors, Liouville and strong parabolicity, are in fact equivalent. See Section 4.4 of \cite{AMR}. For later use we recall that Ahlfors parabolicity expresses as follows: $\Delta_U$ (respectively, $\Delta_f$) is Ahlfors parabolic on $M$ if for each open set $\Omega\subset M$ with $\partial\Omega\neq\emptyset$ and for each non-constant $v\in C^0(\bar\Omega)\cap C^2(\Omega)$ satisfying 
\begin{equation}
\label{ahlfors}
\begin{cases}
\Delta_U v\ge 0 & {\rm (}respectively,\,  \Delta_f v \ge 0{\rm )} \quad \mbox{\rm on } \quad \Omega\\
\sup_\Omega v <+\infty, &
\end{cases}
\end{equation} 
we have 
\[
\sup_\Omega v = \sup_{\partial\Omega} v.
\]

Of course, for parabolicity,  Khas'minskii  type test still applies appropriately stated in the following mildly stronger form.  

Let $M$ be complete and assume the existence of $v\in C^2(M)$ such that $v(x)\to+\infty$ as $x\to \infty$ in $M$, and
\begin{equation}
\label{K-test-par}
\begin{cases}
\Delta_U v \le 0  & \mbox{ if }\quad\,\,\, U \equiv 0 \mbox{ and }  
\\
\Delta_U v <0 
\quad & \mbox{ if } \quad\,\,\, U\not\equiv 0 \mbox{ on } M\backslash K,
\end{cases}
\end{equation} 
for some  compact set $K\subset M$. Then the operator $\Delta_U$ is parabolic on $M$ (and similarly for $\Delta_f$). 

\begin{remark}
\label{Luis-remark-parabolicity}
Always in case $M$ is complete a sufficient condition can be given for parabolicity with respect to $\Delta_f$ in terms of the growth of a weighted volume of the boundary of geodesic balls. More precisely, having fixed an origin $o\in M$ let
\begin{equation}
\label{weighted-volume}
{\rm vol}_f (\partial B_r) = \int_{\partial B_r} e^{-f}\, {\rm d}M,
\end{equation}
where $\partial B_r$ is the boundary of the geodesic ball $B_r$ centered at $o$ and of radius $r$. If
\[
\frac{1}{{\rm vol}_f (\partial B_r)} \notin L^1(+\infty)
\]
then $M$ is $\Delta_f$-parabolic. Other more elaborated results can be found in Chapter 4 of {\rm \cite{AMR}}.
\end{remark}

Here we give a simple specific sufficient condition for the validity of the weak maximum principle (or equivalently stochastic completeness) for the operator $\Delta_{-c\eta}$ on a complete mean curvature flow soliton in $\bar M = I\times_h P$. For $\psi: M^m \to \bar M^{n+1} = I\times_h P$ we recall that
\begin{equation}
\eta(x) = \int_{t_0}^{\pi(\psi(x))} h(s)\, {\rm d}s
\end{equation}
with $t_0\in I$.

\begin{theorem}
\label{cond-wmp} 
Let $\psi: M^m\to \bar M^{n+1} = I\times_h P$ be a complete mean curvature flow soliton with respect to $X= h(t)\partial_t$.  Assume that $P$ has constant sectional curvature $\kappa$. Let $II$ be the second fundamental form of the immersion and suppose that
\begin{equation}
\label{cond-wmp-1}
\Lambda = \frac{1}{m-1}\sup_M (|II|^2+ch'(\pi\circ\psi)-(m-1)\varkappa(\pi\circ\psi))<+\infty,
\end{equation}
where
\begin{equation}
\label{varkappa}
\varkappa(t) = \min\bigg\{-\frac{h''(t)}{h(t)}, \frac{\kappa}{h^2(t)}-\frac{h'^2(t)}{h^2(t)}\bigg\}, \quad t\in I.
\end{equation}
Then the weak maximum principle holds for the operator $\Delta_{-c\eta}$ on $M$.
\end{theorem}

\noindent \emph{Proof.}  Let $\{x^i\}_{i=1}^m$ be local coordinates in $M$ and let ${\bf H}$ be the mean curvature vector field of the immersion $\psi$. Let $\{N_\alpha\}_{\alpha=1}^{n+1-m}$ be a local orthonormal frame in the normal bundle of the immersion. Taking traces in  Gauss equation with respect to the induced metric $g$, we deduce that the local components of the intrinsic Ricci tensor in $M$ are given by
\begin{equation}
\label{gauss-ricci}
R_{ij} = \bar R_{ij} - \sum_\alpha \langle \bar R(N_\alpha, \partial_i)\partial_j, N_\alpha\rangle+  \langle II(\partial_i, \partial_j), {\bf H}\rangle- g^{k\ell} \langle II(\partial_k, \partial_j),  II(\partial_i, \partial_\ell)\rangle,
\end{equation}
where $\bar R_{ij}$ is the ambient Ricci tensor and $II$ is the second fundamental form of $\psi$. 
A lenghty but straightforward calculation using (\ref{riemann-warped}) yields
\begin{eqnarray}
\label{gauss-ricci-2}
& & R_{ij} =  (m-2)\bigg(\frac{h'^2}{h^2}-\frac{h''}{h}-\frac{\kappa}{h^2}\bigg)\theta_i\theta_j- \bigg(\frac{h''}{h}|\partial_t^\top|^2+(m-1-|\partial_t^\top|^2)\bigg(\frac{h'^2}{h^2}-\frac{\kappa}{h^2}\bigg)\bigg) g_{ij} \nonumber\\
& & \,\,\,\,+  \langle II(\partial_i, \partial_j), {\bf H}\rangle- g^{k\ell} \langle II(\partial_k, \partial_j),  II(\partial_i, \partial_\ell)\rangle,
\end{eqnarray}
where $\theta_i$ are the local components of the form ${\rm d}t$ dual to $\partial_t$. 

On the other hand the tensorial soliton equation (\ref{solitonB-2}) for $T =X^\top = h(t) \partial_t^\top = \nabla \eta$ may be written in terms of local components as
\[
\langle II(\partial_i, \partial_j), -{\bf H}\rangle + c\langle\nabla_{\partial_i} \nabla \eta, \partial_j\rangle= c\,h'(\pi\circ \psi) g_{ij},
\]
where we used the identity
\[
\frac{1}{2}\pounds_{T} g = \frac{1}{2}\pounds_{\nabla\eta} g = \nabla\nabla\eta.
\]
Substituting and rearranging some of the terms in (\ref{gauss-ricci-2}), we obtain 
\begin{eqnarray*}
& & R_{ij} -c \nabla_{\partial_i}\nabla_{\partial_j}\eta = -ch' g_{ij}- (m-1)\bigg(\frac{h'^2}{h^2}-\frac{\kappa}{h^2}\bigg)+|\partial_t^\top|^2\bigg(\frac{h'^2}{h^2}-\frac{h''}{h}-\frac{\kappa}{h^2}\bigg)g_{ij}\\
& & \,\,\,\, + (m-2)\bigg(\frac{h'^2}{h^2}-\frac{h''}{h}-\frac{\kappa}{h^2}\bigg)\theta_i\theta_j- g^{k\ell} \langle II(\partial_k, \partial_j),  II(\partial_i, \partial_\ell)\rangle
\end{eqnarray*}
Now
\[
g^{k\ell} \langle II(\partial_k, \partial_j),  II(\partial_i, \partial_\ell)\rangle \le |II|^2 g_{ij}
\]
and a simple algebraic  argument  yields
\begin{eqnarray*}
R_{ij} -c \nabla_{\partial_i}\nabla_{\partial_j}\eta \ge (-ch'- |II|^2+ (m-1)\varkappa ) g_{ij} 
\end{eqnarray*}
where
\begin{equation}
\varkappa = \min\bigg\{-\frac{h''}{h}, \frac{\kappa}{h^2}-\frac{h'^2}{h^2}\bigg\}\cdot
\end{equation}
Using (\ref{cond-wmp-1})  we conclude  that
\begin{equation}
\label{cond-wmp-2}
{\rm Ric}_M -c \nabla\nabla \eta \ge -(m-1)\Lambda g.
\end{equation}
By Proposition 8.6 of \cite{AMR} there exists a geodesic ball $B_{R_0} = B_{R_0}(o)\subset M$ centered at $o\in M$ with sufficiently small radius $R_0>0$ and a constant $C=C(B_{R_0})>0$ such that
\begin{equation}
\label{delta-vol}
\Delta_{-c\eta} r(x) \le C+ (m-1)\Lambda_+ r(x)\,\, \mbox{ on } \,\, M\backslash B_{R_0},
\end{equation}
where $r(x) = {\rm dist}_M (x,o)$. Setting 
\[
{\rm vol}_{-c\eta} (B_r) = \int_{B_r} e^{c\eta}\, dM
\] 
to denote the weighted volume of the geodesic ball $B_r$, using  Proposition 8.11 of \cite{AMR} we deduce 
\[
{\rm vol}_{-c\eta} (B_r) \le \int_{0}^r e^{C\tau + (m-1)\Lambda_+ \frac{\tau^2}{2}}\, d\tau+D
\]
for some constant $D>0$ and $C$ as above. Hence, since 
\[
\int_0^r e^{C\tau + (m-1)\Lambda_+ \frac{\tau^2}{2}}\, d\tau + D \sim \frac{e^{Cr+ (m-1)\Lambda_+ \frac{r^2}{2}}}{(m-1)\Lambda_+r+C}\,\, \mbox{ as }\,\, r\to+\infty
\]
we deduce that
\[
\liminf_{r\to \infty} \frac{\log{\rm vol}_{-c\eta} (B_r)}{r^2}<+\infty.
\]
Observing that
\[
\Delta_{-c\eta} = e^{-c\eta}\, {\rm div}(e^{c\eta}\nabla\cdot),
\]
we apply Theorem 4.1 of \cite{AMR} with the choices $T=g$, $\varphi(x,s) = se^{c\eta(x)}$, $A(x) = e^{c\eta(x)}$, $b(x)=1$ to deduce the validity of the weak maximum principle for the operator $\Delta_{-c\eta}$ on $M$. \hfill $\square$

\vspace{3mm}

The following is an immediate consequence of Theorem \ref{cond-wmp}.

\begin{corollary}
\label{cor-cond-wmp} Let $\psi: M\to \bar M = I\times_h P$ be a complete mean curvature flow soliton with respect to $X= h(t)\partial_t$ contained in a slab $[a,b]\times P$.  Suppose that $P$ has constant sectional curvature $\kappa$ and that $\sup_M |II|^2 <\infty$. Then the weak maximum principle holds for the operator $\Delta_{-c\eta}$ on $M$. 
\end{corollary}

Furthermore, it is not difficult to see that $\bar M = I\times_h P$ has constant sectional curvature $\bar\kappa$ if and only if $P$ has constant sectional curvature $\kappa$ and $h$ is a solution of the following differential equation
\[
-\frac{h''(t)}{h(t)}=\frac{\kappa}{h^2(t)}-\frac{h'^2(t)}{h^2(t)}=\bar\kappa.
\]
Therefore, as another consequence of Theorem 4 we have the following
\begin{corollary}
\label{cor-cond-wmp-bis}
Let $\psi: M\to \bar M = I\times_h P$ be a complete mean curvature flow soliton with respect to $X= h(t)\partial_t$.  Suppose that $\bar M$ has constant sectional curvature and that $\sup_M (|II|^2+ch'(\pi\circ\psi)) <+\infty$. Then the weak maximum principle holds for the operator $\Delta_{-c\eta}$ on $M$. 
\end{corollary}

The validity of the weak maximum principle in Theorem \ref{cond-wmp} has been detected via the growth of the volume of geodesic balls. We are now going to give an upper estimate which is particularly efficient in case of warped product targets that are space forms. 

First we need the following result that shows, together with Remark \ref{ricci-mcf}, a further similarity with Ricci solitons.

\begin{proposition}
\label{ode-eta}
Let $\psi: M^m \to \bar M^{n+1} = I\times_h P$ be a mean curvature flow soliton with respect to $X= h(t)\partial_t$. Suppose that for some $\chi\in C^1(\mathbb{R})$, $\chi\ge 0$ on $I$, the function $\hat\eta$ defined in {\rm (\ref{inth})} satisfies the differential equation
\begin{equation}
\label{hateta-2}
\hat\eta' = \chi(\hat\eta)^{1/2} \,\, \mbox{ on }\,\, I.
\end{equation}
Then, there exists a constant $C\in \mathbb{R}$ such that
\begin{equation}
\label{eta-ode}
|{\bf H}|^2 +c^2|\nabla \eta|^2 - c^2\chi(\eta) = C \,\, \mbox{ on }\,\, M,
\end{equation}
where $\eta(x) = \hat \eta((\pi\circ\psi)(x)), \, x\in M$.
\end{proposition}

\noindent \emph{Proof.} Since $\bar M$ is a warped product $I\times_h P$ and $X = h(t)\partial_t$ equations (\ref{solitonB-2}) and (\ref{solitonD-2})  can be respectively written in the form
\begin{equation}
\label{soliton-gradB}
II_{-{\bf H}} +c\, \nabla\nabla \eta = c\,\hat\eta''(\pi\circ\psi) g 
\end{equation}
and
\begin{equation}
\label{soliton-gradC}
\frac{1}{2} \nabla |{\bf H}|^2 = c\,II_{-{\bf H}} (\nabla\eta, \cdot)^\sharp.
\end{equation}
Having fixed a vector field $U$ on $M$ we then have
\begin{eqnarray*}
c^2\,\hat\eta''(\pi\circ\psi) \langle \nabla \eta, U\rangle = c\,II_{-{\bf H}} (\nabla \eta, U) +c^2\, \langle\nabla_U\nabla \eta, \nabla \eta\rangle= \frac{1}{2}\langle \nabla |{\bf H}|^2 + c^2\nabla |\nabla \eta|^2, U\rangle.
\end{eqnarray*}
On the other hand, using (\ref{hateta-2}) we have
\begin{eqnarray*}
\langle \nabla \chi(\eta), U\rangle =\chi'(\eta) \langle \nabla \eta, U\rangle = 2\hat\eta''\langle \nabla\eta, U\rangle.
\end{eqnarray*}
Hence, 
\begin{eqnarray*}
\nabla |{\bf H}|^2 + c^2\nabla |\nabla \eta|^2 -c^2\nabla\chi(\eta)\equiv 0 \,\, \mbox{ on } \,\, M.
\end{eqnarray*}
This finishes the proof of the proposition. \hfill $\square$

\vspace{3mm}

For instance let $\bar M^{n+1} = \mathbb{R}\times_{e^t}\mathbb{R}^n = \mathbb{H}^{n+1}$. Then with the choice $\chi(\tau) = \tau^2$ equation (\ref{hateta-2}) is satisfied. Similarly if we have $\bar M = \mathbb{R}^+ \times_{\sinh t} \mathbb{S}^n =\mathbb{H}^{n+1}$ the choice $\chi(\tau) = \tau^2-1$ works. For $\bar M^{n+1} = \mathbb{R}\times_t \mathbb{S}^n= \mathbb{R}^{n+1}$ we choose $\chi(\tau)= 2\tau$ and finally for $\bar M^{n+1} = (0,\pi)\times_{\sin t}\mathbb{S}^n=\mathbb{S}^{n+1}$ we choose $\chi(\tau) = 2\tau-\tau^2$. Observe that in this case $\hat \eta(t) = 1-\cos t$.

We observe that equation (\ref{hateta-2}) only involves the warping function $h$ of $\bar M = I\times_h P$  and therefore the previous examples work whatever factor $P$ we choose. For instance we can consider $\bar M = \mathbb{R}\times_{e^t} P$ with the choice $\chi(\tau) = \tau^2$. In this case $\bar M$ is not necessarily of constant sectional curvature. 

\begin{remark}
It is worth to note that we can add to $\chi$ any constant so that  in {\rm (\ref{eta-ode})} we can always assume to have normalized our choice in such a way that 
\begin{equation}
\label{ode-eta-2}
|{\bf H}|^2 +c^2|\nabla \eta|^2 - c^2\chi(\eta) = 0 \,\, \mbox{ on }\,\, M.
\end{equation}
\end{remark}

\vspace{3mm}

\begin{proposition}
\label{int-eta}
Let $\psi:M^m \to \bar M^{n+1} = I\times_h P$ be a complete mean curvature flow soliton with respect to $X= h(t)\partial_t$. Assume the validity of {\rm (\ref{hateta-2})} in Proposition \ref{ode-eta}. Then, there exists a constant $C\ge 0$ such that
\begin{equation}
\label{eta-est-1}
|\nabla \eta(x)| \le \int_{0}^{r(x)} |h'(\pi\circ\psi \circ\gamma(\tau))|\, {\rm d}\tau + C,
\end{equation}
where  $r(x) ={\rm dist}_M (o,x)$ and  $\gamma:[0, r(x)]\to M$ is any minimizing unit speed geodesic connecting a fixed point $o\in M$ to $x\in M$. In particular, if $h'$ is bounded 
\begin{equation}
\label{eta-est-2}
|\nabla\eta(x)|\le B\,r (x) + C
\end{equation}
for some constants $B, C\ge 0$. 
\end{proposition}

\noindent \emph{Proof.} As observed above we can assume the validity of (\ref{ode-eta-2}). From here we deduce
\begin{equation}
\label{eta-est-3-bis}
|\nabla \eta|^2 \le \chi(\eta) \,\, \mbox{ on }\,\, M.
\end{equation} 
Having fixed the origin $o\in M$, we let $\gamma$ be a unit speed minimizing geodesic from $o$ to $x\in M$ parametrized on $[0, r(x)]$, where $r(x) = {\rm dist}_M(o,x)$. Define
\[
k(\tau) = \eta(\gamma(\tau)), \quad \tau\in [0, r(x)].
\]
Denoting derivatives with respect to $\tau$ by $\cdot$ we have
\[
\dot k(\tau) = \langle \nabla \eta(\gamma(\tau)), \dot\gamma(\tau)\rangle = h((\pi\circ\psi)(\gamma(\tau)))\langle \partial_t, \dot\gamma(\tau)\rangle
\]
which yields
\begin{equation}
\label{k-est}
|\dot k(\tau)| \le |\nabla \eta(\gamma(\tau))| \le \sqrt{\chi(\eta(\gamma(\tau)))},
\end{equation}
so that
\[
|\dot k(\tau)| \le \sqrt{\chi(k(\tau))}.
\]
Hence, 
\begin{equation}
\label{k-der}
\frac{d}{d\tau}\sqrt{\chi(k(\tau))} \le \left| \frac{\chi'(k(\tau))}{2\sqrt{\chi(k(\tau))}} \dot k (\tau)\right|\le \frac{1}{2}|\chi'(k(\tau))|.
\end{equation}
From (\ref{hateta-2}) since $h>0$ we have
\begin{equation}
\label{chi-der}
\chi'(k(\tau)) = 2\hat\eta''(\pi\circ\psi\circ\gamma(\tau)) = 2 h'(\pi\circ\psi\circ\gamma(\tau)).
\end{equation}
Integrating (\ref{k-der}) on $[0,r(x)]$ and using (\ref{chi-der}) together with (\ref{k-est}) we obtain
\[
|\nabla \eta(x)| \le \int_{0}^{r(x)} |h'(\pi\circ\psi \circ\gamma(\tau))|\, {\rm d}\tau + C
\]
with $C = \sqrt{\chi(\eta(o))}\ge 0$. We have thus proved (\ref{eta-est-1}). \hfill $\square$

\vspace{3mm}

We now give a lower bound on $|\nabla \eta|$ at least in case $\bar M = I\times_h P$ with  $P$ of constant sectional curvature.

\begin{proposition}
\label{lower-bound}
Let $\psi: M^m\to \bar M^{n+1} = I\times_h P$ be a complete mean curvature flow soliton  with respect to $X= h(t)\partial_t$ and $c<0$. Assume  that $P$ has constant sectional curvature $\kappa$. Furthermore suppose that 
\begin{equation}
\label{ric-est}
{\rm Ric}_M \le \Lambda g
\end{equation}
for some constant $\Lambda$ and 
\begin{equation}
\sup_M (|II|^2-(m-1)\varkappa) <+\infty,
\end{equation}
where $\kappa$ is defined in (\ref{varkappa}). Then 
\begin{equation}
\label{eta-est-3}
|\nabla\eta(x)|\ge \frac{1}{c}\,\sup_M (|II|^2-(m-1)\varkappa)\,r(x) +\int_0^{r(x)} h'(\pi\circ\psi\circ\gamma(\tau) ) \,{\rm d}\tau + A
\end{equation}
for some constant $A$.
\end{proposition}

\noindent \emph{Proof.}  Again let $o\in M$ be a fixed origin and $\gamma:[0, r(x)]\to M$ a unit speed minimizing geodesic with $\gamma(0)=o$ and $\gamma(r(x))= x$, where $r(x) = {\rm dist}_M (o,x)$. As above, for $\tau\in [0, r(x)]$, we define 
\[
k(\tau) = \eta(\gamma(\tau)).
\]
Denoting by $\cdot$ derivatives with respect to $\tau$, we recall that
\[
\dot k(\tau) = \langle \nabla \eta(\gamma(\tau)), \dot\gamma(\tau)\rangle.
\]
Since $\gamma$ is geodesic and using (\ref{soliton-gradB}) we obtain
\begin{equation}
\label{k-der-2}
\ddot k(\tau)  = \langle \nabla_{\dot\gamma(\tau)}\nabla \eta, \dot\gamma(\tau)\rangle = h'(\pi\circ\psi\circ\gamma(\tau)) +\frac{1}{c}\, II_{{\bf H}} (\dot\gamma(\tau), \dot\gamma(\tau)).
\end{equation}
From Gauss equation (\ref{gauss-ricci-2}) we have
\begin{equation}
\label{gauss-ricci-3}
{\rm Ric}_M 
(\dot\gamma, \dot\gamma) \ge (m-1)\varkappa+  II_{{\bf H}}(\dot\gamma, \dot\gamma)- |II|^2.
\end{equation}
On the other hand, under the assumption (\ref{ric-est}), by the second variation formula for  arclenght, see equation (8.163) in \cite{AMR}, we obtain
\begin{equation}
\label{int-ric}
\int_0^{r(x)} {\rm Ric}_M(\dot\gamma, \dot\gamma) \le 2(m-1)+ 2\Lambda.
\end{equation}
Next, we integrate (\ref{k-der-2}) on $[0,r(x)]$ and we use (\ref{gauss-ricci-3}) and (\ref{int-ric}) to deduce
\begin{eqnarray*}
& &\dot k(r(x)) - \dot k(0)  =   \int_0^{r(x)} \Big(h'(\pi\circ\psi\circ\gamma (\tau)) +\frac{1}{c} II_{{\bf H}} (\dot\gamma (\tau), \dot\gamma (\tau))\Big)\, {\rm d}\tau\\
& & \,\,\ge   \int_0^{r(x)} \Big(h'(\pi\circ\psi\circ\gamma ) -\frac{1}{c} (m-1)\varkappa(\pi\circ\psi\circ\gamma(\tau)) +\frac{1}{c} |II|^2+\frac{1}{c}\,{\rm Ric}_M 
(\dot\gamma, \dot\gamma)\Big)\\
&  & \,\, \ge  \frac{1}{c}\sup_M (|II|^2-(m-1)\varkappa) \,r(x) +\frac{1}{c} (2(m-1)+2\Lambda)+\int_0^{r(x)} \big(h'(\pi\circ\psi\circ\gamma(\tau) ) \,{\rm d}\tau.
\end{eqnarray*}
From (\ref{k-est}) we then obtain 
\[
|\nabla\eta(x)|\ge \frac{1}{c}\sup_M (|II|^2-(m-1)\varkappa) \,r(x) +\frac{1}{c} (2(m-1)+2\Lambda)+\int_0^{r(x)} h'(\pi\circ\psi\circ\gamma(\tau) ) \,{\rm d}\tau + \dot k(0),
\]
that is, (\ref{eta-est-3}).\hfill $\square$

\begin{proposition}
\label{vol-quad} 
In the assumptions of Proposition \ref{lower-bound} let {\rm (\ref{hateta-2})} be satisfied. Furthermore, assume
\begin{equation}
\label{lim-eta}
\lim_{x\to\infty} \eta(x) = +\infty
\end{equation}
and
\begin{equation}
\label{ii-est}
\sup_M ( |II|^2-(m-1)\varkappa) < -c\inf_M  h'(\pi\circ\psi)\le -c\sup_M h'(\pi\circ\psi)<+\infty.
\end{equation}
Then, there exist constants $\mu, C>0$ such that
\begin{equation}
\label{vol-est-eta}
{\rm vol}_{c\eta} (B_R(o)) \le C R^\mu\,\, \mbox{ for } \,\, R\gg 1.
\end{equation}
\end{proposition}

\noindent \emph{Proof.}  The argument is by now standard (see \cite{CZ}); we report it here for the sake of completeness. Define
\begin{equation}
\varrho(x) = 2 \sqrt{-c\eta(x)}
\end{equation}
so that
\begin{equation}
\label{der-rho}
\nabla \varrho = \sqrt{-c}\,\frac{\nabla \eta}{\sqrt\eta} \quad \mbox{ and } \quad \frac{\nabla\varrho}{|\nabla\varrho|} = \frac{\nabla\eta}{|\nabla \eta|}\cdot
\end{equation}
For $R\gg 1$ we set
\[
D_R =\{x\in M: \varrho(x) <R\}.
\]
Completeness of $M$  and assumption (\ref{lim-eta}) imply that $D_R$ is relatively compact in $M$. We can therefore define
\[
V(R) = {\rm vol}_{c\eta} (D_R) = \int_{D_R} e^{c\eta}\, {\rm d}M.
\]
By the co-area formula we deduce
\[
V(R) =\int_0^R \int_{\partial D_R} \frac{e^{c\eta}}{|\nabla \varrho|} 
\]
and
\[
V'(R) = \int_{\partial D_R} \frac{e^{c\eta}}{|\nabla \varrho|} = \frac{R}{2(-c)}\int_{\partial D_R} \frac{e^{c\eta}}{|\nabla \eta|}\cdot
\]
Next, by equation (\ref{delta-eta}) and the condition $c<0$ we have
\[
\Delta \eta = mh'(\pi\circ\psi) +\frac{1}{c}|{\bf H}|^2 \le  mh'(\pi\circ\psi)
\]
so that 
\[
{\rm div} (e^{c\eta}\nabla\eta) = e^{c\eta} (\Delta\eta+c|\nabla\eta|^2) \le e^{c\eta} (mh'(\pi\circ\psi)+c h^2(\pi\circ\psi)).
\]
Integrating on $D_R$ we obtain
\[
\int_{D_R} (mh'(\pi\circ\psi) +c h^2(\pi\circ\psi))\, e^{c\eta}\, {\rm d}M \ge \int_{D_R} {\rm div} (e^{c\eta} \nabla \eta)\, {\rm d}M = \int_{\partial D_R}
\langle \nabla\eta, \nu\rangle e^{c\eta},
\]
where $\nu$ is the outward unit normal to $\partial D_R$. From (\ref{der-rho})
\[
\nu = \frac{\nabla \varrho}{|\nabla \varrho|} = \frac{\nabla \eta}{|\nabla \eta|},
\]
and having set
\[
C_0 = \sup_M (mh'(\pi\circ\psi) +c h^2(\pi\circ\psi))
\]
from the above we obtain  
\begin{equation}
\label{vol-est-R}
C_0\, V(R)  \ge \int_{\partial D_R} |\nabla \eta| e^{c\eta}, \quad R\gg 1.
\end{equation}
Next, we observe that from the first inequality in (\ref{ii-est}) and (\ref{eta-est-3}), we have
\[
|\nabla\eta(x)|^2 \ge C_1\, r(x)^2
\]
while from (\ref{eta-est-2}) integrating 
\begin{equation}
\label{eta-est-4}
\eta(x) \le C_2\, r(x)^2
\end{equation}
for some constants $C_1, C_2>0$ and $r(x)\gg 1$. It follows that, for the constant $C_3 = C_1/C_2>0$  
\begin{equation}
\label{eta-est-5}
|\nabla \eta|^2 \ge C_3\, \eta \,\, \mbox{ on }\,\, D_R
\end{equation}
for $R\gg 1$ because of (\ref{lim-eta}). On the other hand, using (\ref{eta-est-5}) we have
\begin{eqnarray*}
& & \int_{\partial D_R} |\nabla\eta|e^{c\eta}= \frac{R}{2(-c)}\int_{\partial D_R} \frac{2(-c)}{R} |\nabla \eta|^2\frac{e^{c\eta}}{|\nabla \eta|} \ge
C_3\frac{R}{2(-c)}\int_{\partial D_R}\frac{2(-c)}{R}\eta \frac{e^{c\eta}}{|\nabla \eta|} \\
& &\,\,\,\,= C_3 \frac{R^2}{4(-c)}\int_{\partial D_R}\frac{e^{c\eta}}{|\nabla \eta|}=C_3 \frac{R}{2} V'(R)
\end{eqnarray*}
and therefore, from (\ref{vol-est-R}) we arrive at the differential inequality
\[
2\frac{C_0}{C_3} \,V(R) \ge R V'(R)\quad \mbox{for}\quad R\gg 1.
\]
As a consequence, for some constant $C>0$, 
\begin{equation}
\label{vol-est-R-2}
V(R)\le  C R^{2C_0/C_3}\quad \mbox{for}\quad R\gg 1.
\end{equation}
On the other hand, from (\ref{eta-est-4}) we have
\[
B_R(o) \subset D_{2\sqrt{C_2} R} \quad \mbox{for}\quad R\gg 1.
\]
and (\ref{vol-est-eta}) follows at once. \hfill $\square$

\section{First applications of the weak maximum principle}
\label{sec7}
Let $\bar M^{n+1}= I\times_h P$ be a warped product space as in Example \ref{warped-example}. Let $\phi: M^m\to P^n$ be  an immersion of a $m$-dimensional connected manifold $M^m$, with $m<n$, into the Riemannian manifold $(P^{n},g_0)$.  Let us denote by
$g_M$ the Riemannian metric induced on $M$ via $\phi$, that is,
$g_M=\phi^{*}(g_0)$.

For a fixed $\bar{t}\in I$, let $\phi_{\bar{t}}:M^m\to\bar M^{n+1}= I\times_h P$ be the map given by
\[
\phi_{\bar{t}}(x)=\iota_{\bar t}(\phi(x)), \quad \text{ for every } x\in M,
\]
where $\iota: P \to I\times P$ is defined by $\iota_{\bar t} (y) = (\bar t, y)$ for $y\in P$.
It is not difficult to see that $\phi_{\bar{t}}$ is an immersion of $M$ into $\bar M$ with codimension $n-m+1\geq 2$, which is contained in the leaf $P_{\bar{t}}=\{\bar{t}\}\times P$, and that the metric induced on $M$
via $\phi_{\bar{t}}$ from the warped metric 
\[
\bar g=\textrm{d}t^2+h^2(t)g_0(x)
\]
is simply
\begin{equation}
\label{Luis1}
g_{\bar{t}}=h^2(\bar{t})g_M.
\end{equation}
Conversely, given an immersion $\psi:M^m\to\bar M^{n+1}= I\times_h P$ with codimension $n-m+1\ge 2$ such that the submanifold $\psi(M)$  is contained in a leaf $P_{\bar{t}}$, the projection $\phi=\pi_P\circ\psi:M^m\to P^n$ is an immersed submanifold such that $\psi(x)=\iota_{\bar{t}}(\phi(x))=\phi_{\bar{t}}(x)$. 

It follows from (\ref{Luis1}) that, intrinsically, $(M,g_{\bar{t}})$ is homothetic to $(M,g_M)$ with scale factor $h(\bar{t})$. Our objective now is to express the extrinsic geometry of the submanifold $\phi_{\bar{t}}(M)\subset  \bar M= I\times_h P$ in terms of the extrinsic geometry of  $\phi(M)\subset P$. 

In order to compute the second fundamental form $II_{\bar{t}}$ of the submanifold $\phi_{\bar{t}}(M)$, let $\{\xi_\alpha\}_{\alpha=1}^{n-m}$ be a (locally defined) orthonormal frame for the normal bundle of  $\phi$. It follows that the vector fields
\begin{equation}
\label{Luis4}
N_\alpha(\phi_{\bar t}(x))=\frac{1}{h(\bar t)}\iota_{\bar t*}\xi_\alpha (\phi(x)), \quad 1\leq\alpha\leq n-m
\end{equation}
and 
\[
N_{n-m+1}(\phi_{\bar t}(x))=\partial_t|_{\phi_{\bar t}(x)}
\]
defines a local orthonormal frame in the normal bundle of $\phi_{\bar t}$.  The Weingarten map $A_\alpha$ of $\phi_{\bar t}$ in the direction of $N_\alpha$ is given by
\begin{equation}
\label{Luis2}
A_{\alpha}U=\frac{1}{h(\bar t)}A^P_{\alpha}U, \quad U\in\Gamma(TM).
\end{equation}
where $A^{P}_{\alpha}$ stands for the Weingarten operator of  $\phi:M^m\to P^n$ with respect to the normal direction $\xi_\alpha$.
On the other hand, it follows from (\ref{warped-conformal}) that
\begin{equation}
\label{Luis3}
A_{n-m+1}U = -\bar\nabla_U \partial_t=-\frac{h'(t)}{h(t)}(U+\langle U,\partial_t\rangle\partial_t)=-\frac{h'(t)}{h(t)}U
\end{equation}
for every vector field $U\in\Gamma(TM)$. 
%
From (\ref{Luis2}) and (\ref{Luis3}), the second fundamental $II_{\bar t}$ of the submanifold $\phi_{\bar t}$ can be written, for tangent vector fields $U,V\in\Gamma(TM)$, as
\begin{equation}
\label{Luis5}
II_{\bar t}(U,V)=\sum_{\alpha=1}^{n-m+1}\langle A_\alpha U,V\rangle N_\alpha=
\frac{1}{h(\bar t)}\sum_{\alpha=1}^{n-m}\langle A^P_\alpha U,V\rangle N_\alpha
-\frac{h'(\bar t)}{h(\bar t)}\langle U,V\rangle\partial_t.
\end{equation}
Observe that, from (\ref{Luis1}) and (\ref{Luis4}), one has
\[
\langle A^P_\alpha U,V\rangle N_\alpha=h(\bar t)g_0(A^P_\alpha U,V)\xi_\alpha
\]
for every $1\leq \alpha\leq n-m$, so that (\ref{Luis5}) becomes
\begin{equation}
\label{Luis5-bis}
II_{\bar t}(U,V)=\sum_{\alpha=1}^{n-m}g_0(A^P_\alpha U,V)\xi_\alpha
-\frac{h'(\bar t)}{h(\bar t)}\langle U,V\rangle\partial_t=
II_\phi(U,V)-\frac{h'(\bar t)}{h(\bar t)}\langle U,V\rangle\partial_t,
\end{equation}
where $II_\phi$ stands for the second fundamental form of the immersion $\phi:M^m\to P^n$.  Taking traces in both sides with respect to the 
metric $g_{\bar t}$ we deduce that the (non-normalized) mean curvature vector field ${\bf H}_{\bar t}$ of  $\phi_{\bar t}$ is
\begin{equation}
\label{Luis6}
{\bf H}_{\bar t}=\frac{1}{h^2(\bar t)}{\bf H}_\phi-m\frac{h'(\bar t)}{h(\bar t)}\partial_t.
\end{equation}
where ${\bf H}_\phi={\rm tr}_{g_M} II_{\phi}$ is the (non-normalized) mean curvature vector field of  $\phi(M)\subset P$.

On the other hand, since $\partial_t$ is normal to $\phi_{\bar t}(M)$, we also have
\begin{equation}
\label{Luis7}
X=X^\perp=h(\bar t)\partial_t
\end{equation}
along $\phi_{\bar t}$. It then follows from (\ref{Luis6}) and (\ref{Luis7}) that 
$\phi_{\bar t}:M^m\to\bar M^{n+1}=I\times_hP$ is a mean curvature flow soliton with respect to $X=h(t)\partial_t$ if and only if $\phi:M^m\to P^n$ is a minimal submanifold (that is, ${\bf H}_\phi\equiv 0$) and $\bar t\in\mathbb{R}$ is given implicity by
\begin{equation}
\label{Luis8}
mh'(\bar t)+ch^2(\bar t)=0
\end{equation}
or, equivalently, by 
\begin{equation}
\label{Luis9}
nh'(\bar t)+\frac{nc}{m}h^2(\bar t)=0
\end{equation}
Recall from Example \ref{warped-example-2} that (\ref{Luis9}) geometrically means that the leaf $P_{\bar t}$ is a codimension one mean curvature flow soliton of $\bar M$ with respect to $X$ with $\hat{c}=nc/m$. 

We can summarize all the above in the following:
\begin{proposition}
\label{Luisprop1}
Let $\psi:M^m\to \bar M^{n+1}= I\times_h P$ be a submanifold with $m<n$ and assume that $\psi(M)$ is contained in a leaf $P_{\bar t}=\{{\bar t}\}\times P$. Then the submanifold $\psi(M)$ is a mean curvature flow soliton with respect to the vector field $X=h(t)\partial_t$ if and only if the projection $\phi=\pi_P\circ\psi:M^m\to P^n$ is a minimal submanifold and ${\bar t}\in I$ is implictly given by
\begin{equation}
\label{Luis10}
mh'(\bar t)+ch^2(\bar t)=0,
\end{equation}
which means that the leaf $P_{\bar t}$ is a mean curvature flow soliton of $\bar M$ with respect to $X$ with $\hat{c}=nc/m$. 
\end{proposition}

Since the function 
\begin{equation}
\zeta(t) = mh'(t) +ch^2(t), \quad t\in I,
\end{equation}
will repeatedly appear in the sequel we reserve for it the name of \textit{soliton function}. The terminology is justified by what follows. 

Now we are ready to give the first main result in this section, which is a direct consequence of Proposition \ref{laplace-eta} and Definition \ref{parabolic}.
\begin{theorem}
\label{sol-function}
Let $\psi:M^m\to \bar M^{n+1}= I\times_h P$ be a mean curvature flow soliton with respect to the vector field $X=h(t)\partial_t$. Consider the soliton function 
\[
\zeta(t) = mh'(t) +ch^2(t), \quad t\in I,
\]
and let  $\zeta(\pi\circ\psi): M \to \mathbb{R}$ be non-negative {\rm(}resp., non-positive{\rm)} on $M$. Assume that $\pi\circ\psi$ {\rm(}equivalently, $\eta = \hat \eta(\pi\circ\psi)${\rm)} is bounded above {\rm(}resp. bounded below{\rm)} on $M$. If $M$ is parabolic with respect to the operator $\Delta_{-c\eta}$ then $\psi(M)$ is contained in the leaf $P_{\bar t}=\{\bar t\}\times P$ given implicitly by
\begin{equation}
\label{Luis11}
mh'(\bar t)+ch^2(\bar t)=0.
\end{equation} 
\end{theorem}
For the proof of Theorem \ref{sol-function} simply observe that equation (\ref{delta-eta-X}) can be written as $\Delta_{-c\eta}\eta=\zeta(\pi\circ\psi)$. Then, the $\Delta_{-c\eta}$-parabolicity of $M$ implies that $\eta$ is constant, so that $\pi\circ\psi$ is constant and $\psi(M)$ is contained in a leaf $P_{\bar t}$. Condition (\ref{Luis11}) follows from the fact that $\Delta_{-c\eta}\eta=\zeta({\bar t})=0$ (see also Proposition \ref{Luisprop1}). 

\

Motivated by the above considerations, we introduce the following terminology.
\begin{definition}\label{def-hs} 
Suppose that there exists $\bar t \in I$ such that 
\begin{equation}
\zeta(\bar t) = mh'(\bar t) + ch^2(\bar t)=0.
\end{equation}
We say that the immersion $\psi:M^m \to \bar M^{n+1} = I\times_h P$ is on one side of the leaf $P_{\bar t}$ if the function 
\[
\zeta(\pi\circ\psi)= mh'(\pi\circ\psi)+ch^2(\pi\circ\psi)
\]
has constant sign on $M$.
\end{definition}
As a consequence of Theorem \ref{sol-function} we have the next
\begin{corollary}
\label{half-space-B} Let $\psi: M^m\to \bar M^{n+1} =I\times_h P$ be a complete mean curvature flow soliton with respect to $X=h(t)\partial_t$  and such that 
\begin{equation}
\label{par-sol}
\frac{1}{{\rm vol} (\partial B_r)}\notin L^1(+\infty),
\end{equation}
where $\partial B_r$ is the boundary of the geodesic ball $B_r$ of radius $r$ centered at a fixed origin $o\in M$.
Assume 
\begin{equation}
\label{sol-slab}
\psi(M) \subset [a,b]\times P, 
\end{equation}
for some interval $[a,b]\subset I$ and that $\psi(M)$ is on one side of $P_{\bar t}$. Then $\psi(M)\subset P_{\bar t}$. 
\end{corollary}

\noindent \emph{Proof.} First of all observe that (\ref{sol-slab}) together with completeness and (\ref{par-sol}) imply that $M$ is $\Delta_{-c\eta}$-parabolic. Actually, since $\eta=\hat\eta(\pi\circ\psi)$ is bounded on $M$, we have $\sup_M e^{c\eta}=C<+\infty$, regardless the sign of $c$, and 
\[
{\rm vol}_{-c\eta}(\partial B_r)=\int_{\partial B_r} e^{c\eta}\, {\rm d}M\le C{\rm vol}(\partial B_r).
\]
Hence
\[
\frac{1}{{\rm vol}_{-c\eta} (\partial B_r)}\notin L^1(+\infty)
\]
and therefore $M$ is $\Delta_{-c\eta}$-parabolic (see Remark \ref{Luis-remark-parabolicity}). Now the result follows directly from Theorem \ref{sol-function}. \hfill $\square$

\begin{remark}
In case $M$ is compact the conclusion of Corollary \ref{half-space-B} holds obviously without requiring assumptions {\rm (\ref{par-sol})} and {\rm (\ref{sol-slab})}.
\end{remark}

%

As another consequence of Proposition \ref{laplace-eta} we have the next mean curvature estimates that will give rise both to rigidity results and height estimates.
\begin{theorem}
\label{half-space-A-Luis} Let $\psi: M^m\to \bar M^{n+1}= I\times_h P$ be a stochastically complete mean curvature flow soliton with respect to $X= h(t) \partial_t$ on $\bar M$.
\begin{enumerate}
\item[{\rm i)}] Assume 
\begin{equation}
\label{halfspace-Luis}
t_* = \inf_M \,(\pi\circ \psi)>-\infty, \quad t_*\in I.
\end{equation}
Then the mean curvature vector field of $M$ satisfies
\begin{equation}
\label{Luis12}
\inf_M|{\bf H}|^2\le -mch'(t_*) \quad \mbox{ if } c<0
\end{equation}
and 
\begin{equation}
\label{Luis13}
\sup_M|{\bf H}|^2\ge -mch'(t_*) \quad \mbox{ if } c>0.
\end{equation}
\item[{\rm ii)}] Assume
\begin{equation}
\label{halfspace}
t^* = \sup_M \,(\pi\circ \psi) <+\infty, \quad t^*\in I.
\end{equation}
Then the mean curvature vector field of $M$ satisfies
\begin{equation}
\label{Luis14}
\sup_M|{\bf H}|^2\ge -mch'(t^*) \quad \mbox{ if } c<0
\end{equation}
and 
\begin{equation}
\label{Luis15}
\inf_M|{\bf H}|^2\le -mch'(t^*) \quad \mbox{ if } c>0.
\end{equation}
\end{enumerate}
\end{theorem}
\noindent \emph{Proof.} Since $\hat{\eta}'=h>0$, the function
$\eta=\hat{\eta}(\pi\circ\psi)$ can be thought as a reparametrization of $\pi\circ\psi$. Therefore, under assumption (\ref{halfspace-Luis}) it holds that
\[
\eta_*=\inf_M\eta=\hat{\eta}(t_*)>-\infty,
\]
and, by the stochastic completeness of $M$, we may apply the weak maximum principle for $\Delta$ to the function $\eta$. Hence, there exists a sequence of points $\{x_k\}_{k\in\mathbb{N}}$ in $M$ such that
\[
\eta(x_k)<\eta_*+\frac{1}{k} \quad \mbox{ and } \quad \Delta\eta(x_k)>-\frac{1}{k}.
\]
Let $t_k=(\pi\circ\psi)(x_k)$ and observe that $\lim_{k\to+\infty}t_k=t_*$ because $\hat{\eta}$ is strictly increasing. Using (\ref{delta-eta}) we have 
\[
\Delta\eta(x_k)=mh'(t_k)+\frac{1}{c}|{\bf H}|^2(x_k)>-\frac{1}{k},
\]
and thus
\[
-\frac{1}{k}-mh'(t_k)<\frac{1}{c}|{\bf H}|^2(x_k)\le \sup_M\left(\frac{1}{c}|{\bf H}|^2\right).
\]
Letting $k\to+\infty$ we get
\[
-mh'(t_*)\le \sup_M\left(\frac{1}{c}|{\bf H}|^2\right)=
\begin{cases}
\frac{1}{c}\inf_M|{\bf H}|^2 & \mbox{ if } c<0,\\
{} & {} \\
\frac{1}{c}\sup_M|{\bf H}|^2 & \mbox{ if } c>0.
\end{cases}
\]
In other words,
\[
\inf_M|{\bf H}|^2\le -mch'(t_*) \quad \mbox{ if } c<0
\]
and 
\[
\sup_M|{\bf H}|^2\ge -mch'(t_*) \quad \mbox{ if } c>0.
\]
Similarly, under assumption (\ref{halfspace}) and working now at $t^*$ we get
\[
-mh'(t^*)\ge \inf_M\left(\frac{1}{c}|{\bf H}|^2\right)=
\begin{cases}
\frac{1}{c}\sup_M|{\bf H}|^2 & \mbox{ if } c<0,\\
{} & {} \\
\frac{1}{c}\inf_M|{\bf H}|^2 & \mbox{ if } c>0,
\end{cases}
\]
which yields
\[
\sup_M|{\bf H}|^2\ge -mch'(t^*) \quad \mbox{ if } c<0
\]
and 
\[
\inf_M|{\bf H}|^2\le -mch'(t^*) \quad \mbox{ if } c>0.
\]
\hfill $\square$

We now analyze some rigidity consequences of Theorem \ref{half-space-A-Luis}.
For the compact case we get the following corollary that will be used in Corollary \ref{Luis-compact-2}.
\begin{corollary}
\label{Luis-compact} Let $\psi: M^m\to \bar M^{n+1}= I\times_h P$ be a compact mean curvature flow soliton with respect to $X= h(t) \partial_t$ on $\bar M$. If $c<0$, the mean curvature vector field of $M$ satisfies
\begin{equation}
\min_M|{\bf H}|^2\le -mch'(t_*) \quad \text{and} \quad \max_M|{\bf H}|^2\ge -mch'(t^*),
\end{equation}
where $t_*=\min_M(\pi\circ\psi)$ and $t^*=\max_M(\pi\circ\psi)$. Similarly, if $c>0$ the mean curvature vector field of $M$ satisfies
\begin{equation}
\min_M|{\bf H}|^2\le -mch'(t^*) \quad \mbox{ and } \quad \max_M|{\bf H}|^2\ge -mch'(t_*).
\end{equation}
In particular, if
\[
h''(t)\ge 0\,\, \mbox{ for } t\in I, \,\, t_*\le t\le t^*, 
\]
then
\begin{enumerate}
\item[{\rm i)}] If $c<0$ then $\min_M|{\bf H}|^2\le -mch'(t_*)\le -mch'(t^*)\le \max_M|{\bf H}|^2$.
\item[{\rm ii)}] If $c>0$ then $\min_M|{\bf H}|^2\le -mch'(t^*)\le -mch'(t_*)\le \max_M|{\bf H}|^2$.
\end{enumerate}
\end{corollary}
Obviously, the final statement of Corollary \ref{Luis-compact} follows directly from the fact that $h'(t_*)\le h'(t^*)$ under the additional assumption $h''\ge 0$. As a consequence of the above result we get
\begin{corollary}
\label{Luis-compact-2} Let $\psi: M^m\to \bar M^{n+1}= I\times_h P$ be a compact mean curvature flow soliton with respect to $X= h(t) \partial_t$ on $\bar M$ and assume that 
\[
h''(t)\ge 0\,\, \mbox{ for } t\in I, \,\, t_*\le t\le t^*
\]
with $t_*$ and $t^*$ as in Corollary \ref{Luis-compact}. 
If $|{\bf H}|$ is constant then $\psi(M)\subset P_{\bar t}=\{\bar t\}\times P$, where $\bar t\in I$ is given implicitly by 
\begin{equation}
mh'(\bar t)+ch^2(\bar t)=0.
\end{equation}
Moreover, if $m=n$ then $P$ is necessarily compact and $\psi(M)=P_{\bar t}$. If $m<n$ then $\psi(M)=\{\bar t\}\times M_0$ with $M_0\subset P$ a compact minimal submanifold of $P$.
\end{corollary}
\noindent \emph{Proof.} Since $h''\ge 0$ on $[t_*,t^*]$ and $|{\bf H}|$ is constant, we know from i) and ii) in Corollary \ref{Luis-compact} that
\[
h'(t_*)=h'(t^*)=\frac{-1}{mc}|{\bf H}|^2=\text{ constant}.
\]
Since $h'(t)$ is non-decreasing on $[t_*,t^*]$, from here it follows  that 
\[
h'(t)=\frac{-1}{mc}|{\bf H}|^2=\text{ constant on } [t_*,t^*]
\]
and therefore
\[
mh'(\pi\circ\psi)+\frac{1}{c}|{\bf H}|^2=0 \text{ on } M.
\]
In other words 
\[
\Delta\eta =mh'(\pi\circ\psi)+\frac{1}{c}|{\bf H}|^2=0
\]
on the compact manifold $M$, which implies that $\eta$ and, equivalently, $\pi\circ\psi$, is constant on $M$; that is $\psi(M)$ is contained in a leaf $P_{\bar t}$. The rest of the proof follows directly from Proposition \ref{Luisprop1}. \hfill $\square$

The next result extends Corollary \ref{Luis-compact-2} to the case of stochastically complete mean curvature flow solitons, with the help of Theorem \ref{half-space-A-Luis} 
\begin{corollary}
Let $\psi: M^m\to \bar M^{n+1}= I\times_h P$ be a stochastically complete mean curvature flow soliton with respect to $X= h(t) \partial_t$ on $\bar M$. Assume that 
\[
\psi(M)\subset[a,b]\times P
\]
for some interval $[a,b]\subset I$ on which 
\[
h''(t)\ge 0.
\]
If $|{\bf H}|$ is constant and $h'(t)$ is not locally constant on $[a,b]$ {\rm(}in other words, the equality $h''=0$ holds only at isolated points{\rm)} then $\psi(M)\subset P_{\bar t}=\{\bar t\}\times P$, where $\bar t\in I$ is given implicitly by 
\begin{equation}
mh'(\bar t)+ch^2(\bar t)=0.
\end{equation}
Moreover, if $m=n$ then $P$ is necessarily stochastically complete and $\psi(M)=P_{\bar t}$. If $m<n$ then $\psi(M)=\{\bar t\}\times M_0$ with $M_0\subset P$ a stochastically complete minimal submanifold of $P$.
\end{corollary}
\noindent \emph{Proof.} Let $t_*=\inf_M(\pi\circ\psi)\ge a>-\infty$ and $t^*=\sup_M(\pi\circ\psi)\le b<+\infty$. Since $h''\ge 0$ on $[t_*,t^*]$, we have $h'(t_*)\le h'(t^*)$ and, since $|{\bf H}|$ is constant, we obtain from Theorem \ref{half-space-A-Luis} similarly to the compact case
\[
h'(t_*)=h'(t^*)=\frac{-1}{mc}|{\bf H}|^2=\text{ constant}.
\]
The hypothesis on $h'(t)$ implies now that $h'(t)$ is strictly increasing on $[t_*,t^*]$. Therefore, $t_*=t^*$ and $\pi\circ\psi$ is constant on $M$; that is $\psi(M)$ is contained in a leaf $P_{\bar t}$. The rest of the proof follows again directly from Proposition \ref{Luisprop1}. \hfill $\square$

\

The general result given in Theorem \ref{half-space-A-Luis} provides some particular height estimates in case of space forms described as warped product spaces that we now illustrate.
\begin{corollary}
\label{corollary-hs-A}
Let $\psi: M^m \to  (0,\infty) \times_t \mathbb{S}^n = \mathbb{R}^{n+1}\backslash\{o\}$ be a stochastically  complete mean curvature flow soliton with respect to $X= t\partial_t$ on $\mathbb{R}^{n+1}\backslash\{o\}$ and $c<0$. Let
\[
t^* = \sup_M \,(\pi\circ\psi).
\]
Then
\begin{equation}
\label{hs-1}
t^*\ge\sqrt{-m/c}.
\end{equation}
\end{corollary}
\noindent \emph{Proof.}  If $t^*=+\infty$ there is nothing to prove. Otherwise, from (\ref{solitonA-2})
\[
|{\bf H}|^2 \le c^2|X\circ \psi|^2 \le  c^2(\pi\circ\psi)^2\le c^2(t^*)^2 \,\, \mbox{ on } \,\, M.
\]
Therefore, since $h(t)=t$ from (\ref{Luis14}) we obtain
\[
-mc\le\sup_M|{\bf H}|^2\le c^2(t^*)^2,
\]
so that $(t^*)^2\ge m/-c$ or, equivalently, (\ref{hs-1}). \hfill $\square$

\vspace{3mm}

Probably Corollary \ref{corollary-hs-A} can be better stated in the more geometrical form
\begin{corollary}\label{corollary-hs-B}
Let $\psi:M^m \to (0,\infty)\times_t \mathbb{S}^n = \mathbb{R}^{n+1}\backslash\{o\}$ be a stochastically complete mean curvature flow soliton with respect to $X= t\partial_t$ on $\mathbb{R}^{n+1}\backslash\{o\}$ and $c<0$. If $\psi(M)\subset (0,b)\times \mathbb{S}^n$ then $b\ge\sqrt{-m/c}$.
\end{corollary}

The next corollary addresses the case of hyperbolic space foliated by horospheres. 
\begin{corollary}
\label{corollary-hs-A2}
Let $\psi: M^m \to  \mathbb{R} \times_{e^t} \mathbb{R}^n = \mathbb{H}^{n+1}$ be a stochastically  complete mean curvature flow soliton with respect to $X= e^t\partial_t$ on $\mathbb{H}^{n+1}$ and $c<0$. Let
\[
t^* = \sup_M \,(\pi\circ\psi).
\]
Then
\begin{equation}
\label{hs-3}
t^*\ge\log (-m/c).
\end{equation}
\end{corollary}
Assuming $t^*<+\infty$, otherwise there is nothing to show, the proof is similar to that of Corollary \ref{corollary-hs-A}, observing that, again by the soliton equation, we deduce now that
\[
|{\bf H}|^2 = c^2|(X\circ\psi)^\perp|^2  \le c^2|X\circ \psi|^2 = c^2e^{2 \pi(\psi(x))} \le c^2e^{2t^*} 
\]
on $M$. Since in this case $h(t) = e^t$, from (\ref{Luis14}) we obtain  
\[
-mce^{t^*}\le\sup_M|{\bf H}|^2\le c^2e^{2t^*},
\]
which yields $e^{t^*}\ge -m/c$ or, equivalently, (\ref{hs-3}).

\

Let us now consider a different realization of $\mathbb{H}^{n+1}$ as $(0,\infty)\times_{\sinh t}\mathbb{S}^n= \mathbb{H}^{n+1}\backslash\{o\}$. We have
\begin{corollary}
\label{hyp-sphere}
Let $\psi: M^m \to(0,\infty)\times_{\sinh t}\mathbb{S}^n= \mathbb{H}^{n+1}\backslash\{o\}$ be a stochastically  complete mean curvature flow  soliton with respect to $X=\sinh t\partial_t$ on $\mathbb{H}^{n+1}$ and $c<0$. Let 
\[
t^* = \sup_M \,(\pi\circ\psi).
\]
Then
\begin{equation}
\label{hs-4}
t^*\ge \arg\cosh\left(-\frac{m+\sqrt{m^2+4c^2}}{2c}\right).
\end{equation}
\end{corollary}
Again the proof is similar to the previous ones. In this case, when $t^*<+\infty$ from (\ref{solitonA-2}) we have
\[
|{\bf H}|^2 = c^2|(X\circ\psi)^\perp|^2  \le c^2|X\circ \psi|^2 = c^2\,\sinh^2 ( \pi(\psi(x))) \le c^2\,\sinh^2t^*=c^2\,\cosh^2t^*-c^2
\]
on $M$. On the other hand, using (\ref{Luis14}) we get
\[
-mc\cosh{t^*}\le\sup_M|{\bf H}|^2\le c^2\,\cosh^2t^*-c^2,
\]
which yields 
\[
\cosh^2t^*+\frac{m}{c}\cosh t^*-1\ge 0.
\]
This implies 
\[
\cosh t^*\ge-\frac{m+\sqrt{m^2+4c^2}}{2c}
\]
or, equivalently, (\ref{hs-4}).

\

Another result in this direction is obtained, needless to say, by a further application of the weak maximum principle. However we need to restrict ourselves to the codimension one case. 

\begin{theorem}
\label{hs-D}
Let $\psi:M^m \to \bar M^{m+1} = I\times_h P$ be an orientable, codimension one, mean curvature flow soliton with respect to $X= h(t)\partial_t$. Fix an orientation $N$ of $\psi(M)$ such that  $H = \langle {\bf H}, N\rangle\ge 0$. Suppose that
\begin{equation}
\label{sup-tt}
 t^* = \sup_M \,(\pi\circ\psi)<+\infty
\end{equation}
and let 
\[
\mathcal{H}(t) = m\frac{h'(t)}{h(t)}\cdot
\]
Assume that, for some $\bar t\in I$, 
\begin{equation}
\label{HH}
 \sup_M H < \mathcal{H}(\bar t)  \\
\end{equation}
and
\begin{equation}
\label{supH-dot}
  \mathcal{H}'(t) \ge 0 \,\,\, {\rm  for  } \,\,\, t\ge \bar t.   
\end{equation}
If $M$ is stochastically complete then
\begin{equation}
\label{ineq-tt}
t^*\le \bar t.
\end{equation}
\end{theorem}

\vspace{3mm}

\noindent \emph{Proof.} We reason by contradiction and we suppose that $t^*>\bar t$. We observe that $\pi\circ\psi$ cannot be constant on $M$; otherwise $\pi(\psi(x)) = t^*>\bar t$ and
\[
\psi(M)\subseteq \{t^*\}\times P.
\]
If this is the case, because of (\ref{supH-dot}) and our choice of $t^*$ the scalar mean curvature of the immersion would satisfy
\[
H = \mathcal{H}(t^*) \ge \mathcal{H}(\bar t)>0,
\] 
contradicting (\ref{HH}). Since $\pi\circ\psi$ is non-constant on $M$ we can choose a regular value $\gamma$ of $\pi\circ\psi$ with 
\[
\bar t < \gamma < t^*
\]
and the property that $\partial\Omega_\gamma \neq \emptyset$ where 
\[
\Omega_\gamma =\{x\in M: (\pi\circ\psi)(x)>\gamma\}.
\]
Next we define on $M$ the angle function
\[
\Theta = \langle \partial_t, N\rangle = \frac{1}{h(\pi\circ\psi)}\langle X, N\rangle.
\]
This enables us to transform equation (\ref{delta-eta}) of Proposition \ref{laplace-eta} as follows
\begin{eqnarray*}
\Delta\eta & = & mh'(\pi\circ\psi) +\frac{1}{c} |{\bf H}|^2 = 
h(\pi\circ\psi) \Big(\mathcal{H}(\pi\circ\psi)+\frac{1}{ch(\pi\circ\psi)}\langle {\bf H}, N\rangle^2\Big)\\
{} & = &  h(\pi\circ\psi) (\mathcal{H}(\pi\circ\psi)+\langle N, \partial_t\rangle\langle {\bf H}, N\rangle)\\
{} & = &  h(\pi\circ\psi)(\mathcal{H}(\pi\circ\psi)+\Theta H).
\end{eqnarray*}
Because of the assumptions on $\mathcal{H}$ and $\mathcal{H}'$ we infer 
\[
\mathcal{H}(t)\ge\mathcal{H}(\bar t)>\sup_MH\ge 0 \text{ for } t\ge\bar t. 
\]
In particular $h$ is increasing for $t\ge\bar t$ and
\[
h(\pi\circ\psi)> h(\gamma) \,\, \mbox{ on } \,\, \Omega_\gamma.
\]
Furthermore, since $\mathcal{H}$ is non-decreasing for $t\ge \bar t$ we also have
\[
\mathcal{H}(\pi\circ\psi)\ge \mathcal{H}(\gamma)\ge \mathcal{H}(\bar t)>\sup_MH\ge H \,\, \mbox{ on} \,\, \Omega_\gamma.
\]
From $H=\langle {\bf H},N\rangle \ge 0$ and $\Theta \ge -1$ we deduce $\Theta H \ge - H$ and
\[
\mathcal{H}(\pi\circ\psi)+ \Theta H \ge \mathcal{H}(\pi\circ\psi) - H\ge \mathcal{H}(\gamma)- \sup_M H>0
\]
on $\Omega_\gamma$, so that
\[
\Delta\eta  \ge h(\gamma) (\mathcal{H}(\gamma)- \sup_M H)>0
\]
on $\Omega_\gamma$. Since $\eta=\hat\eta(\pi\circ\psi)$ and $\hat\eta$ is increasing 
\[
\Lambda_{\hat\eta(\gamma)} = \{x\in M: \eta(x)>\hat\eta(\gamma)\} = \Omega_\gamma
\]
and
\[
\partial\Lambda_{\hat\eta(\gamma)} = \partial\Omega_\gamma.
\]
We set $\Omega = \Lambda_{\hat\eta(\gamma)}$ and 
\[
v = \eta|_{\bar\Omega}
\]
to deduce, from the above,
\begin{equation}
\begin{cases}
\Delta v\ge h(\gamma) (\mathcal{H}(\gamma)- \sup_M H)>0\,\, \mbox { on } \,\, \Omega, & \\
\sup_\Omega v = \hat\eta(t^*)<+\infty. & 
\end{cases}
\end{equation}
Stochastic completeness of $M$ enables us to apply the open form of the weak maximum principle to infer
\[
\hat\eta(t^*) = \sup_\Omega v = \sup_{\partial\Omega} v = \hat\eta(\gamma).
\]
Since $\eta(t^*)>\hat\eta(\gamma)$ we obtain the desired contradiction completing  the proof of Theorem \ref{hs-D}. \hfill $\square$

\begin{remark}
Since every parabolic manifold is stochastically complete, Theorem \ref{hs-D} remains valid for parabolic manifolds. 
\end{remark}

Let us consider the following special case of the theorem to illustrate the result.
\begin{corollary}
Let $\psi: M \to \mathbb{H}^{m+1} = \mathbb{R}\times_{e^t} \mathbb{R}^m$ be an orientable, codimension one,  mean curvature flow soliton with respect to $X= e^t\partial_t$. Fix an orientation $N$ of $\psi(M)$ such that  $H = \langle {\bf H}, N\rangle\ge 0$
and assume that
\begin{align}
\label{HH3}
& \sup_M H<m.  
\end{align}
If $M$ is stochastically complete, then $t^* = \sup_M (\pi\circ\psi)=+\infty$.
\end{corollary}
Clearly this result can be interpreted as a ''half-space'' theorem similar to that of Hoffman and Meeks \cite{HMe} for minimal surfaces in $\mathbb{R}^3$.

\vspace{3mm}

Next results are obtained with a different technique. We begin with
\begin{theorem}
\label{thm-zeta}
Let $\psi: M^m \to \bar M^{n+1} = I\times_h P$ be a complete mean curvature flow soliton with respect to $X= h(t)\partial_t$ with $c<0$ and such that ${\rm vol} (M) = +\infty$. Assume that, having fixed some origin $o\in M$,
\begin{equation}
\label{logvol}
\limsup_{r\to +\infty} \frac{\log {\rm vol} \,B_r}{r}=\alpha\in\mathbb{R}^+,
\end{equation}
where $B_r\subset M$ is the geodesic ball of radius $r$ centered at $o$.
For $[a,b]\subset I$ let
\[
\zeta_* = \inf_{[a,b]}\zeta,
\]
where $\zeta(t) = mh'(t)+ch^2(t)$ is the soliton function. Furthermore suppose  that, for some compact $K\subset M$
\begin{equation}
\label{slab-t}
\psi(M\backslash K) \subseteq [a,b]\times P
\end{equation}
and let $\hat\eta$ be the function defined in {\rm (\ref{inth})}. Then 
\begin{equation}
\label{zeta-eta}
\zeta_* \le \frac{\alpha^2}{4}\big(\hat\eta(b)-\hat\eta(a)\big).
\end{equation}
\end{theorem}

\noindent \emph{Proof.} Since $\hat\eta$ is increasing
\[
\Lambda = \hat\eta(b)-\hat\eta(a) = \int_a^b h(s)\, {\rm d}s\ge 0,
\]
hence (\ref{zeta-eta}) is certainly satisfied if $\zeta_*\le 0$. Thus without loss of generality we can suppose $\zeta_*> 0$ so that $\zeta(t)> 0$ for $t\in [a,b]$. Furthermore, we can also assume that $\hat\eta$ is defined by
\[
\hat\eta(t) = \int^t_a h(s) \, {\rm d}s.
\]
Recall from the soliton equation (\ref{solitonA-2}) and from (\ref{Luis16}) that, along the immersion $\psi$, 
\[
X=h(t)\partial_t=X^\top+X^\perp=\nabla\eta+\frac{1}{c}{\bf H}
\]
with $\eta= \hat \eta (\pi\circ\psi)$, so that
\[
\frac{1}{c}|{\bf H}|^2=ch^2(\pi\circ\psi)-c|\nabla\eta|^2.
\]
Since $c<0$, from equation (\ref{delta-eta}) we then know that $\eta$ satisfies
\[
\Delta \eta = mh'( \pi\circ \psi) +\frac{1}{c} |{\bf H}|^2\ge mh'(\pi\circ\psi) +ch^2(\pi\circ\psi)
\]
on $M$. Fix $R>0$ sufficiently large such that for some fixed origin $o\in M$, we have $K\subset B_R$ and let $\varepsilon>0$. Then from  (\ref{slab-t})  we deduce
\[
\eta(x)=\hat\eta((\pi\circ\psi)(x)) \le \int_a^b h(s)\, {\rm d}s < \Lambda+\varepsilon
\]
on $M\backslash B_R$. We define
\[
v = \Lambda+\varepsilon -\eta>0
\]
on $M\backslash B_R$. By  a simple generalization of Barta's theorem, \cite{Barta}, we deduce
\[
\lambda_1^\Delta(M\backslash B_R) \ge \inf_{M\backslash B_R} \Big(-\frac{\Delta v}{v}\Big) = \inf_{M\backslash B_R}\frac{\Delta \eta}{v} \ge \inf_{M\backslash B_R} \frac{1}{v}\big(mh'(\pi\circ\psi)+ch^2(\pi\circ\psi)\big)\ge \frac{\zeta_*}{\Lambda+\varepsilon},
\]
where $\lambda_1^\Delta(M\backslash B_R)$ is the first eigenvalue of the Laplace-Beltrami operator $\Delta$ on $M\backslash B_R$.
Letting $\varepsilon \searrow 0^+$ we finally obtain
\[
\lambda_1^\Delta (M\backslash B_R) \ge \frac{\zeta_*}{\Lambda}\cdot
\]
On the other hand, using  (\ref{logvol}) and the assumption  ${\rm vol}(M)=+\infty$,  from Theorem 3.23 and Persson formula (2.88) of  \cite{BMR},  we obtain 
\[
\lambda_1^\Delta (M\backslash B_R) \le \frac{\alpha^2}{4}\cdot
\]
Comparing the last two inequalities we obtain the validity of (\ref{zeta-eta}). \hfill $\square$

\

More interestingly we have the following consequence.
\begin{corollary}
\label{corollary-hs-C1}
Let $\psi: M^m \to \bar M^{n+1} = I\times_h P$ be a complete mean curvature flow soliton with respect to $X= h(t)\partial_t$ with $c<0$, subexponential volume growth, infinite volume and such that $\psi(M\backslash K)\subset [a,b]\times P$, for some compact $K\subset M$ and  $[a,b]\subset I$. Suppose that the soliton function $\zeta(t)$ satisfies $\zeta(\hat{t})\ge 0$ for some $\hat{t}\in [a,b]$. Then there exists $\bar t\in [a,b]$ such that the  corresponding leaf  $P_{\bar t}= \{\bar t\}\times P$ is a mean curvature flow soliton of $\bar M$ with respect to $X$ with $\hat c=nc/m$. 
\end{corollary}

\noindent \emph{Proof.}  We observe that, since $M$ has subexponential volume growth,  (\ref{logvol}) holds for each $\alpha>0$. Hence letting $\alpha \searrow 0^+$ we obtain $\zeta_*\le 0$. Since $\zeta(\hat{t})\ge 0$ for some $\hat{t}\in [a,b]$, by continuity there exists $\bar t\in[a,b]$ such that $\zeta(\bar t)=0$, which means that the  corresponding leaf  $P_{\bar t}= \{\bar t\}\times P$ is a mean curvature flow soliton of $\bar M$ with respect to $X$ with $\hat c=nc/m$. \hfill $\square$

\

We describe here some other consequences of Theorem \ref{thm-zeta}. 
For instance, we apply it with $\bar M = \mathbb{R}\times_{e^t}\mathbb{R}^n = \mathbb{H}^{n+1}$. Then $\zeta(t) = e^t(m+ce^t)$ and $\zeta(t)\ge 0$ on $(-\infty, \log(-m/c)]$. However the only zero of $\zeta(t)$ is at $\log(-m/c)$; it follows that there are no complete mean curvature flow solitons with respect to $X=e^t\partial_t$ with $c<0$, subexponential volume growth, infinite volume and such that $\psi(M\backslash K) \subset [a,b]\times \mathbb{R}^n$ for some compact $K\subset M$ and $b<\log(-m/c)$.

Reasoning in a similar way and representing $\mathbb{H}^{n+1}\backslash\{o\}$ as $\mathbb{R}^+\times_{\sinh t}\mathbb{S}^n$, we deduce that there are no complete mean curvature flow solitons with respect to  the vector field $X=\sinh t\partial_t$ with $c<0$, subexponential volume growth, infinite volume and such that $\psi(M\backslash K) \subset [a,b]\times \mathbb{S}^n$ for some compact subset $K\subset M$ and $\cosh b<-(m+\sqrt{m^2+4c^2})/2c$.

\

Finally we explicitly state the next result, again a further direct consequence of Theorem \ref{thm-zeta}.
\begin{corollary}
\label{corollary-hs-C2} 
Let $\psi: M^m\to \bar M^{n+1}= (0,+\infty)\times_t P$ be a complete mean curvature flow soliton with respect to $X = t\partial_t$ and $c<0$. Assume ${\rm vol} \,M = +\infty$ and
\[
\limsup_{r\to +\infty} \frac{\log {\rm vol} \,B_r}{r}= \alpha \in \mathbb{R}^+
\]
Moreover, suppose that for some compact $K\subset M$
\[
\psi(M\backslash K) \subset [a,b]\times P,
\]
for some $[a,b]\subset (0,+\infty)$. Then
\begin{equation}
\label{ineq-b}
b\ge \sqrt{\frac{8m+\alpha^2 a^2}{\alpha^2-8c}}\cdot
\end{equation}
\end{corollary}

\vspace{3mm}

\noindent \emph{Proof. }  Note that in this case $\zeta(t) = m+ct^2$ so that
\[
\zeta_* = \inf_{[a,b]}\zeta(t) = m+cb^2.
\]
Furthermore, 
\[
\hat\eta(t) =\int_a^t s\, {\rm d}s = \frac{1}{2}(t^2-a^2).
\] 
Therefore (\ref{zeta-eta}) of Theorem \ref{thm-zeta} yields
\[
m+cb^2 \le \frac{\alpha^2}{8} \big(b^2-a^2),
\]
from which (\ref{ineq-b}) follows at once. \hfill $\square$

\

Note that Corollary \ref{corollary-hs-C2} applies in particular to $\bar M=\mathbb{R}^{n+1}\backslash\{o\}=(0,+\infty)\times_{t}\mathbb{S}^n$. In case 
$\psi: M^m\to \mathbb{R}^{n+1}\backslash\{o\}=(0,+\infty)\times_{t}\mathbb{S}^n$ is proper and $c<0$, estimate (\ref{ineq-b}) becomes $b\ge \sqrt{-m/c}$. Indeed, because of properness the volume growth of the geodesic balls of $M$ is at most polynomial.

\section{A further fundamental equation and some characterization results}
\label{char}

The aim of this section is to provide a useful expression for the weighted Laplacian of
the squared norm of ${\bf H}$, and infer some interesting geometric consequences.
\begin{proposition}
\label{laplacian-H2}
Let $\psi: M^m \to \bar M^{n+1}$ be a mean curvature flow soliton  with respect to a closed conformal vector field $X\in \Gamma(T\bar M)$ on $\bar M$. Let ${\bf H}$ be its mean curvature vector field. Then 
\begin{equation}
\label{deltaH2} \frac{1}{2}\Delta_{-cT} |{\bf H}|^2= -c\varphi |{\bf
H}|^2 -|II_{\bf H}|^2+ |(\bar\nabla_{(\cdot)}{\bf H})^\perp|^2-mc
\langle\bar\nabla\bar\varphi, {\bf H}\rangle - {\rm tr}_M \, \bar R
(\cdot, {\bf H}){\bf H},
\end{equation}
where $T$ is the tangential component of $X$ along the immersion $\psi$, 
\begin{equation}
\bar\varphi = \frac{1}{n+1}\,{\rm div}_{\bar M } X,
\end{equation}
and $\varphi=\bar\varphi\circ\psi$.
\end{proposition}

\noindent \emph{Proof.} From (\ref{solitonD-2}) we have
\begin{eqnarray}
\label{grad-H2-bis}  \frac{1}{2}\langle \nabla |{\bf H}|^2, U\rangle
 =-c\,II_{\bf
H}(T, U), \,\,\, \mbox{ for all }\,\, U\in \Gamma(TM).
\end{eqnarray}
Differentiating (\ref{grad-H2-bis}), we deduce
\begin{eqnarray*}
& & \frac{1}{2}\langle \nabla_U\nabla |{\bf H}|^2, V\rangle  =
\frac{1}{2}U\langle \nabla |{\bf H}|^2, V\rangle -
\frac{1}{2}\langle \nabla |{\bf H}|^2,
\nabla_U V\rangle  \\
& &\,\,=-c\,U ( II_{\bf H}(T,V)) +c \, II_{\bf H}(T, \nabla_U V) =-c\,\nabla_U II_{\bf H} (T,V)-c\, II_{\bf H}(\nabla_U T, V).
\end{eqnarray*}
On the other hand Codazzi's equation may be written in the form
\begin{eqnarray*}
& &\nabla_U II_{\bf H} (T,V) - \nabla_T II_{\bf H}(U,V) =
-\langle \bar R(\psi_* U,\psi_*T){\bf H},\psi_*V\rangle \\
& &\,\,+\langle\bar\nabla_{\psi_* U}{\bf H}, II(T,V)\rangle
-\langle\bar\nabla_{\psi_* T}{\bf H}, II(U,V)\rangle.
\end{eqnarray*}
We note that the soliton equation $cX^\perp ={\bf H}$ implies that
\begin{eqnarray*}
& & \langle\bar\nabla_{\psi_* U}{\bf H}, II(T,V)\rangle
=c\langle\bar\nabla_{\psi_* U}(X-\psi_* T), II(T,V)\rangle\\
& & \,\, =c\varphi \langle \psi_* U,
II(T,V)\rangle-c\langle\bar\nabla_{\psi_* U}\psi_* T,
II(T,V)\rangle=-c\,\langle II(T,U), II(T,V)\rangle.
\end{eqnarray*}
Therefore,
\begin{eqnarray*}
& & \frac{1}{2}\langle \nabla_U\nabla |{\bf H}|^2, V\rangle
=-c\nabla_T II_{\bf H} (U,V)+c\langle \bar R (\psi_* U,
\psi_*T){\bf H}, \psi_* V)\\
& & \,\,\,\, +c^2\langle II(T,U), II(T,V)\rangle+c\langle\bar\nabla_{\psi_* T}{\bf H}, II(U,V)\rangle-c
II_{\bf H}(\nabla_U T, V).
\end{eqnarray*}
Taking traces and using (\ref{solitonA-2}) again we get
\begin{eqnarray*}
& & \frac{1}{2}\Delta |{\bf H}|^2  
= -c\langle\nabla |{\bf H}|^2, T\rangle+cg^{ij}\langle \bar R
(\psi_* \partial_i, X){\bf H}, \psi_* \partial_j) -g^{ij}\langle
\bar R (\psi_* \partial_i,{\bf H}){\bf H}, \psi_* \partial_j)\\
& &\,\,\,\, +c^2|II(T,\cdot)|^2+\frac{c}{2}\langle \nabla |{\bf
H}|^2, T\rangle -c g^{ij}(II_{{\bf H}})_{j}^k
\langle\bar\nabla_{\psi_*\partial_i} \psi_*T, \psi_*\partial_k
\rangle
\\
& &\,\, = -\frac{c}{2}\langle\nabla |{\bf H}|^2,
T\rangle+cg^{ij}\langle \bar R (\psi_* \partial_i, X){\bf H},
\psi_* \partial_j) -g^{ij}\langle \bar R (\psi_* \partial_i,{\bf
H}){\bf H}, \psi_*
\partial_j)
\\
& & \,\,+c^2|II(T,\cdot)|^2 -c g^{ij}(II_{{\bf H}})_{j}^k
\langle\bar\nabla_{\psi_*\partial_i} X, \psi_*\partial_k \rangle +
g^{ij}(II_{{\bf H}})_{j}^k \langle\bar\nabla_{\psi_*\partial_i}
{\bf H}, \psi_*\partial_k \rangle
\\
& &\,\, = -\frac{c}{2}\langle\nabla |{\bf H}|^2,
T\rangle+cg^{ij}\langle \bar R (\psi_* \partial_i, X){\bf H},
\psi_* \partial_j) -g^{ij}\langle \bar R (\psi_* \partial_i,{\bf
H}){\bf H}, \psi_*
\partial_j)
\\
& & \,\, +c^2|II(T,\cdot)|^2  -c \varphi\, g^{ij}(II_{{\bf H}})_{j}^k
g_{ik} - g^{ij}(II_{{\bf H}})_{j}^k
\langle\bar\nabla_{\psi_*\partial_i} \psi_*\partial_k, {\bf H}
\rangle,
\end{eqnarray*}
and we conclude that
\begin{eqnarray*}
& & \frac{1}{2}\Delta |{\bf H}|^2   = -\frac{c}{2}\langle\nabla
|{\bf H}|^2, T\rangle+cg^{ij}\langle \bar R (\psi_* \partial_i,
X){\bf H}, \psi_* \partial_j) -g^{ij}\langle \bar R (\psi_*
\partial_i,{\bf H}){\bf H}, \psi_* \partial_j)
\\
& & \,\, +c^2|II(T,\cdot)|^2 -c\varphi |{\bf H}|^2 -|II_{\bf H}|^2.
\end{eqnarray*}
Now, using (\ref{grad-H2-bis}) we obtain
\begin{eqnarray*}
& &  \frac{1}{2}\Delta_{-cT} |{\bf H}|^2  =  \frac{1}{2}\Delta |{\bf H}|^2  +\frac{c}{2}\langle \nabla |{\bf H}|^2, T\rangle \\
&  & \,\,\,\,=-c\varphi |{\bf H}|^2 + |(\bar\nabla_{(\cdot)}{\bf
H})^\perp|^2-|II_{\bf H}|^2 +cg^{ij}\langle \bar R (\psi_*
\partial_i, X){\bf H}, \psi_* \partial_j) -g^{ij}\langle \bar R
(\psi_* \partial_i,{\bf H}){\bf H}, \psi_* \partial_j).
\end{eqnarray*}
On the other hand (\ref{closedconf}) yields
\begin{eqnarray*}
\bar R(\psi_* \partial_j, {\bf H}) X =\bar\nabla_{\psi_*\partial_j} \bar\varphi {\bf H} - \bar\nabla_{\bf H} \bar\varphi \psi_*\partial_j - \bar\varphi [\psi_*\partial_j, {\bf H}] = \langle \bar\nabla \bar\varphi, \psi_*\partial_j\rangle {\bf H} - \langle\bar\nabla\bar\varphi, {\bf H}\rangle \psi_*\partial_j.
\end{eqnarray*}
Then, we deduce
\begin{eqnarray*}
g^{ij}\langle \bar R (\psi_* \partial_i, X){\bf H}, \psi_*
\partial_j) =g^{ij}\langle \bar R (\psi_* \partial_j, {\bf H})X, \psi_*
\partial_i) = -m \langle \bar\nabla\bar\varphi, {\bf H}\rangle.
\end{eqnarray*}
Finally, note that
\[
g^{ij}\langle \bar R (\psi_* \partial_i,{\bf H}){\bf H}, \psi_*
\partial_j)= {\rm tr}_M \bar R(\cdot, {\bf H}){\bf H}.
\]
This finishes the proof of (\ref{deltaH2}). \hfill $\square$

\

In the particular case where $\bar M^{n+1}=I\times_h P$ is a warped product and $X=h(t)\partial_t$, we know that $T=\nabla\eta$ and $\bar\varphi=h'(t)$. Therefore, $\Delta_{-cT}=\Delta_{-c\eta}$, $\bar\varphi=h'(t)$, $\varphi=h'(\pi\circ\psi)$ and 
\[
\bar\nabla\bar\varphi=h''(t)\partial_t=\frac{h''(t)}{h(t)}X.
\]
In what follows, to simplify the writing, we will denote with $h$, $h'$ and $h''$ the functions $h(t)$, $h'(t)$ and $h''(t)$ along the immersion $\psi$, that is,
\begin{equation}
\label{notation}
h=h(\pi\circ\psi), \quad h'=h'(\pi\circ\psi), \quad h''=h''(\pi\circ\psi).
\end{equation}
Observe that, along the immersion and using the soliton equation (\ref{solitonA-2}), one has
\begin{equation}
\label{luis3}
\langle {\bf H},\partial_t\rangle=\frac{1}{h}\langle {\bf H},X^\perp\rangle=\frac{1}{ch}|{\bf H}|^2.
\end{equation}
This gives
\begin{equation}
\label{luis1}
c\langle\bar\nabla\bar\varphi, {\bf H}\rangle=
\frac{h''}{h}|{\bf H}|^2.
\end{equation}
Moreover, if $P$ has constant sectional curvature $\kappa$, from (\ref{riemann-warped}) we deduce
\begin{eqnarray}
\label{luis2}
{\rm tr}_M \bar R(\cdot, {\bf H}){\bf H} & = & m\Big(\frac{\kappa}{h^2}-\frac{h'^2}{h^2} \Big) |{\bf
H}|^2\\
\nonumber {} & {} & 
-\Big(\frac{\kappa}{h^2}+\frac{1}{h^2}\big(h''h-h'^2\big)\Big)
\big(m\langle {\bf H},
\partial_t\rangle^2 +|\partial_t^\top|^2 |{\bf H}|^2\big).
\end{eqnarray}
Therefore, using (\ref{luis3}), (\ref{luis1}) and (\ref{luis2}), the two last terms in (\ref{deltaH2}) can be easily written in the form
\begin{equation}
\label{luis4}
-mc\langle\bar\nabla\bar\varphi,{\bf H}\rangle-{\rm tr}_M \bar R(\cdot, {\bf H}){\bf H}=-(m-1)\Big(\frac{\kappa}{h^2}+\frac{1}{h^2}\big(h''h-h'^2\big)\Big)|{\bf H}|^2
\Big(1-\frac{|{\bf H}|^2}{c^2h^2}\Big).
\end{equation}
Furthermore, observe that $\partial_t=\partial_t^\top+\partial_t^\perp$
where $\partial_t^\top=\nabla(\pi\circ\psi)$ and, by the soliton equation (\ref{solitonA-2}),
\[
\partial_t^\perp=\frac{1}{h}X^\perp=\frac{1}{ch}{\bf H}.
\]
Hence,
\[
|\nabla(\pi\circ\psi)|^2=|\partial_t^\top|^2=1-|\partial_t^\perp|^2=1-\frac{|{\bf H}|^2}{c^2h^2}\geq 0.
\]
Summing up, we have proved the following
\begin{corollary}
\label{laplacian-H2-bis}
Suppose that $\bar M^{n+1}=I\times_h P$ is a warped product space with $P$ of constant sectional curvature $\kappa$, and let $\psi: M^m \to \bar M^{n+1}$ be a mean curvature flow soliton  with respect to $X=h(t)\partial_t$. Let ${\bf H}$ be its mean curvature vector field. Then
\begin{eqnarray}
\label{deltaH2-kappa} 
\frac{1}{2}\Delta_{-c\eta} |{\bf H}|^2 & = & -c h'
|{\bf H}|^2 -|II_{\bf H}|^2+ |(\bar\nabla_{(\cdot)}{\bf
H})^\perp|^2\\
& {} & -(m-1)\Big(\frac{\kappa}{h^2}+\frac{1}{h^2}\big(h''h-h'^2\big)\Big)|{\bf H}|^2
\Big(1-\frac{|{\bf H}|^2}{c^2h^2}\Big),\nonumber
\end{eqnarray} 
where
\begin{equation}
\label{luis5}
1-\frac{|{\bf H}|^2}{c^2h^2}=|\nabla(\pi\circ\psi)|^2\geq 0.
\end{equation}
\end{corollary}
Recall that $\bar M = I\times_h P$ has constant sectional curvature $\bar\kappa$ if and only if $P$ has constant sectional curvature $\kappa$ and $h$ is a solution of the differential equations
\[
-\frac{h''(t)}{h(t)}=\frac{\kappa}{h^2(t)}-\frac{h'^2(t)}{h^2(t)}=\bar\kappa.
\]
Therefore, as a direct consequence of Corollary \ref{laplacian-H2-bis}, we obtain
\begin{corollary}
\label{laplacian-H2-bis-bis}
Suppose that $\bar M^{n+1}=I\times_h P$ is a warped product space with constant sectional curvature and let $\psi: M^m \to \bar M^{n+1}$ be a mean curvature flow soliton  with respect to $X=h(t)\partial_t$. Let ${\bf H}$ be its mean curvature vector field. Then
\begin{equation}
\label{deltaH2-kappa-bis} 
\frac{1}{2}\Delta_{-c\eta} |{\bf H}|^2=-c h'
|{\bf H}|^2 -|II_{\bf H}|^2+ |(\bar\nabla_{(\cdot)}{\bf
H})^\perp|^2
\end{equation}
\end{corollary}

As a first consequence of Corollary \ref{laplacian-H2-bis} we prove the following
\begin{theorem}
\label{theorem30}
Let $\psi: M^m \to \bar M^{n+1} = I\times_h P$ be a complete mean curvature flow soliton with respect to $X = h(t)\partial_t$. Assume that $P$ has constant sectional curvature $\kappa$ with
\begin{equation}
\label{luis8}
\kappa\le\inf_M(h'^2-h''h)
\end{equation} 
\begin{equation}
\label{kappa-bounded}
\inf_M\left(\frac{\kappa-h'^2}{h^2}\right)>-\infty,
\end{equation}
and
\begin{equation}
\label{hgA}
ch'+|II|^2 \le 0 \,\,\, \mbox{on}\,\,\, M.
\end{equation} 
If $|{\bf H}|\in L^2(M, e^{c\eta}\, {\rm d}M)$ then the mean curvature vector is parallel and either ${\bf H}\equiv 0$ or otherwise $|II|^2\equiv-ch'$ on $M$.
\end{theorem}

\noindent \emph{Proof.}  First of all we observe that (\ref{luis8}) is equivalent to
\[
\frac{\kappa-h'^2}{h^2}\le-\frac{h''}{h} \quad \text{ on } \quad M
\]
and therefore, with the notations of Theorem \ref{cond-wmp}, we have
\[
\varkappa=
\min\bigg\{-\frac{h''}{h}, \frac{\kappa-h'^2}{h^2}\bigg\}
=\frac{\kappa-h'^2}{h^2} \quad \text{ on } \quad M.
\]
Then, (\ref{kappa-bounded}) is equivalent to $\inf_M\varkappa>-\infty$ and, again with the notations of Theorem \ref{cond-wmp}, we have
\begin{equation}
\label{phi-A}
\Lambda + \inf_M\varkappa\le \frac{1}{m-1}\sup_M (ch'+|II|^2)\le 0. 
\end{equation}
Hence, denoting by $r(x)={\rm dist}_M(o,x)$ the distance in $M$ from a fixed origin $o\in M$, we know from (\ref{delta-vol}) of Theorem \ref{cond-wmp}  that
\begin{equation}
\label{lapl-r}
\Delta_{-c\eta} r(x) \le C_0 + (m-1)\Lambda_+ r(x) \quad \text{on} \quad M\setminus B_{R_0}
\end{equation}
for a sufficiently small radius $R_0>0$, where the constant $C_0=C_0(B_{R_0})>0$ is fixed once $R_0$ is so. 

Fix $T$ and $R$ such that $R_0<T<R$ and let $\alpha \in C^1(\mathbb{R}^+_0)\cap C^2([0,R))$ satisfy $\alpha(r)\ge 0$ on $\mathbb{R}^+_0$, 
\[
\alpha(r) \equiv 1 \,\, \mbox{ on }\,\, [0,T], \quad \alpha(r) \equiv 0\,\, \mbox{ on }\,\, [R, +\infty)
\]
and
\[
|\alpha'(r)|\le \frac{C}{R-T} \,\, \mbox{ on }\,\, [0,R], \quad |\alpha''(r)|\le \frac{C}{(R-T)^2}\,\, \mbox{ on } \,\, [0,R],
\]
for some constant $C>0$ independent of $R$ and $T$. Next we define the cut-off function 
\[
f(x) = \alpha(r(x)).
\]
Then we have 
\[
f\ge 0\,\, \text{ on } \,\, M, \quad f(x) \equiv 1 \,\, \mbox{ on }\,\, B_T, \quad f(x) \equiv 0\,\, \mbox{ on }\,\, M\backslash B_R
\]
and
\[
\nabla f =0  \,\, \mbox{ on }\,\, \partial B_R\cup B_T \cup (M\backslash B_R), \quad |\nabla f|\le \frac{C}{R-T}\cdot
\]
We also have
\[
\Delta_{-c\eta} f \equiv 0 \,\, \mbox{ on } \,\, B_T \cup (M\backslash B_R).
\]
Since 
\[
\Delta_{-c\eta} f = \alpha'(r) \Delta_{-c\eta} r + \alpha''(r)\,\,\mbox{ on } \,\, \bar B_R\backslash B_T
\]
using (\ref{lapl-r}) and $\alpha'(r)\equiv 0$ on $[0,T]$ we deduce
\begin{equation}
\label{lapl-phi}
\Delta_{-c\eta} f \le C\frac{R}{R-T}+D\frac{1}{R-T}+ E\frac{1}{(R-T)^2} \quad \text{on} \quad M
\end{equation}
for some constants $C, D, E\ge 0$, independent of $R$ and $T$. 

Let $u\ge 0$ be any $C^2(M)$ function. Then
\begin{eqnarray*}
{\rm div} (f e^{c\eta}\nabla u) & = & f \, {\rm div} (e^{c\eta}\nabla u)+\langle \nabla f, \nabla u\rangle e^{c\eta}\\
{} & = & 
f\,e^{c\eta}\Delta_{-c\eta}u+\langle \nabla f, \nabla u\rangle e^{c\eta}
\end{eqnarray*}
so that, since $f \equiv 0$ on $\partial B_R$ and $\nabla f\equiv 0$ on $B_T$, by the divergence theorem we have
\begin{equation}
\label{div-eta}
\int_{B_R} f\,e^{c\eta}\Delta_{-c\eta}u + \int_{B_R\backslash B_T} \langle \nabla f, \nabla u\rangle e^{c\eta} = 0.
\end{equation}
Now computing 
\begin{equation}
{\rm div} (ue^{c\eta} \nabla f) = u e^{c\eta} \Delta_{-c\eta} f +e^{c\eta}
\langle \nabla u, \nabla f\rangle
\end{equation}
and applying again the divergence theorem, observing that
\[
\nabla f \equiv 0 \,\, \mbox{ on } \,\, \partial B_R \cup \partial B_T,
\]
we deduce 
\[
\int_{B_R\backslash B_T} u e^{c\eta}\Delta_{-c\eta} f  +\int_{B_R\backslash B_T} \langle \nabla u, \nabla f\rangle e^{c\eta}=0.
\]
Inserting into (\ref{div-eta}) we obtain
\begin{equation}
\label{lapl-eta-int}
\int_{B_R} f  e^{c\eta}\Delta_{-c\eta} u  = \int_{B_R\backslash B_T} u e^{c\eta}\Delta_{-c\eta} f .
\end{equation}

We now let $u=|{\bf H}|^2$ and let $T=R/2$. It follows from Corollary \ref{laplacian-H2-bis}, (\ref{luis5}), (\ref{luis8}) and (\ref{hgA}) that
\begin{eqnarray}
\label{luis9}
\nonumber \Delta_{-c\eta}u & = & -ch'u-|II_{\bf H}|^2+
|\bar\nabla^\perp{\bf H}|^2-(m-1)\Big(\frac{\kappa}{h^2}+\frac{1}{h^2}\big(h''h-h'^2\big)\Big)u(1-\frac{u}{c^2h^2})\\
\nonumber & \ge & -ch'u-|II_{\bf H}|^2+ |(\bar\nabla_{(\cdot)}{\bf
H})^\perp|^2 \\
{} & \ge & -(ch'+|II|^2)u+|(\bar\nabla_{(\cdot)}{\bf
H})^\perp|^2\\
\nonumber {} &  \ge & 0 \quad \text{on} \quad M,
\end{eqnarray}
where we have used the fact that $|II_{\bf H}|^2\le|II|^2|{\bf H}|^2$ with equality if and only if $M$ is pseudo-umbilical. In particular, since $f\ge 0$ and $\Delta_{-c\eta}u\ge 0$ on $M$ and $f\equiv 1$ on $B_{R/2}$ we obtain from (\ref{lapl-eta-int}) 
\begin{eqnarray}
\label{luis10}
\nonumber \int_{B_{R/2}} e^{c\eta}\Delta_{-c\eta} u & \le & 
\int_{B_{R/2}} e^{c\eta}\Delta_{-c\eta} u +
\int_{B_R\setminus B_{R/2}} fe^{c\eta}\Delta_{-c\eta} u \\
{} & = & \int_{B_R} fe^{c\eta}\Delta_{-c\eta} u \\
\nonumber {} & = &\int_{B_R\backslash B_{R/2}} u e^{c\eta}\Delta_{-c\eta} f .
\end{eqnarray}
Using (\ref{lapl-phi}) we deduce from here
\begin{equation}
\label{div-eta-R2}
\int_{B_{R/2}} e^{c\eta}\Delta_{-c\eta} u  \le C 
\bigg (1+\frac{1}{R}+\frac{1}{R^2}\bigg)\int_{B_R\backslash B_{R/2}} u e^{c\eta} 
\end{equation}
for some $C>0$ sufficiently large and independent on $R$. Using inequalities  (\ref{luis9}) in (\ref{div-eta-R2}) yields
\[
0\le -\int_{B_{R/2}} (ch'+|II|^2)ue^{c\eta}+
\int_{B_{R/2}}|\bar\nabla^\perp{\bf H}|^2e^{c\eta}\le \frac{C}{2}\bigg (1+\frac{1}{R}+\frac{1}{R^2}\bigg)\int_{B_R\backslash B_{R/2}}ue^{c\eta}.
\]
Since $|{\bf H}|\in L^2(M, e^{c\eta}\, {\rm d}M)$, letting $R\to+\infty$ from the above we infer
$|\bar\nabla^\perp{\bf H}|\equiv 0$ so that $u=|{\bf H}|^2$ is constant and if $u\not\equiv 0$ then 
\[
|II|^2 \equiv -ch',
\] 
completing the proof. 

Observe that in the case where $|{\bf H}|^2$ is a positive constant one gets $\Delta_{-c\eta}u=0$ and all the inequalities in (\ref{luis9}) are equalities.
It follows from here that $M$ is pseudo-umbilical, $ch'+|II|^2\equiv 0$ and, in case inequality (\ref{luis8}) is strict, we also have equality in (\ref{luis5}). This means that $M$ is contained in a slice $P_{\bar t}$ such that $|{\bf H}|^2=c^2h^2({\bar t})$ and $|II|^2=-ch'({\bar t})$.  \hfill $\square$

\vspace{3mm}

Next result in some sense extends work of Cao and Li \cite{CL} and Ding and Xin \cite{DX}. We obtain the same conclusion of Theorem \ref{theorem30} but under different assumptions.
\begin{theorem}
\label{theorem30.1}
Let $\psi: M^m \to \bar M^{n+1} = \mathbb{R}\times_h P$ be a complete mean curvature flow soliton with respect to $X = h(t)\partial_t$. Assume that $P$ has constant sectional curvature $\kappa$ and that {\rm(\ref{luis8})} and {\rm(\ref{kappa-bounded})} are satisfied. Furthermore, let 
\begin{equation}
\label{160.0}
h(t)\notin L^1(+\infty)
\end{equation}
with 
\begin{equation}
\label{160.2}
\lim_{t\to +\infty}\zeta(t)=-\infty,
\end{equation}
where $\zeta(t)=mh'(t) +ch^2(t)$ is the soliton function, and
\begin{equation}
\label{160.1}
\pi\circ\psi(x)\to +\infty \quad \text{as} \quad x\to\infty \,\,\, \text{in} \,\,\, M.
\end{equation}
Set $\Omega_\gamma=\{x\in M : |{\bf H}|^2(x)>\gamma \}$ and suppose that, for some $\gamma\in\mathbb{R}$, $\Omega_\gamma\neq\emptyset$ and
\begin{equation}
\label{160.3}
|II|^2(x)\le\sup_{\Omega_\gamma}(-ch'(\pi\circ\psi))\in\mathbb{R}.
\end{equation} 
Then the mean curvature vector is parallel and either ${\bf H}\equiv 0$ or otherwise $|II|^2\equiv-ch'$ on $M$.
\end{theorem}

\noindent \emph{Proof.}  
As in the proof of Theorem \ref{theorem30}, with the aid of Corollary \ref{laplacian-H2-bis} and (\ref{luis8}) and (\ref{kappa-bounded}), we arrive at the differential inequality 
\begin{equation}
\label{160.4}
\Delta_{-c\eta}u \ge -(ch'+|II|^2)u+|(\bar\nabla_{(\cdot)}{\bf H})^\perp|^2 \quad \text{on} \quad M,
\end{equation}
with $u=|{\bf H}|^2$. Next we observe that because of (\ref{160.1}) and (\ref{160.0})
\[
\eta(x)=\int_{t_0}^{\pi\circ\psi(x)}h(s)ds\to+\infty \quad \text{as} \quad x\to\infty \,\,\, \text{in} \,\,\, M.
\]
Furthermore, from (\ref{160.2}) and equation (\ref{delta-eta-X}) of Proposition \ref{laplace-eta}
\[
\Delta_{-c\eta}\eta\to-\infty \quad \text{as} \quad x\to\infty \,\,\, \text{in} \,\,\, M.
\]
Since $M$ is complete, according to the discussion after Definition \ref{parabolic} the function $\eta$ satisfies the Khas'minskii test for $\Delta_{-c\eta}$-parabolicity. From (\ref{160.4}) we then obtain
\[
\Delta_{-c\eta}u\geq 0 \quad \text{on} \,\,\, \Omega_\gamma.
\]
Now because of (\ref{160.3}) $u$ is bounded above on $\Omega_\gamma$. Using Ahlfors parabolicity we deduce that if $\Omega_\gamma\not\equiv M$, and therefore $\partial\Omega_\gamma\neq\emptyset$,
\[
\sup_{\Omega_\gamma}u=\sup_{\partial\Omega_\gamma}u=\gamma,
\]
which is a contradiction. Hence $\Omega_\gamma=M$ and $u$ is constant on $M$. From (\ref{160.4}) it then follows that ${\bf H}$ is parallel and either ${\bf H}\equiv 0$ or $-ch'(\pi\circ\psi)\equiv |II|^2$ on $M$. \hfill $\square$

\begin{remark}
\label{remark8.1}
Suppose that $\bar M=\mathbb{R}^{n+1}$. Without loss of generality we can suppose that $\psi(M)\subset\mathbb{R}^{n+1}\backslash\{o\}$ so that the latter can be represented as the warped product $(0,+\infty)\times_{t}\mathbb{S}^n$ and $\psi:M\to (0,+\infty)\times_{t}\mathbb{S}^n$. Now we have $X=t\partial_t$ the position vector field. Then $\eta(x)=|\psi|^2(x)$ so that assuming $\psi$ to be a proper immersion {\rm (\ref{160.1})} is satisfied. Clearly, since $h(t)=t$ {\rm(\ref{160.0})} is satisfied too and {\rm (\ref{160.2})} is satisfied if $c<0$. As for the requirement {\rm (\ref{160.3})}, it means
\[
|II|^2\le -c \quad \text{on} \,\,\, \Omega_\gamma.
\]
In this case, for instance using Theorem 4 of {\rm \cite{Lawson}} and Theorem 1 of {\rm \cite{Yau}} the conclusion of Theorem \ref{theorem30.1} yields that either the immersion is totally geodesic or $\psi(M)$ is a product of the type $\mathbb{S}^{k}(\sqrt{-k/c})\times\mathbb{R}^{m-k}$ for some $1\le k\le m$. This result compares directly with Theorem 4 if {\rm\cite{CL}} {\rm(}see also {\rm\cite{DX}}{\rm)}. We observe that Theorem \ref{theorem30.1} compares also with Theorem \ref{theorem44} below.
\end{remark}

For the sake of completeness, we include here the following proof of Remark \ref{remark8.1}, adapted from \cite{CL} to our situation. If ${\bf H}=0$ then $\langle X, N_\alpha\rangle=0$ for all $N_\alpha$ in a local orthonormal frame in the normal bundle. This means that the position vector field is tangent to the submanifold. Then, its integral curves are contained in the submanifold which is, therefore, totally geodesic. If $|II|^2=-c$ and ${\bf H}\neq 0$ with $X^\top=T\equiv 0$ then ${\bf H}=cX^\perp = cX$ what implies that $|\psi|^2=|X|^2=c^2|{\bf H}|^2$. Since $\nabla^\perp {\bf H}=0$ we have 
$\nabla |{\bf H}|^2=0$ and then $|\psi|^2={\rm constant}$. In this case, $\psi(M)$ is immersed into a sphere. In general we observe that (\ref{solitonC-2}) may be written in local coordinates as
\[
\nabla^\perp_{\partial_i} {\bf H}^\alpha = -c\, h^\alpha_{ij} T^j.
\]
Differentiating both sides we get
\begin{eqnarray*}
& & -\nabla^\perp_{\partial_i}\nabla^\perp_{\partial_j} {\bf H} =  -\nabla^\perp_{\partial_i}(\nabla^\perp_{\partial_j} {\bf H}) + \nabla^\perp_{\nabla_{\partial_i}\partial_j} {\bf H}= c\,\nabla_{\partial_i}^\perp (II (T, \partial_j)) - c\, II(T, \nabla_{\partial_i}\partial_j)\\
& & \,\, = c (\nabla^\perp_{\partial_i} II) (T, \partial_j) + c\,II(\nabla_{\partial_i} T, \partial_j) + c\, II(T, \nabla_{\partial_i}\partial_j) - c\, II(T, \nabla_{\partial_i}\partial_j)\\
& & \,\, = c (\nabla^\perp_{\partial_i} II) (T, \partial_j) + c\,II(\nabla_{\partial_i} T, \partial_j). 
\end{eqnarray*}
However
\begin{eqnarray*}
& &  c\,II(\nabla_{\partial_i} T, \partial_j) = c (\nabla_{\partial_i} T)^k II(\partial_k, \partial_j) = c (\bar\nabla_{\partial_i} X-\bar\nabla_{\partial_i} X^\perp)^k II(\partial_k, \partial_j) \\
& & \,\,= c\varphi\delta_i^k II(\partial_k, \partial_j) -(\bar\nabla_{\partial_i} {\bf H})^k II(\partial_k, \partial_j)  = c\varphi II(\partial_i, \partial_j)-g^{ k\ell }II_{-{\bf H}} (\partial_i, \partial_\ell) II (\partial_k, \partial_j).
\end{eqnarray*}
Therefore
\begin{equation}
\nabla^\perp_{\partial_i}\nabla^\perp_{\partial_j} {\bf H}  =- c (\nabla^\perp_{\partial_i} II) (T, \partial_j) -c\varphi II(\partial_i, \partial_j)+g^{ k\ell }II_{-{\bf H}} (\partial_i, \partial_\ell) II (\partial_k, \partial_j).
\end{equation}
In terms of components
\[
\nabla^\perp_{\partial_i}\nabla^\perp_{\partial_j} {\bf H}^\alpha = -c\, T^k \nabla^\perp_{\partial_i } h^\alpha_{jk} -c\varphi h^\alpha_{ij} + (II_{-{\bf H}})^k_i  h^\alpha_{kj}. 
\]
Since $\nabla^\perp {\bf H}=0$ we conclude that
\begin{equation*}
II(T,\cdot) =0
\end{equation*}
and
\begin{equation}
c (\nabla^\perp_{\partial_i} II) (T, \partial_j) +c\varphi II(\partial_i, \partial_j)-g^{ k\ell }II_{-{\bf H}} (\partial_i, \partial_\ell) II (\partial_k, \partial_j) =0.
\end{equation}
In coordinates
\[
c T^k \nabla^\perp_{\partial_i } h^\alpha_{jk} +c\varphi h^\alpha_{ij}- (II_{-{\bf H}})^k_i  h^\alpha_{kj}=0.
\]
Now, differentiating
\[
|II|^2=-c
\]
in the direction of $T$ and with respect to the normal connection  yields (using Codazzi's equation for space forms)
\begin{eqnarray*}
& & 0 = \frac{c}{2} \langle \nabla |II|^2, T\rangle = \frac{c}{2}\sum_{\alpha} T \big( (h^\alpha)^{ij} h^\alpha_{ij}\big) = c(h^\alpha)^{ij}\nabla^\perp_T h^\alpha_{ij} = c(h^\alpha)^{ij} T^k \nabla^\perp_{\partial_k} h^\alpha_{ij}\\
& & \,\, = c(h^\alpha)^{ij} T^k \nabla^\perp_{\partial_i} h^\alpha_{jk} = - c\varphi (h^\alpha)^{ij} h^\alpha_{ij} + (II_{-{\bf H}})_{ik}  (h^\alpha)^k_{j}(h^\alpha)^{ij}.
\end{eqnarray*}
We conclude that
\begin{equation}
\label{IIH}
(II_{-{\bf H}})_{ik}  (h^\alpha)^k_{j}(h^\alpha)^{ij}  = c|II|^2 = -c^2,
\end{equation}
that is,
\[
{\bf H}^\beta h^\beta_{ik} (h^\alpha)^k_j (h^\alpha)^{ij} = -c|II|^2 = c^2.
\]
Now we use Simons' equation (see for instance Proposition 1.4 in \cite{AMR})
\begin{eqnarray*}
& & 0 = \frac{1}{2}\Delta |II|^2 = |\nabla^\perp II|^2 + (h^\alpha)^{ij} ({\rm tr} II_\alpha)_{;ij} - (h^\alpha)^{ij} (h^\beta)_{ij}  (h^\alpha)_{k\ell} (h^\beta)^{k\ell} +(h^\alpha)_{ij} (h^\alpha)^i_{k}  (h^\beta)^{jk} (h^\beta)_{\ell}^\ell\\
& &\,\, + 2 (h^\alpha)_{ij} (h^\alpha)^{k\ell}  (h^\beta)^{jk} (h^\beta)^{i}_\ell -  2 (h^\alpha)_{ik} (h^\alpha)^{kj}  (h^\beta)^{i\ell} (h^\beta)_{j\ell}.
\end{eqnarray*}
In an open subset $U\subset M$ where ${\bf H}\neq 0$ we set 
\[
N_{n-m+1} = \frac{{\bf H}}{|{\bf H}|}
\]
Therefore for all $\alpha $ we have
\[
{\rm tr} \,II_\alpha =  g^{ij}\langle II (\partial_i, \partial_j), N_\alpha\rangle = \underbrace{\langle  {\bf H}, N_\alpha\rangle}_{={\bf H}_\alpha} = |{\bf H}| \delta_{\alpha}^{n-m+1}
\]
Since $|{\bf H}|={\rm constant}$ it follows that 
\[
({\rm tr} II_\alpha)_{;ij} =0.
\]
Moreover
\[
 g^{ik}g^{j\ell}\langle II (\partial_i, \partial_j), {\bf H}\rangle \langle II (\partial_k, \partial_\ell), {\bf H}\rangle \le |II|^2 |{\bf H}|^2.
\]
In coordinates
\[
\sum_{\alpha, \beta}  h^\alpha_{ij} (h^\beta)^{ij}  {\bf H}_\alpha  {\bf H}_\beta\le |II|^2 |{\bf H}|^2.
\]
However
\begin{eqnarray*}
& & \nabla_T |II-{\bf H}g|^2 = \langle \nabla^\perp_T II - \nabla^\perp_T {\bf H} g, II-{\bf H}g\rangle = \frac{1}{2} T|II|^2 - {\bf H}^\alpha g^{ij} \nabla^\perp_T h^\alpha_{ij}\\
& &\,\, = - {\bf H}^\alpha \nabla^\perp_T g^{ij}h^\alpha_{ij} = -{\bf H}^\alpha \nabla^\perp_T {\bf H}_\alpha  = - \nabla_T |{\bf H}|^2 =0.
\end{eqnarray*}
Therefore
\[
II = {\bf H}g
\]
what implies that
\[
II_\alpha =0
\]
for all $\alpha\neq n-m+1$ and
\[
-c=|II|^2 = (h^{n-m+1})^{ij} h^{n-m+1}_{ij}
\]
with
\[
h^{n-m+1}_{ij} = |{\bf H}| g_{ij}.
\]
Using all this information in Simons' formula we get
\begin{eqnarray*}
|\nabla^\perp II|^2  =0
\end{eqnarray*}
on $U\subset M$, that is,
\[
\nabla^\perp II =0 \quad \mbox{on} \quad U.
\]
The proof now follows as in \cite{CL} at the end of the proof of Theorem 1.1. 

\vspace{3mm}

\subsection{First applications in the codimension one case: Einstein ambient spaces}
\label{subsectionEinstein}
In what follows we will be mainly concerned with the codimension one case, where the equations above become simpler.

Indeed, let $\psi: M^m \to \bar M^{m+1}$ be a codimension one mean curvature flow soliton with respect to a closed conformal vector field $X\in\Gamma(T\bar M)$ on $\bar M$. Having chosen a local unit normal $N$ to the immersion, we have 
\[
{\bf H}=HN, \quad |{\bf H}|^2=H^2, \quad |II_{\bf H}|^2=H^2|A|^2 \quad \text{and} \quad 
|(\bar\nabla_{(\cdot)}{\bf H})^\perp|^2=|\nabla H|^2,
\]
where 
\begin{equation}
\label{luis7}
H =\langle {\bf H}, N\rangle = c\,\langle X, N\rangle
\end{equation}
and $A$ denotes the Weingarten operator in the unit direction $N$. From (\ref{closedconf-bis}) and (\ref{luis7}) it follows directly that
\begin{equation}
\label{luis23}
\nabla H=-cAT
\end{equation}
where we recall here that $T=X^\top$ is the tangential component of $X$. Observe that (\ref{closedconf-bis}) implies also that
\begin{equation}
\label{luis21}
\nabla_UT=\varphi U+\frac{H}{c}AU
\end{equation}
for any $U\in\Gamma$, where $\varphi=\bar{\varphi}\circ\psi$. Hence, using Codazzi equation, for any $U,V\in\Gamma(TM)$ we have 
\begin{eqnarray}
\label{luis22}
\langle \nabla_U \nabla H, V\rangle & = &  -c\langle (\nabla_U A)T,
V\rangle -c\varphi\langle AU, V\rangle -
H\langle AU, AV\rangle\\
\nonumber {} & = & -c\langle (\nabla_T A)U,
V\rangle +c\langle\bar{R}(U,T)N,V\rangle
-c\varphi\langle AU, V\rangle -H\langle AU, AV\rangle.
\end{eqnarray}
Therefore
\begin{equation*}
\Delta H=-c\langle\nabla H,T\rangle+c{\rm Ric}_{\bar M}(T,N)-(c\varphi+|A|^2)H.
\end{equation*}
In other words
\begin{equation}
\label{kk1}
\Delta_{-cT}H=c{\rm Ric}_{\bar M}(T,N)-(c\varphi+|A|^2)H.
\end{equation}
In particular, if $\bar M$ is Einstein we deduce
\begin{equation}
\label{deltaXH-einstein}
\Delta_{-cT} H= -(c\varphi+|A|^2)H,
\end{equation}
and
\begin{equation}
\label{delta-eta-H2-bis}
\Delta_{-cT} H^2= -2(c\varphi+|A|^2)H^2+2|\nabla H|^2.
\end{equation}
As a first application of the latter equation we have the following result, which compares with Theorem 1.1 in \cite{cheng-peng}. 

\begin{theorem}
\label{gap-einstein}
Let $\psi: M^m \to \bar M^{m+1}$ be a codimension one mean curvature flow soliton with respect to a closed conformal vector field $X\in\Gamma(T\bar M)$ and $c<0$. Suppose that $\bar M$ is Einstein. Assume the validity of the weak maximum principle for the operator $\Delta_{-cT}$ on $M$, where $T = X^\top$. Finally suppose that
\begin{equation}
\label{bound-H}
\sup_M H^2 <+\infty.
\end{equation} 
Then either 
\begin{equation}
\label{gap-e1}
\sup_M (|A|^2 +c\varphi)\ge 0,
\end{equation}
where $\varphi = \frac{1}{m+1}{\rm div}_{\bar M}X\circ\psi$, or
\begin{equation}
\label{gap-e2}
H\equiv 0 \,\, \mbox{ on } \,\, M.
\end{equation}
\end{theorem}

\begin{remark}
We note that {\rm (\ref{gap-e1})} in particular implies 
\begin{equation}
\label{gap-e3}
\sup_M |A|^2 \ge -c\inf_M \varphi.
\end{equation}
\end{remark}

\noindent \emph{Proof.} Since $\bar M$ is Einstein 
from (\ref{delta-eta-H2-bis}) we have
\begin{equation}
\label{laplace-H2k}
\frac{1}{2}\Delta_{-cT} H^2 = (-c\varphi-|A|^2)H^2+|\nabla H|^2.
\end{equation}
If 
\[
\sup_M (|A|^2+c\varphi)\ge 0
\]
there is nothing to prove. Otherwise, suppose that
\[
\sup_M (|A|^2+c\varphi) < 0.
\]
Then $|A|^2< -c\varphi$ on $M$, so that
\begin{equation}
\label{chain1}
H^2 \le m |A|^2 <-mc\varphi \le -mc \sup_M \varphi.
\end{equation}
Furthermore,
\begin{equation}
\label{chain2}
\inf_M (-c\varphi-|A|^2)=C>0.
\end{equation}
By the weak maximum principle for the operator $\Delta_{-cT}$ on $M$, there exists a sequence
$\{x_k\}_{k=1}^\infty$ in $M$ such that
\begin{equation}
H^2(x_k) >\sup_M H^2-\frac{1}{k} \,\,\, \mbox{ and } \,\,\, \Delta_{-cT} H^2 (x_k) <\frac{1}{k}\cdot
\end{equation}
Using  (\ref{laplace-H2k}) it follows that
\begin{eqnarray*}
& &  2 (-c\varphi(x_k)-|A|^2(x_k))H^2(x_k)\le 2|\nabla H|^2(x_k)+2
(-c\varphi(x_k)-|A|^2(x_k))H^2(x_k)\\
& &\,\,= \Delta_{-cT} H^2(x_k)< \frac{1}{k}
\end{eqnarray*}
from which we infer
\[
C\,\bigg(\sup_M H^2-\frac{1}{k}\bigg)
<\inf_M(-c\varphi(x_k)-|A|^2(x_k))\, H^2(x_k)<\frac{1}{2k}\cdot
\]
Passing to the limit as $k\to +\infty$ we conclude that
\[
\sup_M H^2=0,
\]
that is,  $H\equiv 0$  on $M$.  \hfill $\square$

Similarly, we also get the following result,  which compares with Theorem 1.2 in \cite{cheng-peng}.
\begin{theorem}
\label{gap-F}
Let $\psi:M^m\to \bar M^{m+1}$ be a codimension one, oriented mean curvature flow soliton with respect to a closed conformal vector field $X\in\Gamma(T\bar M)$ and $c<0$. Suppose that $\bar M$ is Einstein. Assume the validity of the weak maximum principle for the operator $\Delta_{-cT}$ on $M$. Finally suppose that 
$\inf_M H^2 >0$ and $\sup_M |A|^2<+\infty$.
Then 
\begin{equation}
\label{Avarphi}
\inf_M (|A|^2+c\varphi)\le 0,
\end{equation}
where
\[
\varphi = \frac{1}{m+1}\,{\rm div}_{\bar M} X \circ\psi.
\]
\end{theorem}

\begin{remark}
We note that {\rm (\ref{Avarphi})} in particular implies
\begin{equation}
\label{Avarphi-bis}
\inf_M |A|^2 \le -c\sup_M \varphi.
\end{equation}
\end{remark}

\noindent \emph{Proof of Theorem \ref{gap-F}.} Up to choosing an appropriate unit normal $N$ to the hypersurface we can suppose that for $H=\langle {\bf H}, N\rangle$,
\[
\inf_M H = C_1 >0.
\]
We set $u=-H$ and we use (\ref{deltaXH-einstein}) to obtain
\begin{eqnarray*}
& & \Delta_{-cT} u= (c\varphi+|A|^2)H \ge \inf_M (|A|^2+c\varphi) H.
\end{eqnarray*}
Thus assuming by contradiction that
\[
\inf_M (|A|^2+c\varphi)=C_2>0
\]
we deduce 
\[
\Delta_{-cT} u \ge C_1C_2>0.
\]
An application of the weak maximum principle directly yields the desired contradiction. \hfill $\square$

\


%

On the other hand, since 
\[
4 H^2 |\nabla H|^2=|\nabla H^2|^2,
\]
equation (\ref{delta-eta-H2-bis}) yields
\begin{eqnarray}
\label{lapl-Hsq-bis}
\nonumber H^2 \Delta_{-cT} H^2 & = & -2(c\varphi+|A|^2) H^4 + 2 |\nabla H|^2H^2\\
{} & = & -2(c\varphi+|A|^2)H^4+
\frac{1}{2}|\nabla H^2|^2.
\end{eqnarray}
In other words, in the codimension one case and if $\bar M$ is Einstein the function $u = H^2\ge 0$ satisfies
\begin{equation}
\label{kato-1-bis}
u \Delta_{-cT} u+2 (c\varphi+|A|^2) u^2=\frac{1}{2}|\nabla u|^2 \quad \text{on} \,\,\, M.
\end{equation}
In particular, when $\bar M=I\times_hP$ is a warped product space we have $T=\nabla\eta$ and $\varphi=h'$, and (\ref{kato-1-bis}) becomes
\begin{equation}
\label{kato-1}
u \Delta_{-c\eta} u+2 (ch'+|A|^2) u^2=\frac{1}{2}|\nabla u|^2  \quad \text{on} \,\,\, M.
\end{equation}
As a consequence of the above equations we have
\begin{theorem}
\label{gap-Eb}
Let $\psi: M^m \to \bar M^{m+1} = I\times_h P$ be a complete codimension one  mean curvature flow soliton with respect to $X= h(t)\partial_t$, with $\bar M$ Einstein. Assume that, for some $\mu>1/2$,
\[
\lambda_1^{L_\mu}(M)\ge 0,
\]
where $L_\mu$ is the operator
\[
L_\mu = \Delta_{-c\eta}+2 \mu (ch'+|A|^2).
\]
Suppose that, for some $1\le q\le 2\mu$
\begin{equation}
\label{H-p}
\bigg(\int_{\partial B_r}H^{2q} e^{c\eta}\, {\rm d}M\bigg)^{-1} \notin L^1(+\infty).
\end{equation}
Then, $H$ is constant on $M$. Furthermore, if $|A|^2\not\equiv -ch'$ then $H\equiv 0$.
\end{theorem}


In order to prove Theorem \ref{gap-Eb} we need to recall the next result, whose proof is a minor variation of that of Theorem 4.5 in \cite{PRSgreenbook}.
%
%
\begin{theorem}
\label{29.1}
Let $(M,\langle,\rangle,e^{-f}dM)$ be a complete weighted manifold with $f\in C^\infty(M)$ and $a(x)\in L^\infty_{\rm loc}(M)$. Let $u\in{\rm Lip}_{\rm loc}(M)$ satisfy the differential inequality
\begin{equation}
\label{29.2}
u\Delta_fu+a(x)u^2+\alpha|\nabla u|^2\ge 0 \quad \text{weakly on } M
\end{equation}
for some $\alpha\in\mathbb{R}$. Let $v\in{\rm Lip}_{\rm loc}(M)$ be a positive solution of 
\begin{equation}
\label{29.3}
\Delta_fv+\mu a(x)v\le 0 \quad \text{weakly on } M
\end{equation}
for some $\mu$, and suppose that
\begin{equation}
\label{29.4}
\alpha+1\le \mu, \quad  \mu>0.
\end{equation}
If 
\begin{equation}
\label{29.6}
\left(\int_{\partial B_r}u^{2(\beta+1)}e^{-f}\right)^{-1}\notin L^1(+\infty)
\end{equation}
holds for some $\beta$ satisfying
\begin{equation}
\label{29.5}
\alpha\le \beta\le \mu -1,\quad \beta>-1.
\end{equation}
then there exists a constant $C\in\mathbb{R}$ such that
\begin{equation}
\label{29.7}
Cv={\rm sgn} u\, |u|^\mu.
\end{equation}
Furthermore,
\begin{enumerate}
\item[(i)] If $\alpha+1<\mu$, then $u$ is constant on $M$, and if in addition $a(x)$ does not vanish identicall, then $u\equiv 0$. 
\item[(ii)] If $\alpha+1=\mu$ and $u$ does not vanish identically, then $v$ and therefore $|u|^\mu$ satisfy (\ref{29.3}) with equality sign.
\end{enumerate}
\end{theorem}

\begin{remark}
\label{29.8}
Consider the operator
\begin{equation}
\label{29.9}
L^\mu_f=\Delta_f+\mu a(x)
\end{equation}
and set
\begin{equation}
\label{29.10}
\lambda_1^{L^\mu_f}=\inf\bigg\{\frac{\int_M(|\nabla w|^2-\mu a(x)w^2)e^{-f}}{\int_Mw^2e^{-f}}: 
w\in C^\infty_c(M), w\not\equiv 0\bigg\}
\end{equation}
to denote the spectral radius of the operator $L^\mu_f$ on $(M,\langle,\rangle,e^{-f}dM)$. We recall that, according to a minor variation of a result of Fischer-Colbrie and Schoen {\rm \cite{FCS}}, the requirement 
\begin{equation}
\label{29.11}
\lambda_1^{L^\mu_f}\ge 0
\end{equation}
is equivalent to the existence of a positive solution $v$ of 
\begin{equation}
\label{29.12}
L^\mu_fv=\Delta_fv+\mu a(x)v=0,
\end{equation}
that is, $v$ solves {\rm (\ref{29.3})} with equality sign. 
\end{remark}

\noindent \emph{Proof of Theorem \ref{gap-Eb}.} The proof is now an immediate consequence of Theorem \ref{29.1} and Remark \ref{29.8}, and of the differential equality (\ref{kato-1}), with the choices $u=H^2\ge 0$, $f=-c\eta$, $a(x)=2(ch'+|A|^2)$, $\alpha=-\frac{1}{2}$, $\mu>\frac{1}{2}$ and $\beta=\frac{q}{2}-1$, so that $H^{2q}=u^{2(\beta+1)}$.  This completes the proof. \hfill $\square$

\vspace{3mm}

\begin{remark}
Note that condition {\rm (\ref{H-p})} is implied by 
\begin{equation}
\label{H-p2}
H^{2q}\in  L^1(M, e^{c\eta}\, {\rm d}M), \quad 1\le q\le 2\mu
\end{equation}
\end{remark}

\subsection{Further applications}
In what follows we consider the case where $\bar M^{m+1}=I\times_{h}P^m$ is a warped product space, not necessarily Einstein. In this case we know that $T=\nabla\eta$ and $\varphi=h'$, so that (\ref{kk1}) becomes
\begin{equation}
\label{kk2}
\Delta_{-c\eta}H=c\,{\rm Ric}_{\bar M}(\nabla\eta,N)-(ch'+|A|^2)H.
\end{equation}
Observe now that 
\begin{equation}
\label{luis6}
{\rm Ric}_{\bar M}(Z,X)=-m\frac{h''}{h}\langle Z,X\rangle
\end{equation} 
for every vector field $Z\in\Gamma(T\bar M)$. Therefore, using (\ref{luis7}) and the soliton equation we have
\begin{eqnarray*}
c\,{\rm Ric}_{\bar M}(\nabla\eta,N) & = & -{\rm Ric}_{\bar M}({\bf H},N)+c{\rm Ric}_{\bar M}(X,N)\\
{} & = & -H\left({\rm Ric}_{\bar M}(N,N)+m\frac{h''}{h}\right),
\end{eqnarray*}
which jointly with (\ref{kk2}) implies
\begin{equation}
\label{kk3}
\Delta_{-c\eta}H=-\left({\rm Ric}_{\bar M}(N,N)+m\frac{h''}{h}\right)H
-c(h'+|A|^2)H.
\end{equation}
Therefore, 
\begin{eqnarray}
\label{deltaXH-warped-bis} 
\frac{1}{2}\Delta_{-c\eta}H^2 & = & H\Delta_{-c\eta}H+|\nabla H|^2 \\
\nonumber {} & = &  -(ch'+|A|^2)H^2+ |\nabla H|^2-H^2
\left({\rm Ric}_{\bar M}(N,N)+m\frac{h''}{h}\right).
\end{eqnarray}
As a consequence of this equation we have the following

\begin{theorem}
\label{gap-E}
Let $\psi: M^m \to \bar M^{m+1} = (a,+\infty)\times_h P$, $a\ge-\infty$, be a complete, orientable codimension one mean curvature flow soliton with respect to the vector field $X= h(t)\partial_t$ with $h\notin L^1(+\infty)$. Let $N$ be a chosen unit normal vector field along $\psi$ and  let $H$ be given by ${\bf H}=HN$. Assume
\begin{itemize}
\item[{\rm (i)}] ${\rm Ric}_{\bar M} (N,N) \le K_{\bar M}$ for some constant $K_{\bar M}\in \mathbb{R}$
\item[{\rm (ii)}] $\lim_{x\to \infty} (\pi\circ\psi)(x)= +\infty$
\item[{\rm (iii)}] $\limsup_{t\to +\infty}\zeta(t) <0$
\end{itemize}
where
\[
\zeta(t) = mh'(t)+ch^2(t)
\]
is the soliton function. Furthermore, suppose that there exists $\gamma\in \mathbb{R}$ such that
\begin{itemize}
\item[{\rm (iv)}] $\Omega_\gamma = \{x\in  M: H^2(x)>\gamma\}\neq \emptyset$
\item[{\rm (v)}] $\Gamma =\inf_{[\gamma, +\infty)} \Big(-ch'-m\frac{h''}{h}\Big) - K_{\bar M}\ge 0.$
\item [{\rm (vi)}] $|A|^2\le \Gamma$ on $\Omega_\gamma$. 
\end{itemize}
Then either $H\equiv 0$ or $|A|^2=\Gamma$ and $H$ is constant on $M$.
\end{theorem}

\begin{remark}
Note that the conditions on the completeness of $M$  together with assumptions {\rm ii)} and {\rm iii)} in the statement of Theorem \ref{gap-E}
will be used only to show that $M$ is $\Delta_{-c\eta}$-parabolic.  Clearly, for $\gamma>0$ the alternative $H\equiv 0$ cannot happen. 
\end{remark}

\vspace{3mm}

\noindent \emph{Proof.}  By assumptions ii) and iii), using equation (\ref{delta-eta-X}) we have that
\[
\Delta_{-c\eta}\eta=\zeta(\pi\circ\psi)<0 \,\, \mbox{ on }\,\, M\backslash K
\]
for some compact $K\subset M$ where $\eta(x) = \hat\eta((\pi\circ\psi)(x))$. From $h\notin L^1(+\infty)$, assumption ii)  and the definition of $\eta$ we deduce that
\[
\eta(x) \to +\infty\,\, \mbox{ as }\,\, x\to\infty \,\, \mbox{ in } \,\, M. 
\]
Since $M$ is complete, from Theorem 4.12 of \cite{AMR} we deduce that $M$ is $\Delta_{-c\eta}$-strongly parabolic. In particular, it is parabolic in the usual sense and Ahlfors parabolic (see page 246 of \cite{AMR}). 
Because of i), v), and vi) from (\ref{deltaXH-warped-bis}) we obtain
\begin{eqnarray}
\label{sub-H2}
\nonumber \frac{1}{2}\Delta_{-c\eta} H^2 & \ge & 
\left(-ch'-m\frac{h''}{h}-K_{\bar M}-|A|^2\right)H^2 \\
{} & \ge & (\Gamma-|A|^2) H^2 \ge 0 \quad \mbox{ on } \quad \Omega_\gamma.
\end{eqnarray}
If $\Omega_\gamma \not\equiv M$ then $\partial\Omega_\gamma \neq \emptyset$. Furthermore, because of vi) 
\begin{equation}
\label{supH2}
\sup_{\Omega_\gamma} H^2 < +\infty.
\end{equation}
Hence, by Ahlfors parabolicity of $\Delta_{-c\eta}$, (\ref{sub-H2}) and (\ref{supH2}) yield
\[
\sup_{\Omega_\gamma} H^2 = \sup_{\partial\Omega_\gamma} H^2 = \gamma,
\]
which is a contradiction because $\Omega_\gamma\neq\emptyset$. It follows that $\Omega_\gamma \equiv M$ and $H^2$ is bounded above on $M$. Since $M$ is $\Delta_{-c\eta}$-parabolic from (\ref{sub-H2}), now on $M$, it follows that $H^2$ is constant. If $H\not\equiv 0$, from (\ref{sub-H2}) we deduce $|A|^2=\Gamma$  and from (\ref{deltaXH-warped-bis}) it results that $|\nabla H|=0$ completing the proof. \hfill $\square$

\vspace{3mm}

It is interesting to analyze Theorem \ref{gap-E} in some special cases. For instance,
\begin{corollary}
\label{cor-gap-E1} 
Let $\psi: M^m \to \mathbb{R}^{m+1}$ be a complete, orientable, codimension one self-shrinker {\rm(}that is, $c<0${\rm)}. Assume that,  for some fixed ${\sf v}\in \mathbb{S}^m$
\begin{equation}
\label{sol-cone-cond}
\lim_{x\to \infty}\langle \psi(x), {\sf v}\rangle = +\infty.
\end{equation}
Suppose that there exists $\gamma\in \mathbb{R}$ such that
\[
\Omega_\gamma =\{x\in M: H^2(x)>\gamma\}\neq \emptyset
\]
and 
\[
|A|^2 \le -c \quad \mbox{ on } \quad \Omega_\gamma.
\]
Then, either $H\equiv 0$ or $|A|^2\equiv -c$ and $H$ is constant. In the first case $\psi$ is a totally geodesic hyperplane while in the second $\psi(M)$ is a product of the form $\mathbb{S}^{n}(\sqrt{-n/c})\times \mathbb{R}^{m-n}$ for some $1\le n\le m-1$. 
\end{corollary}

\noindent \emph{Proof of Corollary \ref{cor-gap-E1}.} Since $\psi(M)$ is necessarily different from $\mathbb{R}^{m+1}$ because of dimensional reasons we can fix an origin $o\notin\psi(M)$ in $\mathbb{R}^{m+1}$. By definition $\psi$ is a mean curvature flow soliton with respect to $X=t\partial_t$ and $c<0$, where we are representing $\mathbb{R}^{m+1}\backslash\{o\}$ as $(0,+\infty)\times_{t} \mathbb{S}^m$ having $(0,+\infty)$ in the direction determined by ${\sf v}$. Then (\ref{sol-cone-cond}) yields $(\pi\circ\psi)(x)\to +\infty$ as $x\to \infty$. In this case, $h(t)=t$, $K_{\bar M}=0$, $\zeta(t)=m+ct^2$ with $c<0$ and $\Gamma=-c>0$.
It then follows from Theorem \ref{gap-E} that either $H\equiv 0$ or $|A|^2=-c$. The last statement follows immediately from \cite{CL}. \hfill $\square$


\begin{remark} Since $H^2 \le m |A|^2$ and $|A|^2\le -c$ on $\Omega_\gamma$  necessarily $\gamma\in (-\infty,-mc)$.
\end{remark}

Corollary \ref{cor-gap-E1} is a slight extension of a result of Cao and Li in \cite{CL}, where the codimension is arbitrary. As we shall see next, Theorem  \ref{gap-E} indeed extends to arbitrary codimension giving, in favourable circunstaces, interesting geometric conclusions. Specifically, we can state the following general result.
\begin{theorem}
\label{gap-E-codimension}
Suppose that $\bar M^{n+1}=(a,+\infty)\times_h P$, $a\ge-\infty$, is a warped product space with constant sectional curvature and let $\psi: M^m \to \bar M^{n+1}$ be a complete mean curvature flow soliton with respect to the vector field $X= h(t)\partial_t$ with $h\notin L^1(+\infty)$. Assume
\begin{itemize}
\item[{\rm (i)}] $\lim_{x\to \infty} (\pi\circ\psi)(x)= +\infty$
\item[{\rm (ii)}] $\limsup_{t\to +\infty}\zeta(t) <0$
\end{itemize}
where
\[
\zeta(t) = mh'(t)+ch^2(t)
\]
is the soliton function. Furthermore, suppose that there exists $\gamma\in \mathbb{R}$ such that
\begin{itemize}
\item[{\rm (iii)}] $\Omega_\gamma = \{x\in  M: |{\bf H}|^2(x)>\gamma\}\neq \emptyset$
\item[{\rm (iv)}] $\Gamma =\inf_{[\gamma, +\infty)} \Big(-ch'\Big)\ge 0.$
\item [{\rm (v)}] $|II|^2\le \Gamma$ on $\Omega_\gamma$. 
\end{itemize}
Then either ${\bf H}\equiv 0$ or $|II|^2=\Gamma$ on $M$ and ${\bf H}$ is parallel .
\end{theorem}

\noindent \emph{Proof.} Reasoning as in the proof of Theorem \ref{gap-E}, 
we deduce that $\Delta_{-c\eta}$ is Ahlfors parabolic on $M$. 
Because of iv) and v), from (\ref{deltaH2-kappa-bis}) we obtain using $|II_{\bf H}|^2\le|II|^2|{\bf H}|^2$ 
\begin{equation}
\label{sub-H2-cod}
\frac{1}{2}\Delta_{-c\eta} |{\bf H}|^2\ge
\left(-ch'-|II|^2\right)|{\bf H}|^2
\ge (\Gamma-|II|^2)|{\bf H}|^2 \ge 0 \quad \mbox{ on } \quad \Omega_\gamma.
\end{equation}
If $\Omega_\gamma \not\equiv M$ then $\partial\Omega_\gamma \neq \emptyset$. Furthermore, because of v) and using $|{\bf H}|^2\leq m|II|^2$ we have
\begin{equation}
\label{supH2-cod}
\sup_{\Omega_\gamma} |{\bf H}|^2 < +\infty.
\end{equation}
Hence, by Ahlfors parabolicity of $\Delta_{-c\eta}$, (\ref{sub-H2-cod}) and (\ref{supH2-cod}) yield
\[
\sup_{\Omega_\gamma} |{\bf H}|^2 = \sup_{\partial\Omega_\gamma} |{\bf H}|^2 = \gamma,
\]
which is a contradiction because $\Omega_\gamma\neq\emptyset$. It follows that $\Omega_\gamma \equiv M$ and $|{\bf H}|^2$ is bounded above on $M$. Since $M$ is $\Delta_{-c\eta}$-parabolic, from (\ref{sub-H2-cod}), now on $M$, it follows that $|{\bf H}|^2$ is constant. If ${\bf H}\not\equiv 0$, from (\ref{sub-H2-cod}) we deduce $|II|^2=\Gamma$  and from (\ref{deltaH2-kappa-bis}) it results that $|(\bar\nabla {\bf H})^\perp|=0$ completing the proof. \hfill $\square$

\vspace{3mm}

As a direct consequence of Theorem \ref{gap-E-codimension} we have the following 
\begin{corollary}
\label{gap-Fbis} Let $\psi: M^m \to \mathbb{R}^{n+1}$, $n\ge m$, be a complete self-shrinker (that is, $c<0$).  Assume that, for some ${\sf v}\in \mathbb{S}^n$, 
\begin{equation}
\langle \psi(x), {\sf v}\rangle \to +\infty\,\, \mbox{ as } \,\, x\to \infty \,\, \mbox{ in }\,\, M
\end{equation}
and let 
\begin{equation}
\label{Al1}
|II|^2 \le -c.
\end{equation}
Then either ${\bf H} \equiv 0$ or $|II|^2 \equiv -c$ and ${\bf H}$ is parallel and Remark \ref{remark8.1} applies. 
\end{corollary}

\noindent \emph{Proof.} Observe that $\psi(M)\neq \mathbb{R}^{n+1}$ since $n\ge m$. Fix an origin $o\in \mathbb{R}^{n+1}$ such that $o\notin \psi(M)$. Consider the axis determined by ${\sf v}$ and represent $\mathbb{R}^{n+1}\backslash\{o\}$ as $(0,+\infty)_{t}\times \mathbb{S}^n$ with warping function $h(t)=t$. With these choices, if ${\bf H}\not\equiv 0$ then Theorem \ref{gap-E-codimension} applies directly with $\gamma=0$ to conclude that $|II|^2=-c$ on $M$ and ${\bf H}$ is parallel. 
At this point, Remark \ref{remark8.1} applies.   \hfill $\square$

\

Observe that representing either 
$\mathbb{H}^{n+1}\backslash\{o\} = (0,+\infty)\times_{\sinh t} \mathbb{S}^n$ or $\mathbb{H}^{n+1}=\mathbb{R}\times_{e^t} \mathbb{R}^n$, conclusions that parallel those of Theorem \ref{gap-E-codimension} can be drawn in this case too. For instance we have 
\begin{corollary}
\label{gap-Fbisbis} Let $\psi: M^m \to \mathbb{H}^{n+1}=\mathbb{R}\times_{e^t} \mathbb{R}^n$, $n\ge m$, be a complete self-shrinker with respect to $X=e^t\partial_t$ (that is, $c<0$).  Suppose that
\begin{equation}
\pi\circ\psi(x) \to +\infty\,\, \mbox{ as } \,\, x\to \infty \,\, \mbox{ in }\,\, M
\end{equation}
and 
\begin{equation}
\label{Al1-bis}
|II|^2 \le -c.
\end{equation}
Then either ${\bf H} \equiv 0$ or $|II|^2 \equiv -c$ on $M$ and ${\bf H}$ is parallel. 
\end{corollary}

\section{Applications of a Simons' type formula}
\label{simons}

This section is devoted to deduce some geometric consequences of Simons' formula that we express, for codimension one mean curvature flow solitons, in a form convenient for our purposes. Note that in this section our target manifold is not necessarily always a warped product.

\begin{proposition}
\label{prop-simons-soliton}
 Let $\bar M^{m+1}$ be a Riemannian manifold  of constant
sectional curvature $\kappa$ and let $\psi:M^m\to \bar M^{m+1}$ be
a  codimension one mean curvature flow soliton with respect to a closed conformal vector field $X\in \Gamma(T\bar M)$. Denote with $A$ the  Weingarten operator in the direction of some local unit normal to the immersion. Then
\begin{equation}
\label{simons-soliton-spaceform} \frac{1}{2}\Delta_{-cT} |A|^2=|\nabla
A|^2 -(c\varphi+|A|^2)|A|^2 +(m |A|^2 - H^2)\, \kappa,
\end{equation}
where 
\[
\varphi = \frac{1}{m+1}\, {\rm div}_{\bar M} X \circ\psi.
\]
\end{proposition}

\noindent \emph{Proof.}  Simons' formula for hypersurfaces in   spaces of constant sectional curvature reads as
\begin{eqnarray}
\frac{1}{2}\Delta |A|^2-|\nabla A|^2={\rm tr}(A\circ \nabla\nabla H)+H
{\rm tr}A^3- |A|^4 +\kappa(m|A|^2-H^2).
\end{eqnarray}
For a proof of this formula, we refer to \cite[Lemma 6.1]{AMR}, paying attention that in our definition here we are taking $H=\mathrm{tr}(A)$ instead of $mH=\mathrm{tr}(A)$. From subsection \ref{subsectionEinstein} we recall that, in case of codimension one mean curvature flow solitons, we have
\[
H =c\langle X, N\rangle
\]
that gives (see (\ref{luis23}))
\[
\nabla H = -cAT.
\]
Since $\bar M$ has constant sectional curvature, (\ref{luis22}) becomes 
\begin{equation*}
\langle \nabla_U \nabla H, V\rangle = -c\langle (\nabla_T A)U,
V\rangle -c\varphi\langle AU, V\rangle -
H\langle AU, AV\rangle
\end{equation*}
for any $U, V\in \Gamma(TM)$. Therefore
\begin{equation*}
\langle A, \nabla\nabla H\rangle = -\frac{c}{2}\langle \nabla|A|^2,
T\rangle -c\varphi |A|^2 - H{\rm tr}A^3,
\end{equation*}
and we obtain
\begin{equation*}
\frac{1}{2}\Delta |A|^2-|\nabla A|^2=-\frac{c}{2}\langle
\nabla|A|^2, T\rangle -c\varphi |A|^2 - |A|^4 +\kappa(m|A|^2-H^2).
\end{equation*}
That is,
\begin{equation}
\label{simons-soliton-bis}
\frac{1}{2}\Delta_{-cT} |A|^2=|\nabla A|^2 -(c\varphi+|A|^2)|A|^2
+m\kappa |A|^2 -\kappa H^2,
\end{equation}
as desired. \hfill $\square$

\vspace{3mm}  

Motivated by a result of Huisken \cite{huisken}, and as a simple
application of the classical strong maximum principle, we prove

\begin{theorem}
\label{huisken-cl}
Let $\bar M^{m+1}$ be a Riemannian manifold of constant sectional curvature
$\kappa\ge 0$ and let $\psi:M^m\to \bar M^{m+1}$ be a
\emph{compact}, orientable codimension one mean curvature flow soliton with respect to a closed conformal vector field $X\in \Gamma(T\bar M)$. Assume that $H=\langle {\bf H}, N\rangle$ does not change sign on $M$, where $N$ is a unit normal to the immersion. 
\begin{itemize}
\item[{\rm (i)}] If $\kappa>0$ then either $H\equiv 0$ or the immersion is totally umbilical (and hence $H>0$ is constant).
\item[{\rm (ii)}] If $\kappa=0$ and $\varphi$ does not change sign on $M$, then $H$ is constant. 
\end{itemize}
\end{theorem}

\noindent \emph{Proof.} We recall from (\ref{luis23}) and (\ref{deltaXH-einstein}) that
\begin{equation}
\label{nablaH} \nabla H =-cAT
\end{equation}
and that for a manifold of constant sectional curvature $\kappa$
\begin{equation}
\label{deltaH-huisken} 
\Delta_{-cT} H =-(c\varphi+|A|^2) H
\end{equation}
with $H = \langle {\bf H}, N\rangle$ and  
\[
\varphi = \frac{1}{m+1}\, {\rm div}_{\bar M} X \circ\psi.
\]

Without loss of generality we can assume $H\ge 0$ on $M$. Thus if $H(x) =0$ for some $x\in M$, by the strong maximum principle (see \cite{GT} page 35) and (\ref{deltaH-huisken}) we have $H\equiv 0$ on $M$, in particular $H$ is constant. Hence we can suppose $H>0$ on $M$. In this case a computation gives
\begin{eqnarray*}
& &
\Delta_{-cT}\bigg(\frac{|A^2|}{H^2}\bigg)=\frac{1}{H^2}\Delta_{-cT}|A|^2
-\frac{4}{H^3}\langle \nabla H,\nabla |A|^2\rangle+|A|^2 \Delta_{-cT}
\frac{1}{H^2}\\
& & \,\,=\frac{1}{H^2}\Delta_{-cT}|A|^2 -\frac{4}{H^3}\langle \nabla
H,\nabla |A|^2\rangle - \frac{2}{H^3}|A|^2 \Delta_{-cT} H+
\frac{6}{H^4}|\nabla H|^2 |A|^2.
\end{eqnarray*}
Hence, using (\ref{simons-soliton-spaceform}) and (\ref{deltaH-huisken}) we obtain
\begin{eqnarray*}
\Delta_{-cT}\bigg(\frac{|A^2|}{H^2}\bigg)=
\frac{2}{H^2}\bigg(m-\frac{H^2}{|A|^2}\bigg)\kappa|A|^2
+\frac{2}{H^2}|\nabla A|^2-\frac{4}{H^3}\langle \nabla H,\nabla
|A|^2\rangle  + \frac{6}{H^4}|\nabla H|^2 |A|^2.
\end{eqnarray*}
We rewrite this expression in the form
\begin{eqnarray*}
& &\Delta_{-cT}\bigg(\frac{|A^2|}{H^2}\bigg)=
\frac{2}{H^2}\bigg(m-\frac{H^2}{|A|^2}\bigg)\kappa|A|^2
+\frac{2}{H^4}\big(H^2|\nabla A|^2-2\langle H\nabla H,\nabla
|A|^2\rangle  +3|\nabla H|^2
|A|^2\big)\\
& & \,\, = \frac{2}{H^2}\bigg(m-\frac{H^2}{|A|^2}\bigg)\kappa|A|^2
+\frac{2}{H^4}\big(-\langle H\nabla H,\nabla |A|^2\rangle
+2|\nabla H|^2 |A|^2\big)\\
& &\,\,+\frac{2}{H^4}\big(H^2|\nabla A|^2-\langle H\nabla H,\nabla
|A|^2\rangle  +|\nabla H|^2 |A|^2\big).
\end{eqnarray*}
Next we observe that an easy computation yields
\begin{eqnarray*}
\bigg\langle \nabla\bigg(\frac{|A|^2}{H^2}\bigg), \nabla
H\bigg\rangle=-\frac{1}{H^3}(-\langle H\nabla H, \nabla
|A|^2\rangle+2|\nabla H|^2|A|^2)
\end{eqnarray*}
and
\begin{eqnarray*}
|\nabla H \otimes A-H\nabla A|^2= |\nabla H|^2 |A|^2- H \langle
\nabla H, \nabla |A|^2\rangle +H^2|\nabla A|^2.
\end{eqnarray*}
Substituting into the above equation we finally obtain
\begin{eqnarray*}
\Delta_{-cT}\bigg(\frac{|A^2|}{H^2}\bigg) & = & \frac{2}{H^4}|\nabla H
\otimes A-H\nabla A|^2+
2\bigg(m-\frac{H^2}{|A|^2}\bigg)\kappa\frac{|A|^2}{H^2}\\
{} & {} & -\frac{2}{H}\bigg\langle \nabla\bigg(\frac{|A|^2}{H^2}\bigg),
\nabla H\bigg\rangle\\
{} & = & \frac{2}{H^4}|\nabla H
\otimes A-H\nabla A|^2+
2\bigg(m-\frac{H^2}{|A|^2}\bigg)\kappa\frac{|A|^2}{H^2}\\
{} & {} & -\bigg\langle \nabla\bigg(\frac{|A|^2}{H^2}\bigg),
\nabla \log H^2\bigg\rangle.
\end{eqnarray*}
In other words,
\begin{equation}
\label{delta-huisken}
\Delta_{-cT-\nabla\log H^2}\bigg(\frac{|A^2|}{H^2}\bigg)=
\frac{2}{H^4}|\nabla H
\otimes A-H\nabla A|^2+
2\bigg(m-\frac{H^2}{|A|^2}\bigg)\kappa\frac{|A|^2}{H^2}
\end{equation}
Note that, by Newton's inequalities
\[
m|A|^2-H^2\ge 0,
\]
with equality holding if and only if the immersion is totally umbilical
and therefore, since $\kappa\ge 0$, (\ref{delta-huisken}) yields 
\begin{equation}
\label{delta-huisken-2}
\Delta_{-cT-\nabla\log H^2}\bigg(\frac{|A^2|}{H^2}\bigg)\ge 0\,\, \mbox{ on } \,\, M.
\end{equation}
Since $M$ is compact, the strong maximum principle then implies that $|A|^2/H^2$ is a positive constant on $M$ and hence by (\ref{delta-huisken}) we also have
\begin{equation}
\label{luis24}
|\nabla H\otimes A-H\nabla A|^2\equiv 0 \quad \text{ on } M
\end{equation}
and
\begin{equation}
\label{luis25}
\bigg(m-\frac{H^2}{|A|^2}\bigg)\kappa\equiv 0 \quad \text{ on } M.
\end{equation}
If $\kappa>0$, (\ref{luis25}) directly implies $m|A|^2=H^2$ and the immersion is totally umbilical (and hence $H$ is constant). 

On the other hand, if $\kappa=0$ from (\ref{luis24}) we obtain $\nabla H\otimes A=H\nabla A$ that, by Codazzi equation, is a symmetric tensor. Therefore,
\begin{equation}
\label{luis26}
U(H)AV=V(H)AU
\end{equation}
for any $U,V\in\Gamma(TM)$. In this case we conclude that $H$ is constant. Indeed, otherwise, let $\mathcal{U}=\{ x\in M : \nabla H(x)\neq 0 \}\neq\emptyset$. At a given point $x\in\mathcal{U}$ we may choose a local orthonormal frame $\{E_1,\ldots, E_m\}$ such that $E_1=\nabla H/|\nabla H|$ and $E_i(H)=0$ for every $i\ge 2$. Equation (\ref{luis26}) implies then that $AE_1=HE_1$ and $AE_i=0$ for every $i\ge 2$. In particular, $|A|^2=H^2$ on $\mathcal{U}$ and therefore everywhere on $M$. From (\ref{deltaH-huisken}) we then have
\begin{equation}
\label{luis27}
\Delta H=-(c\varphi+H^2)H-c\langle T,\nabla H\rangle.
\end{equation}
On the other hand, using equation (\ref{luis21}) we also obtain
\begin{equation}
\label{luis28}
\mathrm{div}(HT)=\langle T,\nabla H\rangle+H\mathrm{div}T=
\langle T,\nabla H\rangle+H\bigg(m\varphi+\frac{H^2}{c}\bigg).
\end{equation}
Putting the two equations together we deduce
\begin{equation}
\label{luis29}
\Delta H+c\,\mathrm{div}(HT)=c(m-1)\varphi H.
\end{equation}
Integrating (\ref{luis29}) gives $\int_M\varphi H=0$. Since $H>0$ and $\varphi$ does not change sign on $M$, this implies $\varphi\equiv 0$. Inserting $\varphi\equiv 0$ into (\ref{luis27}) and using the strong maximum principle we get the desired contradiction. 
Summing up, $H$ is constant on $M$ and this finishes the proof. \hfill $\square$

\vspace{3mm}

We now extend the result to the complete, non-compact case, under some expected further assumptions.

\begin{theorem}
\label{huisken-complete}
Let $\bar M^{m+1}= I\times_h P^m$ be of  constant sectional curvature
$\kappa\ge 0$ and let $\psi:M^m\to \bar M^{m+1}$ be a
\emph{complete}, non-compact, orientable  codimension one mean curvature flow soliton with respect to $X = h(t)\partial_t$ and such that  $H = \langle {\bf H}, N\rangle$ does not change sign with respect to the unit normal vector field $N$. 
Assume that the function
\[
\frac{\eta }{|H|}(x)\to +\infty \,\, \mbox{ as } \,\, x\to\infty \,\, \mbox{ in }\,\, M.
\]
and  suppose 
\[
\limsup_{x\to \infty}\frac{1}{|H|}\big(m h'(\pi\circ\psi)+ch^2(\pi\circ\psi)\big)<+\infty.
\]
Furthermore assume
\[
\limsup_{x\to\infty} (ch'(\pi\circ\psi)+|A|^2)<0
\]
and
\begin{equation}
\label{AHsup}
\sup_M \frac{|A|^2}{H^2}<+\infty.
\end{equation}
\begin{itemize}
\item[(i)] If $\kappa>0$ then either $H\equiv 0$ or the immersion is totally umbilical (and hence $H>0$ is constant).
\item[(ii)] If $\kappa=0$ then 
\begin{equation}
\label{Q0}
|\nabla H
\otimes A-H\,\nabla A|^2=0
\end{equation}
\end{itemize}
\end{theorem}

\noindent \emph{Proof.}  As in Theorem \ref{huisken-cl}, since $H = \langle {\bf H}, N\rangle$ does not change sign, then either $H\equiv 0$ or $H$ is never zero on $M$ by the maximum principle applied to (\ref{deltaH-huisken}) with $\varphi=h'(\pi\circ\psi)$. If $H\equiv 0$ then there is nothing to prove; otherwise we can suppose, up to have chosen $N$ appropriately, that $H>0$ on $M$.  With $H>0$ on $M$, considering the weighted Laplacian   
$\Delta_{-c\eta-\log H^2}$ and taking into account that now $T=\nabla\eta$,
we see that we can  rewrite (\ref{delta-huisken}) into the form
\begin{equation}
\label{huisken-w}
\Delta_{-c\eta-\log H^2}\bigg(\frac{|A^2|}{H^2}\bigg)=\frac{2}{H^4}|\nabla H
\otimes A-H\nabla A|^2+
2\bigg(m-\frac{H^2}{|A|^2}\bigg)\kappa\frac{|A|^2}{H^2}.
\end{equation}
where
\[
m|A|^2-H^2\ge 0
\]
with equality holding if and only of $M$ is totally umbilical, 
and therefore, since $\kappa\ge 0$, (\ref{huisken-w}) yields 
\begin{equation}
\label{huisken-w2}
\Delta_{-c\eta-\log H^2}\bigg(\frac{|A^2|}{H^2}\bigg)\ge 0\,\, \mbox{ on } \,\, M.
\end{equation}

The next step is to show that $M$ is parabolic with respect to $\Delta_{-c\eta-\log H^2}$.
Towards this aim we compute
\begin{eqnarray*}
\Delta_{-c\eta}\bigg(\frac{1}{H}\bigg)= -\frac{1}{H^{2}}\Delta_{-c\eta} H + \frac{2}{H^{3}} |\nabla H|^2.
\end{eqnarray*}
Hence, using (\ref{deltaH-huisken})
\begin{eqnarray*}
\Delta_{-c\eta}\bigg(\frac{1}{H}\bigg) & = &  \frac{1}{H} (ch'(\pi\circ\psi)+|A|^2) + \frac{2}{H^{3}} |\nabla H|^2\\
{} & = & \frac{1}{H} (ch'(\pi\circ\psi)+|A|^2) + \frac{2}{H} |\nabla\log H|^2.
\end{eqnarray*}
By Proposition \ref{laplace-eta},
\[
\Delta_{-c\eta} \eta = m\varphi +c |X|^2=mh'+ch^2,
\]
where to simplify notations we write $h$ and $h'$ respectively for denoting $h(\pi\circ\psi)$ and $h'(\pi\circ\psi)$. Using the above equations, we obtain
\begin{eqnarray*}
& & \Delta_{-c\eta}\bigg(\frac{\eta}{H}\bigg) =\frac{1}{H}\Delta_{-c\eta} \eta + \eta \Delta_{-c\eta}\bigg(\frac{1}{H}\bigg)-\frac{2}{H^{2}} \langle \nabla H, \nabla \eta\rangle\\
& & \,\, = \frac{1}{H} (mh'+ch^2)+\frac{1}{H} (ch'+|A|^2)\eta 
+ \frac{2}{H} |\nabla\log H|^2\eta-\frac{2}{H} \langle \nabla \log H, \nabla \eta\rangle.
\end{eqnarray*}
Observing that
\begin{eqnarray*}
\bigg\langle\nabla \log H^2, \nabla\frac{\eta}{H}\bigg\rangle = \frac{2}{H} \langle \nabla \log H, \nabla \eta\rangle-\frac{2}{H}|\nabla \log H|^2\eta,
\end{eqnarray*}
we rewrite the above in the form
\begin{equation}
\Delta_{-c\eta-\log H^2}\bigg(\frac{\eta}{H}\bigg)=\frac{1}{H} (mh'+ch^2)+(ch'+|A|^2)\frac{\eta}{H}.
\end{equation}

Setting
\[
v= \frac{\eta}{H},
\]
the assumptions of the Theorem give
\begin{align*}
v(x) \to +\infty \,\, \mbox{ as }\,\, x\to \infty\,\, \mbox{ in } \,\, M  
\end{align*}
and
\begin{align*}
\Delta_{-c\eta-\log H^2} v < 0
\end{align*}
outside a compact set. Thus, by Theorem 4.12 of \cite{AMR}  the complete manifold $M$ is $\Delta_{-c\eta-\log H^2}$-parabolic.  

Once we know that $M$ is parabolic with respect to the operator $\Delta_{-c\eta-\log H^2}$, since by  (\ref{AHsup}) the function $|A|^2/H^2$ is  bounded above, from (\ref{huisken-w2}) we deduce that $|A|^2/H^2$ is a positive constant, so that from (\ref{huisken-w}) we immediately infer 
\begin{equation}
\label{luis31}
|\nabla H\otimes A-H\nabla A|^2\equiv 0 \quad \text{ on } M
\end{equation}
and
\begin{equation}
\label{luis32}
\bigg(m-\frac{H^2}{|A|^2}\bigg)\kappa\equiv 0 \quad \text{ on } M.
\end{equation}
If $\kappa>0$, (\ref{luis32}) directly implies $m|A|^2=H^2$ and the immersion is totally umbilical (and hence $H$ is constant). If $\kappa=0$ we obtain the validity of (\ref{Q0}). This completes the proof. \hfill $\square$

\vspace{3mm}


Here is another result in the same direction, whose proof is an application of Theorem \ref{29.1}. 
\begin{theorem}
\label{gap-G}
Let $\psi:M^m\to \bar M^{m+1}= I\times_h P$ be a complete, codimension one, orientable mean curvature flow soliton with respect to $X = h(t)\partial_t$  and such that  $H=\langle {\bf H}, N\rangle$ does not change sign with respect to the unit normal vector field $N$. Suppose that $\bar M$ has constant sectional curvature $\kappa\ge 0$, and assume
\begin{equation}
\label{growth-Hp} 
\bigg(\int_{\partial B_r}|A|^2\, e^{c\eta}\, {\rm d}M\bigg)^{-1} \notin L^1(+\infty).
\end{equation}
\begin{itemize}
\item[{\rm (i)}] If $\kappa>0$ then either $H\equiv 0$ or the immersion is totally umbilical (and hence $H>0$ is constant).
\item[{\rm (ii)}] If $\kappa=0$ then 
\begin{equation}
|\nabla H\otimes A-H\,\nabla A|^2=0
\end{equation}
\end{itemize}
\end{theorem}
\begin{remark}
Note that condition (\ref{growth-Hp}) is implied by
\begin{equation}
\label{H2p} 
|A|^2 \in L^1(M, e^{c\eta}\, {\rm d}M).
\end{equation}

\end{remark}
\noindent \emph{Proof.} 
Reasoning as at the beginning of the proof of Theorem \ref{huisken-complete}, either $H\equiv 0$ or $H$ never vanishes on $M$. If $H\equiv 0$ there is nothing to prove, so that without loss of generality we may assume $H>0$ on $M$. 

From (\ref{simons-soliton-spaceform}) in Proposition \ref{prop-simons-soliton} it follows
\begin{eqnarray*}
|A|\Delta_{-c\eta}|A| & = & \frac{1}{2}\Delta_{-c\eta}|A|-|\nabla|A||^2\\
\nonumber {} & = & -(ch'+|A|^2)|A|^2+(m|A|^2-H^2)\kappa+|\nabla A|^2-|\nabla |A||^2,
\end{eqnarray*}
and hence 
\begin{equation}
\label{DeltamodA}
|A|\Delta_{-c\eta}|A|+(ch'+|A|^2)|A|^2= (m|A|^2-H^2)\kappa +|\nabla A|^2-|\nabla |A||^2.
\end{equation}
By Newton and Kato's inequalities the right hand side of the above equation is non-negative provided $\kappa\ge 0$, and setting $u=|A|\ge 0$ and $a(x)=(ch'+|A|^2)$, we finally have
\begin{equation}
\label{DeltamodA-2}
u\Delta_{-c\eta} u+a(x)u^2 \ge 0.
\end{equation}

We let $v=H>0$, and rewrite (\ref{deltaH-huisken}) in the form
\[
\Delta_{-c\eta}v+a(x)v=0.
\]
Next we apply Theorem \ref{29.1} with the choices $u=|A|$, $v=H$, $f=-c\eta$, $\alpha=0$, 
$\mu=1$ and $\beta=0$ to conclude that  if 
\[
\bigg(\int_{\partial B_r}|A|^2\, e^{c\eta}\, {\rm d}M\bigg)^{-1} \notin L^1(+\infty)
\]
then there exists a positive constant $C\in\mathbb{R}$ such that $CH=|A|$ and therefore $|A|^2/H^2$ is a positive constant on $M$. The proof is then completed as in Theorem \ref{huisken-complete}. \hfill $\square$

\begin{remark}
\label{translating-martin}
In the particular case of translating solitons $\psi: M^m\to \mathbb{R}^{m+1}$ with respect to a parallel vector field $X\in \Gamma(T\mathbb{R}^{m+1})$ we conclude from Theorem {\rm{\ref{gap-G}}} that $\psi(M)$ is a grim hyperplane proceeding exactly as in cases 1 and 2 in the proof of Theorem A in {\rm\cite{martin2}}.
\end{remark}

It is worth to observe that Theorem \ref{29.1}, that we have been using in the proof of Theorem \ref{gap-G}, can in fact be extended to the next result.
\begin{theorem}
\label{29.1.ex}
Let $(M,\langle,\rangle,e^{-f}dM)$ be a complete weighted manifold with $f\in C^\infty(M)$ and $a(x)\in L^\infty_{\rm loc}(M)$. Let $u\in{\rm Lip}_{\rm loc}(M)$ satisfy the differential inequality
\begin{equation}
\label{29.2.ex}
u\Delta_fu+a(x)u^2+\alpha|\nabla u|^2\ge 0 \quad \text{weakly on } M
\end{equation}
for some $\alpha\in\mathbb{R}$. Let $v\in{\rm Lip}_{\rm loc}(M)$ be a positive solution of 
\begin{equation}
\label{29.3.ex}
\Delta_fv+\mu a(x)v\le -K\frac{|\nabla v|^2}{v} \quad \text{weakly on } M
\end{equation}
for some $\mu$ and $K>-1$, and suppose that
\begin{equation}
\label{29.4.ex}
\alpha+1\le \mu(K+1), \quad  \mu>0 \quad \mu\ge \alpha+1.
\end{equation}
If 
\begin{equation}
\label{29.6.ex}
\left(\int_{\partial B_r}v^{\frac{\beta+1}{\mu}(2-p)}u^{p(\beta+1)}e^{-f}\right)^{-1}\notin L^1(+\infty)
\end{equation}
holds for some $p>1$ and $\beta$ satisfying
\begin{equation}
\label{29.5.ex}
\alpha\le \beta\le \mu(K+1)-1,\quad \beta>-1.
\end{equation}
then there exists a constant $C\in\mathbb{R}$ such that
\begin{equation}
\label{29.7.ex}
Cv={\rm sgn} u\, |u|^\mu.
\end{equation}
Furthermore,
\begin{enumerate}
\item[(i)] If $\alpha+1<\mu(K+1)$, then $u$ is constant on $M$, and if in addition $a(x)$ does not vanish identically, then $u\equiv 0$. 
\item[(ii)] If $\alpha+1=\mu(K+1)$ and $u$ does not vanish identically, then $v$ and therefore $|u|^\mu$ satisfy (\ref{29.3.ex}) with the equality sign.
\end{enumerate}
\end{theorem}

Then, proceeding as in the proof of Theorem \ref{gap-G} but applying now Theorem \ref{29.1.ex} with the choices 
$u=|A|$, $v=H$, $f=-c\eta$, $\alpha=0$, 
$\mu=1$, $K=0$ and $\beta=0$, we obtain the validity of Theorem \ref{gap-G} replacing assumption (\ref{growth-Hp}) with the more general assumption 
\begin{equation}
\label{growth-Hp.ex} 
\bigg(\int_{\partial B_r}|H|^{2-p} |A|^p\, e^{c\eta}\, {\rm d}M\bigg)^{-1} \notin L^1(+\infty),
\end{equation}
for some $p>1$. 

As a consequence of this extended version of Theorem \ref{gap-G}, in case of Euclidean self-shrinker we obtain
following

\begin{corollary}
\label{gap-self}
Let   $\psi:M^m\to \mathbb{R}^{m+1}$ be a complete,  orientable mean curvature flow soliton such that  $H = \langle {\bf H}, N\rangle$ is not identically zero and does not change sign with respect to the unit normal vector field $N$. 
If for some $p>1$
\begin{equation}
\label{H2p-bis} 
|H|^{2-p} |A|^p \in L^1(M, e^{c\eta}\, {\rm d}M)
\end{equation}
then either $\psi(M)$ is a sphere or a product of a sphere and an Euclidean factor;  or $\psi(M)$ is  cylinder of the form $\gamma \times \mathbb{R}^{n}$ where $\gamma$ is a mean curvature flow soliton in an Euclidean plane.
\end{corollary}

\noindent \emph{Proof.} 
We proceed as in the proof of Theorem A in \cite{martin2}, cases 1 and 2. It follows from
\[
H \nabla A - \nabla H \otimes A = 0
\]
that 
\begin{equation}
\label{cod1}
H \nabla_{E_i} A(E_j, E_k) - \langle \nabla H, E_i \rangle A(E_j, E_k) =0
\end{equation}
for a given orthonormal tangent frame $\{E_i\}_{i=1}^m$. It follows from Codazzi's equation and the relation above that
\[
H \nabla_{E_i} A(E_j, E_k)  = H \nabla_{E_j} A(E_i, E_k) =\langle \nabla H, E_j\rangle  A(E_i, E_k).
\]
Therefore
\begin{equation}
\label{cod2}
\langle \nabla H, E_j\rangle  A(E_i, E_k) -\langle \nabla H, E_i \rangle A(E_j, E_k) =0.
\end{equation}
Case 1:  assume $H$ is constant. Since by assumption $H={\rm cte}\neq 0$, it follows from (\ref{cod1}) that
\[
\nabla_{E_i}A(E_j, E_k) =0
\]
for all $i, j, k$. Therefore $\nabla A =0$ on $M$ what implies $\psi(M)$ is locally isometric to a sphere or a product of a sphere and an Euclidean factor \cite{Lawson}.  Those are examples of self-shrinkers in Euclidean space.  

\vspace{3mm}
\noindent Case 2: Now we suppose that $\nabla H\neq 0$ on an open subset $U\subset M$. Setting
\[
{E}_1 = \frac{\nabla H}{|\nabla H|}
\]
on $U$, it follows from ({\rm\ref{cod2}}) that
\[
A(E_i, E_k) =0
\]
for $i\neq 1$ and any $k$. This implies that $\psi$ has only one non-zero principal curvature (since we are assuming that $H\neq 0$) necessarily in the direction of $E_1$. Then the nullity distribution 
\[
\mathcal{D}(x) = \{{\sf v}\in T_x U: A({\sf v}, \cdot)=0\}
\]
has constant rank. Therefore, it is smooth and integrable with totally geodesic integral leaves (their images by $\psi$ being totally geodesic submanifolds in $\mathbb{R}^{m+1}$), see \cite{Ferus}, \cite{martin2}. Its orthogonal complement
\[
\mathcal{D}^\perp(x) = {\rm span}\{E_1(x)\}, \quad x\in U
\]
is also smooth and integrable and its integral curves are geodesics in $M$. Moreover both distributions are parallel. One concludes  that the submanifold is a product of a curve and an Euclidean factor. Then the curve must is a mean curvature flow soliton  in an Euclidean plane. \hfill $\square$

\vspace{3mm}

We end this section with another application of our Simon's type formula for mean curvature flow solitons given by the following
\begin{theorem} 
\label{theorem44}
Let $\psi: M^m \to \bar M^{m+1} = I\times_h P$ be a complete, codimension one mean curvature flow soliton with respect to $X= h(t)\partial_t$ and suppose that $\bar M$ has constant sectional curvature $\kappa$. Assume 
\begin{equation}
\label{supAhk}
\sup_M (ch'(\pi\circ\psi)+|A|^2)<m\kappa, \quad \text{and} \quad 
\inf_M ch'(\pi\circ\psi)>-\infty.
\end{equation} 
Then the immersion is totally umbilical. 
\end{theorem}

\noindent \emph{Proof.}  Let $B$ be the umbilicity tensor, that is
\[
B=A-\frac{H}{m}I.
\]
Since  
\[
|B|^2= |A|^2 - \frac{1}{m} H^2\ge 0,
\]
and vanishes at the umbilical points of the immersion, it is enough to show that the function $u=|A|^2 -\frac{1}{m}H^2 \ge 0$ vanishes on $M$. Towards this aim we recall from (\ref{delta-eta-H2-bis}) that
\[
\frac{1}{2}\Delta_{-c\eta} H^2 = -(ch'+|A|^2)H^2 + |\nabla H|^2,
\]
where to simplify notation we have set $h'$ in place of $h'(\pi\circ\psi)$.
On the other hand, from (\ref{simons-soliton-spaceform}) we also have
\begin{equation}
\frac{1}{2}\Delta_{-c\eta} |A|^2 = -(ch'+|A|^2)|A|^2 + m\kappa \bigg(|A|^2-\frac{1}{m} H^2\bigg)+|\nabla A|^2.
\end{equation}
Furthermore, a simple computation yields
\[
|\nabla A|^2 - \frac{1}{m}|\nabla H|^2 = \Big|\nabla A - \frac{1}{m} \nabla H \otimes g\Big|^2,
\]
where $g$ is the metric in $M$ induced by the immersion $\psi$. Putting together these last three equations and recalling the definition of $u$ we obtain
\begin{equation}
\label{delta-eta-uu}
\frac{1}{2}\Delta_{-c\eta} u = (m\kappa-ch'-|A|^2) u + \Big|\nabla A - \frac{1}{m}\nabla H \otimes g\Big|^2.
\end{equation}
Under the assumption (\ref{supAhk}) one has $\sup_M (ch'+|A|^2)<+\infty$ and therefore, by Corollary \ref{cor-cond-wmp-bis}, we have the validity of the weak maximum principle for the operator $\Delta_{-c\eta}$ on $M$. Furthermore, (\ref{supAhk}) gives
\[
(m\kappa-ch'-|A|^2)\geq m\kappa-\sup_M (c\varphi+|A|^2)=\frac{\varepsilon}{2}
\]
for some $\varepsilon>0$, 
and thus from (\ref{delta-eta-uu}) we obtain
\begin{equation}
\label{idelta-eta-uu-2}
\Delta_{-c\eta} u \ge \varepsilon u \,\, \mbox{ on } \,\, M. 
\end{equation}
Assumption (\ref{supAhk}) implies also that $\sup_M u<+\infty$, since
$\sup_M u\le\sup_M|A|^2\leq \sup_M(ch'+|A|^2)-\inf_M(ch')<+\infty$. 
Hence, using (\ref{idelta-eta-uu-2}) and the weak maximum principle for $\Delta_{-c\eta}$ we deduce that $u^*=\sup_M u =0$, taht is, $u\equiv 0$ on $M$ as desired. \hfill $\square$

\vspace{3mm}

For instance as consequence of the above Theorem we have 
\begin{corollary}
\label{euc-A1}
Let $\psi: M^m \to \mathbb{R}^{m+1}$ be a complete, codimension one, mean curvature flow soliton with respect to  the position vector field $X(x) = x$ with $c<0$. Let $A$ be the Weingarten operator in the direction of some  unit normal. Assume
\[
\sup_M |A|^2 < -c.
\]
Then $\psi$ is totally umbilical. 
\end{corollary}

\begin{remark}
Corollary \ref{euc-A1} compares with Corollary \ref{gap-Fbis}.
\end{remark}

Obviously similar results hold for other space forms. We quote the following.

\begin{corollary}
\label{hyp-cosh}
Let $\psi:M^m\to \mathbb{H}^{m+1}=\mathbb{R}\times_{\cosh t} \mathbb{H}^m$ be a complete, codimension one, mean curvature flow soliton with respect to $X =\cosh t \partial_t$ and let $A$ be the Weingarten operator in the direction of some unit normal. Assume
\begin{equation}
\label{supAhk2}
\sup_M |A|^2 < -c\inf_M \sinh (\pi\circ\psi) - m, \quad \text{and} \quad
\inf_M c\sinh(\pi\circ\psi)>-\infty.
\end{equation}
Then $\psi$ is totally umbilical.
\end{corollary}

Note that condition (\ref{supAhk2}) obviously depend on the representation of $\mathbb{H}^{m+1}$ and of the corresponding choice of $X$. For instance,

\begin{corollary}
\label{hyp-exp}
Let $\psi:M^m\to \mathbb{H}^{m+1}=\mathbb{R}\times_{e^t} \mathbb{R}^m$ be a complete, codimension one, mean curvature flow soliton with respect to $X =e^t \partial_t$ and $c<0$. Let $A$ be the Weingarten operator in the direction of some unit normal. Assume
\begin{equation}
\psi(M) \subset [a,b] \times \mathbb{R}^m
\end{equation}
and
\begin{equation}
\label{supAhk3}
\sup_M |A|^2 < -ce^a-m.
\end{equation}
Then $\psi$ is totally umbilical.
\end{corollary}

\section{Translating solitons}

In this section we consider translating solitons $\psi: M^m \to \bar M^{n+1}=\mathbb R\times P^n$ mentioned earlier in Example \ref{translating-solitons}. Without loss of generality we may assume that $c<0$.

\begin{theorem}
\label{th-transl-1}
Suppose that the graph $\Gamma_u: P^n \to\mathbb R\times P^n$ is a translating soliton with $\Gamma_u (P) \subset [a,+\infty) \times P$ for some $a\in\mathbb R$. Assume that $(P, g_0)$ is complete and 
\begin{equation}
\label{vol-P}
\frac{1}{{\rm vol}_P (\partial B_r)}\notin L^1(+\infty),
\end{equation}
where $\partial B_r$ is the boundary of the geodesic ball in $(P,g_0)$ of radius $r$ centered at a fixed origin $o$. 
Then $u$ is constant and $\Gamma_u$ is a slice of the natural foliation of $I\times P$.
\end{theorem}

\noindent \emph{Proof.} Let $f = \log W = \log \sqrt{1+|\nabla^P u|^2} \ge 0$. Then
\[
{\rm vol}_{P, f} (\partial B_r) = \int_{\partial B_r} e^{-f}\, {\rm d}P \le {\rm vol}_P (\partial B_r)
\]
and assumption (\ref{vol-P}) implies 
\begin{equation}
\label{vol-P2}
\frac{1}{{\rm vol}_{P,f} (\partial B_r)} \notin L^1(+\infty).
\end{equation}
Completeness of $(P, g_0)$ and (\ref{vol-P2}) imply by Theorem 4.14 of \cite{AMR} that $P$ is parabolic with respect to the operator
\[
e^f \, {\rm div}_P (e^{-f}\nabla^P).
\]
By equation (\ref{pde-translating}) the function $u$ defining the translating soliton $\Gamma_u$ satisfies
\begin{equation}
\label{pde-transl-w}
e^f \, {\rm div}_P (e^{-f}\nabla^P u) =c<0
\end{equation}
and since $u\ge a$ on $P$ we deduce that it is constant completing the proof. \hfill $\square$

\vspace{3mm}

A further, in some sense refined, version of Theorem \ref{th-transl-1}  is given in the next result. 

\begin{theorem}
\label{th-transl-2}
Assume that $(P, g_0)$ is complete and set $r(x) = {\rm dist}_P(o, x)$ for some fixed origin $o\in P$. Suppose that, for some $0\le \sigma < 2$ 
\begin{equation}
\label{vol-transl-r}
\lim_{r\to+\infty} \frac{\log{\rm vol}_P (B_r)}{r^{2-\sigma}} =0.
\end{equation}
Then there are no translating soliton graphs $\Gamma_u : P \to \mathbb{R}\times P$ lying above the graphs
\[
\iota_\beta(x) = (-\beta\, r(x)^\sigma, x)
\]
at infinity for any $\beta\in \mathbb{R}$.
\end{theorem}

\noindent \emph{Proof.}  We reason by contradiction and we we suppose that for $u\in C^\infty(P)$, the graph $\Gamma_u$ is a translating soliton graph in $\mathbb{R}\times P$ satisfying 
\begin{equation}
\label{limur}
\liminf_{r(x)\to +\infty} \frac{u(x)}{r(x)^\sigma}\ge -\beta
\end{equation}
for some $\beta\in \mathbb{R}$. We let
\[
f =\log W = \log \sqrt{1+|\nabla^P u|^2}\ge 0
\]
and we observe that the validity of (\ref{vol-transl-r}) implies 
\begin{equation}
\label{vol-transl-rw}
\lim_{r\to+\infty} \frac{\log{\rm vol}_{P,f} (B_r)}{r^{2-\sigma}} =0
\end{equation}
Furthermore, $u$ satisfies (\ref{pde-transl-w}). From the latter and Theorem 4.4 of \cite{AMR} we obtain
\[
0\le \sup_P e^f\, {\rm div}_P (e^{-f}\nabla^P u) =c<0,
\] 
a contradiction. \hfill $\square$

\section{Variational setting and stability of solitons}
\label{variational}

For the sake of simplicity in this section we restrict ourselves to the case of codimension one mean curvature flow solitons in warped spaces $\bar M^{m+1} = I\times_h P^m$. Recalling the function  
\[
\hat\eta(t) = \int^{t}_{t_0} h(s)\, {\rm d}s, \quad t_0\in I
\]
we introduce, for a fixed $c\in \mathbb{R}$, the weighted volume functional
\begin{equation}
\label{vol-Mw}
\mathcal{A}_{c\eta} [\psi,\Omega]={\rm vol}_{c\eta} (\psi(\Omega)) = \int_\Omega e^{c\eta}\, {\rm d}M, 
\end{equation}
where $\eta(x)=\hat\eta(\pi\circ\psi(x))$, ${\rm d}M$ is the volume element induced in $M$ from a given isometric immersion $\psi: M^m\to \bar M^{m+1}$ and $\Omega$ is a relatively compact domain of $M$. 

We have the following 
\begin{proposition}
\label{first-variation}
Let $\psi: M^m \to \bar M^{m+1}= I\times_h P$ be a codimension one mean curvature flow soliton with respect to the conformal vector field $X= h(t)\partial_t$. Then the equation 
\begin{equation}
\label{sol-scalar-bis}
H = c\,\langle X, N\rangle
\end{equation}
on the relatively compact domain $\Omega \subset M$ is the Euler-Lagrange equation of the functional {\rm (\ref{vol-Mw})}. Moreover, the second variation formula for normal variations is given by
\begin{equation}
\label{second-variation-vol}
\delta^2\mathcal{A}_{c\eta}[\psi, \Omega]\cdot (u,u)=\int_M e^{c\eta}\,
uL_{c\eta}u\, {\rm d}M,\quad u\in C^\infty_0(\Omega),
\end{equation}
where the stability operator $L_{c\eta}$ is defined by
\begin{equation}
\label{11.4}
L_{c\eta} u = \Delta_{-c\eta} u + (|A|^2 +{\rm Ric}_{\bar M,-c\bar\eta}(N, N)) u.
\end{equation}
As before
\begin{equation}
\label{deltaT-op}
\Delta_{-c\eta} = \Delta +c \langle \nabla\eta, \nabla\,\cdot\,\rangle,
\end{equation}
and ${\rm Ric}_{\bar M,-c\bar\eta}$ denotes the Bakry-Emery-Ricci tensor on $\bar M$ associated to the function $-c\bar\eta$, where $\bar\eta=\hat\eta\circ\pi$. That is, 
\begin{equation}
\label{kk4}
{\rm Ric}_{\bar M,-c\bar\eta} = {\rm Ric}_{\bar M} -c\bar\nabla\bar\nabla\bar\eta = {\rm Ric}_{\bar M} - ch'(\pi)\langle \cdot, \cdot\rangle.
\end{equation}
\end{proposition}

\begin{remark}  
If $\bar M^{m+1} = \mathbb{R}^{m+1}=(0,+\infty)\times \mathbb{S}^m$ and $X = t\partial_t$ then
\[
\mathcal{A}_{c\eta} [\psi,\Omega]={\rm vol}_{c\eta} (\Omega) = \int_\Omega e^{-|X|^2/2}\, {\rm d} M
\]
coincides with the usual functional for self-shrinkers and self-expanders in Euclidean space.
\end{remark}

\noindent \emph{Proof.} Given $\varepsilon >0$, let  $\Psi:(-\varepsilon, \varepsilon)\times M \to \bar M$ be a compactly supported variation of $\psi$ in $\Omega\subset M$  with $\Psi(0,\,\cdot\,) = \psi$ and
normal variational vector field
\[
\frac{\partial\Psi}{\partial s}\Big|_{s=0} = uN
\]
for some function $u\in C^\infty_0 (\Omega)$ and $N$ a local unit normal vector field along $\psi$.  Then
\begin{eqnarray*}
& &\frac{{\rm d}}{{\rm d}s}\Big|_{s=0}{\rm vol}_{c\eta}[\Psi_s(\Omega)] =
\int_\Omega e^{c\eta}\,(c\langle X, N\rangle-H)u\, {\rm d}M.
\end{eqnarray*}
Hence, stationary immersions are characterized by the
scalar soliton equation
\[
H-c\,\langle X, N\rangle =0\,\,\, \mbox{ on } \,\,\, \Omega \subset M
\]
which yields (\ref{sol-scalar-bis}). Now we compute the second variation formula. At a stationary immersion
we have
\begin{eqnarray*}
& &\frac{{\rm d}^2}{{\rm d}s^2}\Big|_{s=0}{\rm vol}_{c\eta}[\Psi_s(\Omega)] = 
\int_M e^{c\eta}\,\frac{d}{ds}\Big|_{s=0}(c\langle X_s,
N_s\rangle-H_s)u\, {\rm d}M.
\end{eqnarray*}
Using the fact that
\[
\bar\nabla_{\partial_s} N = -\nabla u,
\]
we compute
\begin{eqnarray*}
\frac{d}{ds}\Big|_{s=0} \langle X, N\rangle = \langle
\bar\nabla_{\partial_s}X, N\rangle
+\langle X, \bar\nabla_{\partial_s} N\rangle=\varphi\langle\partial_s, N\rangle
-\langle X,\nabla u\rangle,
\end{eqnarray*}
where $\varphi=h'(\pi\circ\psi)$, and recalling that $X^\top=\nabla\eta$ we conclude that
\[
\frac{d}{ds}\Big|_{s=0} \langle X, N\rangle =u\varphi -\langle X,
\nabla u\rangle =u\varphi -\langle \nabla\eta,
\nabla u\rangle.
\]
Since
\begin{eqnarray*}
\frac{d}{ds}\Big|_{s=0} H = \Delta u + (|A|^2+{\rm Ric}_{\bar
M}(N,N))u,
\end{eqnarray*}
we obtain
\begin{eqnarray*}
& & \frac{d}{ds}\big(H-c\langle X, N\rangle\big) =\Delta u +c\langle
\nabla\eta, \nabla u\rangle +|A|^2 u + \big({\rm Ric}_{\bar
M}(N,N)-c\varphi\big)u\\
& & \,\, =\Delta_{-c\eta}u + 
|A|^2 u +
{\rm Ric}_{\bar M,-c\bar\eta}(N,N)u,
\end{eqnarray*}
where
\begin{equation}
{\rm Ric}_{\bar M,-c\bar\eta} = {\rm Ric}_{\bar M}
-c\bar\nabla\bar\nabla\bar\eta
\end{equation}
is the Bakry-Emery-Ricci tensor associated to the function $-c\bar\eta$.
This latter follows from the fact that
\begin{eqnarray*}
{\rm Ric}_{\bar M, -c\bar\eta}(N,N) & = &  {\rm Ric}_{\bar
M}(N,N)-c\langle\bar\nabla_N\bar\nabla\bar\eta, N\rangle \\
{} & = & {\rm Ric}_{\bar M}(N,N)-c\langle\bar\nabla_N X, N\rangle\\
{} & = & {\rm Ric}_{\bar M}(N,N)-c\varphi.
\end{eqnarray*}
This finishes the proof of the proposition. \hfill $\square$

\

\

We introduce the obvious
\begin{definition}
The codimension one, mean curvature flow soliton $\psi:M^m \to \bar M^{m+1} = I\times_h P$ with respect to $X=h(t)\partial_t$ is said to be stable, with finite (infinite) index if so is the stability operator $L_{c\eta}$ defined in (\ref{11.4}).
\end{definition}

For instance, using equations {\rm (\ref{kk3})} and (\ref{kk4}) we deduce that 
\begin{equation}
\Delta_{-c\eta} H=-H{\rm Ric}_{\bar M,-c\bar\eta}(N,N)H-|A|^2H
-\left(2ch'+m\frac{h''}{h}\right)H,
\end{equation}
or in other words
\begin{eqnarray}
\label{eigenvalue-L}
L_{c\eta}H & = & \Delta_{-c\eta} H +
\left(|A|^2+{\rm Ric}_{\bar M,-c\bar\eta}(N,N)\right)H \nonumber \\
{} & = &-\left(2ch'+m\frac{h''}{h}\right)H.
\end{eqnarray}


From equation (\ref{eigenvalue-L}) we have the next consequence.
\begin{proposition}
\label{prop-stable}
Let $\psi:M^m \to \bar M^{m+1} = I\times_h P$ be a complete, codimension one mean curvature flow soliton with respect to $X=h(t)\partial_t$. Let $\psi$ be orientable and suppose that we can choose a unit normal $N$ such that $H = \langle {\bf H}, N\rangle>0$ on $M$. Furthermore assume that
\[
\inf_M \bigg(2ch'+m\frac{h''}{h} \bigg)\ge 0.
\]
Then the immersion $\psi$ is stable.
\end{proposition}

For the proof of Proposition \ref{prop-stable} we first recall that from Lemma 3.10 of \cite{PRSgreenbook} with a minor variation of the proof we deduce the validity of the following 
\begin{lemma}
Let $(M,\langle,\rangle)$ be a complete Riemannian manifold, 
$f\in C^\infty(M)$, $\Omega\subset M$ be a domain in $M$ and let $q(x)\in L^\infty_{loc}(\Omega)$. The following facts are equivalent:
\begin{itemize}
\item[(i)] There exists $w\in C^1(\Omega)$, $w>0$, weak solution of
\[
\Delta_fw-q(x)w=0 \quad \text{on} \quad \Omega.
\]
\item[(ii)] There exists $\varphi\in H^1_{loc}(\Omega)$, $\varphi>0$, weak solution of
\[
\Delta_f\varphi-q(x)\varphi\le 0 \quad \text{on} \quad \Omega.
\]
\item[(iii)] $\lambda_1^{L_f}(\Omega)\ge 0$,
where $L_fu=\Delta_fu-q(x)u$.
\end{itemize}
\end{lemma}
Here $\lambda_1^{L_f}(\Omega)$ is the first eigenvalue of the operator $L_f$ on $\Omega$ which is variationally characterized by the Rayleigh quotient
\[
\lambda_1^{L_f}(\Omega)=\inf_{\varphi\in C^\infty_c(\Omega), \varphi\not\equiv 0}\, 
\frac{\int_\Omega\left(|\nabla\varphi|^2+q(x)\varphi^2\right)e^{-f}}{\int_\Omega\varphi^2e^{-f}}.
\]
Now from equation (\ref{eigenvalue-L}) and our assumptions we have
\[
L_{c\eta} H \le 0\,\, \mbox{ on } \,\, M. 
\]
Then, by the Lemma above we immediately conclude that the operator $L_{c\eta}$ is stable.

\

\

In the proof of our next results we will apply the following weighted version of Theorem 7.8 in \cite{BMR}, in the simplest case where (following the notation of Theorem 7.8 in \cite{BMR}) one has $\alpha=1$ and $\beta=0$. The proof is just a simple adaptation of the proof in \cite{BMR} for the standard Laplace-Beltrami operator.
\begin{lemma}
\label{kk5}
Let $M$ be a complete, non-compact Riemannian manifold such that
\begin{equation}
\label{vol-cond-0}
{\rm vol}_{c\eta} (\partial B_r) \le C e^{\alpha r} \,\, \mbox{ for } \,\, r\gg 1
\end{equation}
and some constants $C, \alpha\ge 0$, where $\partial B_r$ is the boundary of a geodesic ball in $M$ of radius $r$ centered at a fixed origin $o\in M$. Then
\begin{equation}
\label{kk6}
\lambda_1^{\Delta_{-c\eta}} (M\backslash B_R) \le \frac{\alpha^2}{4}.
\end{equation}
\end{lemma}

Now we are ready to give our first result.
\begin{theorem}
\label{thm-G}
Let $\psi:M^m \to \bar M^{m+1} = I\times_h P$ be a complete, codimension one mean curvature flow soliton with respect to $X=h(t)\partial_t$. Suppose that $\bar M$ is Einstein and that 
\begin{equation}
\sup_M \bigg(m\frac{h''}{h}+ch'\bigg)(\pi\circ\psi)=\Theta\le 0.
\end{equation}
Assume the volume growth condition (\ref{vol-cond-0}) and that $\psi$ has finite index. Then 
\begin{equation}
\label{al}
\alpha \ge 2\sqrt{-\Theta}.
\end{equation}
\end{theorem}

\noindent \emph{Proof.} First we observe that since $\bar M$ is Einstein, we necessarily have
\[
{\rm Ric}_{\bar M} = -m\frac{h''}{h}\, \langle\cdot, \cdot\rangle
\]
and the stability operator $L_{c\eta}$ becomes
\[
L_{c\eta} = \Delta_{-c\eta} + \bigg(|A|^2 -m\frac{h''}{h}-c h'\bigg).
\]
Therefore, since $\psi$ has finite index, there exists a compact set $K\subset M$ and, by completeness, a geodesic ball $B_R \supset K$, and a solution $u\in C^2(M\backslash B_R)$, $u>0$ of $L_{c\eta}u=0$ on $M\backslash B_R$. That is,
\[
\Delta_{-c\eta} u=-\bigg(|A|^2-m\frac{h''}{h}-c h'\bigg)u \quad \mbox{ on }\quad M\backslash B_R.
\]
Thus, by a variation of Barta's theorem we obtain
\[
\lambda_1^{\Delta_{-c\eta}} (M\backslash B_R) \ge \inf_{M\backslash B_R} \Big(-\frac{\Delta_{-c\eta} u}{u}\Big)= \inf_{M\backslash B_R}\bigg(|A|^2-m\frac{h''}{h}-ch'\bigg)\ge -\Theta\ge 0.
\]

On the other hand, by the volume growth assumption (\ref{vol-cond-0}) and from Lemma \ref{kk5} we have
\[
\lambda_1^{\Delta_{-c\eta}} (M\backslash B_R) \le \frac{\alpha^2}{4}\cdot
\]
Putting the two inequalities together we have
\[
\alpha^2 +4 \Theta\ge 0,
\]
that is, (\ref{al}). \hfill $\square$

\

\begin{remark}
A similar result holds by replacing the assumption that $\bar M$ is Einstein with that of non-negative Ricci curvature. In that case, under the assumption
\[
\sup_M \, ch'(\pi\circ\psi)\le 0
\]
the same reasoning as above gives (\ref{al}) where now $\Theta=\sup_M \, ch'(\pi\circ\psi)$.
\end{remark}

Similarly we have
\begin{corollary}
\label{cor-G1}
Let $\psi: M^m \to \bar M^{m+1} = I\times_h P$ be a complete, codimension one mean curvature flow soliton with respect to $X=h(t)\partial_t$ and oriented by the unit normal $N$. Suppose that $\psi(M) \subseteq [a,\sup I) \times P$. Assume $ch''\le 0$,
\[
{\rm Ric}_{\bar M}(N,N) \ge 0\,\, \mbox{ on } \,\, [a, \sup I) \times P,
\]
and the volume growth request {\rm (\ref{vol-cond-0})}. Then either $\psi$ has infinite index or
\[
\alpha^2 +4ch'(a)\ge 0,
\]
where $\alpha$ is as in (\ref{vol-cond-0}). 
\end{corollary}

Specializing the corollary we obtain 
\begin{corollary}
\label{cor-G2} 
Let $\psi: M^m \to \mathbb{H}^{m+1} =\mathbb{R}\times_{e^t} \mathbb{R}^m$ be a complete, codimension one  mean curvature flow soliton with respect to $X = e^t\partial_t$ and satisfying the volume growth assumption {\rm (\ref{vol-cond-0})}. Assume that  $\psi(M) \subseteq [a,\infty) \times \mathbb{R}^m$. Then either $\psi$ has infinite index or, otherwise, 
\begin{equation}
\label{al2}
\alpha^2+4c e^{a}\ge 0,
\end{equation}
where $\alpha$ is as in (\ref{vol-cond-0}). 
\end{corollary}

\begin{remark}
Inequality {\rm (\ref{al2})} is quite interesting since, for a soliton of finite index, it relates the ``size'' of its image with the growth of the volume of geodesic spheres.
\end{remark}

 In this direction, recalling that a complete, proper self-shrinker $\psi:M^m \to \mathbb{R}^{m+1}$ has at most polynomial volume growth, we deduce

\begin{corollary}
\label{cor-G3}
There are no complete proper self-shrinkers $\psi: M^m\to\mathbb{R}^{m+1}$ with $c<0$ and image contained in a halfspace $[a,\infty)\times \mathbb{R}^m$ for any $a\in \mathbb{R}$. 
\end{corollary}

Clearly this result compares with the Halfspace Theorem by Hoffman and Meeks \cite{HMe}.

We now recall the following general fact, see Theorem 13 of \cite{BMR}

\begin{proposition}
\label{prop-vol-cond}
Let {\rm (}$M, \langle,\rangle, e^{-f}{\rm d}M{\rm )}$ be a complete weighted manifold, $a(x) \in C^0(M)$ and for
\[
v_f(r) ={\rm vol}_f (\partial B_r) = \int_{\partial B_r} e^{-f}\
\]
let 
\begin{equation}
\label{int-mean}
\bar a(r) = \frac{1}{v_f(r)}\int_{\partial B_r} a e^{-f}
\end{equation}
be the weighted spherical mean of $a(x)$. Assume 
\begin{equation}
\label{int-vol}
\frac{1}{v_f(r)}\in L^1(+\infty)
\end{equation}
and for some $r_0>0$, 
\[
A(r) \le \bar a(r) \quad \mbox{ on } \quad [r_0, \infty)
\]
with $A\in C^0(\mathbb{R}^+_0)$ and satisfying 
\begin{equation}
\label{int-div}
\lim_{r\to+\infty} \int^r_{r_0} A(s) v_f(s)\int^{+\infty}_s \frac{dt}{v_f(t)} \, ds=+\infty.
\end{equation}
Then the operator $\Delta_f + a(x)$ has infinite index. 
\end{proposition}

\begin{remark} In Proposition \ref{prop-vol-cond} conditions {\rm (\ref{int-vol})} and {\rm (\ref{int-div})} can be replaced respectively by
\[
\frac{1}{v_f(r)}\notin L^1(+\infty)
\]
and
\[
\lim_{r\to+\infty} \int^r_{r_0} A(s) v_f(s)\bigg(\int^{s}_{r_0} \frac{dt}{v_f(t)}\bigg)^\sigma =+\infty
\]
for some $r_0>0$ and $\sigma\in (0,1)$. Note that the range of $\sigma$ cannot be improved to $(0,1]$. 
\end{remark}

With the aid of the previous proposition we prove 
\begin{theorem}
\label{thm-H}
Let $\psi:M^m\to \bar M^{m+1}= I\times_h P$ be a complete, codimension one mean curvature flow soliton with respect to $X=h(t)\partial_t$ with $c<0$, local unit normal $N$ and such that $\psi(M)\subset [a,b]\times P, [a,b]\subset I$. Suppose  $h'(t)>0$ on $[a,b]$ and 
\begin{equation}
\label{ricci-pos}
{\rm Ric}_{\bar M}(N, N)\ge 0 \,\, \mbox { on } \,\, [a,b]\times P.
\end{equation}
Furthermore assume
\begin{align}
{\rm (i)} &\,\,  \label{vol-cond-1} \frac{1}{{\rm vol} (\partial B_r)} \in L^1(+\infty)  \\
{\rm (ii)} &\,\, \label{vol-cond-2} {\rm vol} (\partial B_r)\int_r^{+\infty} \frac{d\varrho}{{\rm vol} (\partial B_\varrho)}\notin L^1(+\infty)
\end{align}
Then $\psi(M)$ has infinite index. 
\end{theorem}

\noindent \emph{Proof.}  Since $\psi(M) \subset [a,b]\times P$, considered the stability operator $L_{c\eta}$ in (\ref{11.4}) we see that condition (\ref{vol-cond-1}) is equivalent to (\ref{int-vol}).  From the assumption on $h'$ we have
\[
\inf_M h'(\pi\circ\psi) =C>0.
\]
Hence, using (\ref{ricci-pos}) for the coefficient of the linear term of $L_{c\eta}$ we have
\[
|A|^2 + {\rm Ric}_{\bar M,-c\bar{eta}} (N, N)
|A|^2 + {\rm Ric}_{\bar M} (N, N) -c h' (\pi\circ\psi)\ge C
\]
on $M$. Thus, with the choice $A(r) = C$, assumption  (\ref{vol-cond-2}) is equivalent to (\ref{int-div}). Applying Proposition  \ref{prop-vol-cond} we deduce that the soliton has infinite index. \hfill $\square$

\vspace{3mm}

Specializing the theorem we obtain

\begin{corollary}
\label{cor-H1} Let $\psi: M^m \to  \mathbb{H}^{m+1}=\mathbb{R}\times_{e^t}\mathbb{R}^m$ be a complete, codimension one mean curvature flow soliton with respect to $X=e^t\partial_t$ woth $c<0$ and with $\psi(M) \subset [a,b]\times \mathbb{R}^m$.   Assume {\rm (\ref{vol-cond-1}) and (\ref{vol-cond-2})}. Then the soliton has infinite index. 
\end{corollary}

\end{document}